\documentclass[1 leqno,11pt]{amsart}
\usepackage{amssymb, amsmath,amsmath,latexsym,amssymb,amsfonts,amsbsy, amsthm,mathtools,graphicx,CJKutf8,CJKnumb,CJKulem,color}

\numberwithin{equation}{section}

\allowdisplaybreaks

\let\al=\alpha
\let\b=\beta
\let\g=\gamma
\let\d=\delta

\let\f=\frac

\let\om=\omega

\let\Om=\Omega

\let\na=\nabla
\let\th=\theta
\let\pa=\partial
\let\ep=\epsilon

\def\Int{\mathrm{Int}}
\def\rmA{{\mathrm{A}}}
\def\rmB{\mathrm{B}}
\def\rmC{\mathrm{C}}
\def\rmD{\mathrm{D}}
\def\rmE{\mathrm{E}}

\def\cF{{\mathcal F}}
\def\cG{{\mathcal G}}

\def\cL{{\mathcal L}}

\def\cR{{\mathcal R}}

\def\tc{\widetilde{c}}

\def\R{\mathbb R}

\def\N{\mathbb N}
\def\bbT{\mathbb{T}}

\def\eqdef{\buildrel\hbox{\footnotesize def}\over =}

\newcommand{\beq}{\begin{equation}}
\newcommand{\eeq}{\end{equation}}
\newcommand{\ben}{\begin{eqnarray}}
\newcommand{\een}{\end{eqnarray}}
\newcommand{\beno}{\begin{eqnarray*}}
\newcommand{\eeno}{\end{eqnarray*}}

\newtheorem{theorem}{Theorem}[section]
\newtheorem{definition}[theorem]{Definition}
\newtheorem{lemma}[theorem]{Lemma}
\newtheorem{proposition}[theorem]{Proposition}

\newtheorem{remark}[theorem]{Remark}

\begin{document}

\title[Linear inviscid damping for shear flows]{Linear inviscid damping and vorticity depletion for shear flows}

\author{Dongyi Wei}
\address{School of Mathematical Science, Peking University, 100871, Beijing, P. R. China}
\email{jnwdyi@163.com}

\author{Zhifei Zhang}
\address{School of Mathematical Science, Peking University, 100871, Beijing, P. R. China}
\email{zfzhang@math.pku.edu.cn}

\author{Weiren Zhao}
\address{School of Mathematical Science, Peking University, 100871, Beijing, P. R. China}
\email{zjzjzwr@126.com, zjzwr@math.pku.edu.cn}

\date{\today}

\maketitle

\begin{abstract}
In this paper, we prove the linear damping for the 2-D Euler equations around a class of shear flows under the assumption that
the linearized operator has no embedding eigenvalues. For the symmetric flows, we obtain the explicit decay estimates of the velocity,
which is the same as one for monotone shear flows. We confirm a new dynamical phenomena found by Bouchet and Morita: the depletion of the vorticity at the stationary streamlines, which could be viewed as a new mechanism leading to the damping
for the base flows with stationary streamlines.
\end{abstract}

\section{Introduction}

We consider the 2-D incompressible Euler equations in a finite channel $\Om=\big\{(x,y): x\in \bbT, y\in [-1,1]\big\}$:
\beq
\label{eq:Euler}
\left\{
\begin{array}{l}
\pa_tV+V\cdot\nabla V+\nabla P=0,\\
\na\cdot V=0,\\
V^2(t,x,-1)=V^2(t,x,1)=0,\\
V|_{t=0}=V_0(x,y).
\end{array}\right.
\eeq
where $V=(V^1,V^2)$ and $P$ denote the velocity and the pressure of the fluid respectively.  Let $\omega=\pa_xV^2-\pa_yV^1$ be the vorticity, which satisfies
\beq\label{equ:vorticity}
\omega_t+V\cdot\nabla \omega=0.
\eeq

The appearance of large coherent structures is an important phenomena in 2-D flows. The stability of coherent structures has been an active field of fluid mechanics \cite{DR, SH}, which started in the nineteenth century with Rayleigh, Kelvin, Orr, Sommerfeld and many others. Let us mention some of classical results: Rayleigh's inflection theorem \cite{Ray} giving a necessary condition of spectral instability, Howard's semicircle theorem \cite{How}, and  Arnold's criterion for the Lyapunov stability \cite{Arnold}.  We refer to \cite{BGS, Lin, Lin-CMP, SF, LZ-CPAM, KS} and references therein for some recent mathematical studies.

In this paper, we are concerned with the asymptotic stability of the 2-D Euler equations around the shear flow $(u(y),0)$, which is a steady solution of  \eqref{eq:Euler}.  The key step toward this problem is to study the linearized Euler equations around $(u(y),0)$:
\ben\label{eq:Euler-L}
\left\{
\begin{array}{l}
\pa_t\omega+\mathcal{L}\omega=0,\\
\om|_{t=0}=\om_0(x,y),
\end{array}\right.
\een
where $\mathcal{L}=u(y)\pa_x+u''(y)\pa_x(-\Delta)^{-1}$.

For the Couette flow(i.e., u(y)=y), \eqref{eq:Euler-L} is reduced to a passive transport equation
\beno
\pa_t\om+y\pa_x\om=0,\quad \om|_{t=0}=\om_0(x,y).
\eeno
In this case, Orr \cite{Orr} first observed that the velocity will tend to 0 as $t\to \infty$, which is so-called inviscid damping. Recently, Lin and Zeng \cite{LZ} proved that if $\int_{\bbT}w_0(x,y)dx=0$, then
\beno
\|V(t)\|_{L^2_{x,y}}\le \|\om(t)\|_{L^2_xH^{-1}_y}\le \f C{\langle t\rangle}\|\om_0\|_{H^{-1}_xH^1_y}.
\eeno
The mechanism leading to the damping is the vorticity mixing driven by the shear flow. This phenomena is similar to the well-known Landau damping found by Landau in 1946 \cite{Lan}.

Due to possible nonlinear transient growth, it is a challenging task from linear damping to nonlinear damping.
Moreover, nonlinear damping is sensitive to the topology of the perturbation. Indeed, Lin and Zeng \cite{LZ} proved that nonlinear inviscid damping is not true for the perturbation in $H^s$ for $s<\f32$. Motivated by the breakthrough work of Mouhot and Villani on Landau damping \cite{MV}, Beddrossian and Masmoudi \cite{BM1} proved nonlinear inviscid damping and asymptotic stability of the 2-D Euler equations around the Couette flow in Gevrey class in the domain $\Om=\bbT\times \R$. See \cite{BM2, BM3, BM4} and references therein for more relevant works. 

For general shear flow, the linear inviscid damping is also a difficult problem due to the presence of nonlocal operator $u''(y)\pa_x(-\Delta)^{-1}$. In this case, the linear dynamics is associated with the singularities at the critical layers $u=c$ of the solution for the Rayleigh equation
\beno
(u-c)(\phi''-\al^2\phi)-u''\phi=f.
\eeno
In fact, $\text{Ran}\,u$ is just the continuous spectrum of $\mathcal{L}$, whose properties and the non-normality of
$\mathcal{L}$ are related to many important
phenomena such as transient growth \cite{Far}, inviscid damping \cite{Case} and algebraic instabilities \cite{NM}.

Based on the Laplace transform and analyzing the singularity of the solution $\phi$ at the critical layer,
Case \cite{Case} gave a first prediction of linear damping for monotone shear flow.
However, there are few rigorous mathematical results. Rosencrans and Sattinger \cite{Ros} gave $t^{-1}$ decay of the stream function for analytic monotone shear flow.
Stepin \cite{Ste} proved $t^{-\nu}(\nu<\mu_0)$ decay of the stream function for monotone shear flow $u(y)\in C^{2+\mu_0}(\mu_0>\f12)$ without inflection point.
Zillinger \cite{Zill} proved $t^{-1}$ decay of $\|V(t)\|_{L^2}$ for a class of monotone shear flow, which satisfies $L\|u''\|_{W^{3,\infty}}\ll 1$. In a recent work
\cite{WZZ1}, we removed the smallness assumption in \cite{Zill} and showed that
if $u(y)\in C^4$ is monotone and $\mathcal{L}$ has no embedding eigenvalues, then
\beno
\|V(t)\|_{L^2}\leq \frac{C}{\langle t\rangle}\|\omega_0\|_{H^{-1}_xH_y^1},\quad\|V^2(t)\|_{L^2}\leq \frac{C}{\langle t\rangle^2}\|\omega_0\|_{H^{-1}_xH_y^2},
\eeno
for the initial vorticity $\om_0$ satisfying $\int_{\mathbb{T}}\om_0(x,y)dx=0$ and $P_{\mathcal{L}}\om_0=0$,
where $P_{\mathcal{L}}$ is the spectral projection to $\sigma_d\big(\mathcal{L}\big)$.

However, many important base flows such as Poiseuille flow $u(y)=y^2$ and Kolmogorov flow $u(y)=\cos y$ are not monotone.
Bouchet and Morita \cite{BM} conducted the systematic studies for the asymptotic behaviour of the vorticity and the velocity around the base flows with stationary streamlines. Based on Laplace tools and numerical computations, they found a new dynamic phenomena: depletion phenomena of the vorticity at the stationary streamlines. More precisely, they formally proved that for large times,
\beno
\widehat\om(t,\al,y)\sim \om_\infty(y)\exp(-i\al u(y)t)+O(t^{-\gamma}),
\eeno
where $\om_\infty(y_c)=0$ at stationary points $y_c$ of $u(y)$.
Based on this and using stationary phase expansion, they also predicted similar decay rates of the velocity as in the monotonic case.

The goal of this paper is to study the linear damping for general shear flows. In particular, we confirm Bouchet and Morita's prediction about the linear damping and the depletion phenomena of the vorticity for the base flows with stationary streamlines.

 The first result is the linear damping and vorticity depletion for a class of  shear flows denoted by $\mathcal{K}$, which consists of the function $u(y)$ satisfying $u(y)\in H^3(-1,1)$, and  $u''(y)\neq 0$ for critical points(i.e., $u'(y)=0$) and $u'(\pm 1)\neq 0$.

\begin{theorem}\label{thm:general}
Assume that $u(y)\in \mathcal{K}$ and the linearized operator $\mathcal{R}_\al$ defined by \eqref{def:ray ope} has no embedding eigenvalues. Assume that $\widehat{\om}_0(\al,y)\in H^1_y(-1,1)$ and $P_{\mathcal{R}_\al}\widehat{\psi}_0(\al,y)=0$, where $\psi_0$ is the stream function and  $P_{\mathcal{R}_\al}$
is the spectral projection to $\sigma_d(\mathcal{R}_\al)$.  Then it holds that
\beno
\|\widehat{V}(\cdot,\al,\cdot)\|_{L_t^2L_y^2}+\|\pa_t\widehat{V}(\cdot,\al,\cdot)\|_{L_t^2L_y^2}\leq C_{\al}\|\widehat{\om}_0(\al,\cdot)\|_{H^1_y}.
\eeno
In particular, $\displaystyle\lim_{t\to+\infty}\|\widehat{V}(t,\al,\cdot)\|_{L^2_y}=0$.
\end{theorem}

Formally, Theorem \ref{thm:general} implies that
the velocity should have at least $t^{-\f12}$ decay. To obtain the explicit decay rate of the velocity,
we consider a class of symmetric shear flow:
\beno
\text{(S)}\quad u(y)=u(-y),\quad u'(y)>0 \,\,\text{for}\,\, y> 0,\quad u'(0)=0\,\, \text{and}\,\, u''(0)>0.
\eeno
An important example is the Poiseuille flow $u(y)=y^2$.\smallskip

\begin{theorem}\label{thm:main} Assume that  $u(y)\in C^4([-1,1])$  satisfies $(S)$ and the linearized operator $\mathcal{L}$
has no embedding eigenvalues. Assume that $\int_{\mathbb{T}}\om_0(x,y)dx=0$ and $P_{\mathcal{L}}\om_0=0$, where $P_{\mathcal{L}}$ is the spectral projection to $\sigma_d\big(\mathcal{L}\big)$.
Then it holds that
\begin{itemize}
 \item[1.] if $\om_0(x,y)\in H^{-\f12}_xH^1_y$, then
\beno
\|V(t)\|_{L^2}\leq \frac{C}{\langle t\rangle}\|\omega_0\|_{H^{-\f12}_xH_y^1};
\eeno
\item[2.] if $\om_0(x,y)\in H^{\f12}_xH^2_y$, then
\beno
\|V^2(t)\|_{L^2}\leq \frac{C}{\langle t\rangle^2}\|\omega_0\|_{H^{\f12}_xH_y^2};
\eeno
\item[3.] if $\om_0(x,y)\in H^{-\f12+k}_xH^k_y$ for $k=0,1$, there exists $\om_\infty(x,y)\in H^{-\f12+k}_xH^k_y$ such that
\beno
\|\om(t,x+tu(y),y)-\om_\infty\|_{H^{-\f12+k}_xL^2_y}\longrightarrow 0\quad \textrm{as}\quad t\rightarrow +\infty.
\eeno
\end{itemize}
\end{theorem}

Let us give some remarks on our results:

\begin{itemize}

\item[1.] If $u(y)$ has no inflection points, then $\mathcal{L}$ has no embedding eigenvalues(see section 5).

\item[2.]  For non-monotone flow, the decay rate of the velocity is unexpected. Formal prediction by Case \cite{Case} does not work in this case. If we neglect the nonlocal part $u''(y)\pa_x(-\Delta)^{-1}$ of $\mathcal{L}$, one can only obtain $t^{-\f12}$ decay of the velocity, which will be explained in section 2. On the other hand, whether one can obtain the decay rate of the velocity as in Theorem \ref{thm:main} is a very interesting question for the shear flow in $\mathcal{K}$.

\item[3.] In a recent work \cite{LZ-R}, Lin and Zeng also proved the linear damping for general stable steady flows under a similar spectral
assumption in the following weak sense:
\beno
\lim_{T\to +\infty}\f 1 T\int_0^T\|V(t)\|_{L^2}^2dt=0.
\eeno

\item[4.] Our method could be applicable to the Kolmogorov flow. In a separate work \cite{WZZ2}, we will prove the depletion phenomena of the vorticity and the decay rate of the velocity conjectured by Bouchet and Morita \cite{BM} for the Kolmogorov flow.

\item[5.] Nonlinear inviscid damping should be a challenging problem, even for the Couette flow in a finite channel.

\end{itemize}

\section{New dynamic phenomena: vorticity depletion}
If we neglect the nonlocal term in the linearized vorticity equation \eqref{eq:Euler-L},  the vorticity will satisfy a passive transport equation
\beno
\pa_t\om+u(y)\pa_x\om=0,\quad \om(0,x,y)=\om_0(x,y).
\eeno
If $u'(y)>0$ and $\int_\bbT\om_0dx=0$, then we have 
\ben\label{eq:om-decay2}
\|\om(t)\|_{L^2_xH^{-1}_y}\leq \f{C}{\langle t\rangle}\|\om_0\|_{H^{-1}_xH^1_y}.
\een
If $u$ satisfies $(S)$ and $\int_\bbT\om_0dx=0$, then we have
\ben\label{eq:om-decay}
\|\om(t)\|_{L^2_xH^{-1}_y}\leq \f{C}{\langle t\rangle^\f12}\|\om_0\|_{H^{-\f12}_xH^1_y}.
\een

The damping in \eqref{eq:om-decay2} and \eqref{eq:om-decay}
is due to vorticity mixing. So, this can not explain the enhanced damping in Theorem \ref{thm:main} when $u\in (S)$.
In such case, the mechanism leading to the enhanced damping is a new dynamical phenomena found by Bouchet and Morita \cite{BM}: the depletion of the vorticity at the stationary streamlines, which is due to the effect of the perturbation velocity on the background vorticity gradient. More precisely, based on the Laplace transform analysis, they formally proved that
the solution of \eqref{eq:Euler-L} behaves as $t\to \infty$
\ben\label{eq:w-ansatz}
\widehat\om(t,\al, y)\sim \om_\infty(y)e^{-i\al u(y)t}+O(t^{-\gamma})
\een
for $\gamma>0$, where the profile $\om_\infty(y)$ vanishes at the stationary points. \smallskip

In this paper, we confirm the vorticity depletion phenomena.

\begin{theorem}\label{thm:vorticity}
Under the same assumptions as in Theorem \ref{thm:general}, if $u'(y_0)=0$, then
\beno
\lim_{t\to+\infty}\widehat{\om}(t,\al,y_0)=0.
\eeno
\end{theorem}

Assuming the ansatz \eqref{eq:w-ansatz} and $w_\infty(0)=0$,
one can obtain the decay estimate as in  \eqref{eq:om-decay2} by following the following proof of \eqref{eq:om-decay}.\smallskip

Let us prove \eqref{eq:om-decay}.
Taking Fourier transform in $x$ variable, we obtain
\beno
\pa_t\widehat{\om}(t,\al,y)+i\al u(y)\widehat{\om}(t,\al,y)=0.
\eeno
Thus, $\widehat{\om}(t,\al,y)=\widehat{\om}_0(\al,y)e^{-i\al u(y)t}$. Let $\al>0$ and $v$ be given by \eqref{def:v}. For any smooth function $\eta(y)$,
we have
\begin{align*}
\Big|\int_{-1}^1\widehat{\om}(t,\al,y)\eta(y)dy\Big|
&=\Big|e^{-i\al u(0)t}\int_{-1}^1\eta(y)\widehat{\om}_0(\al,y)e^{-i\al v(y)^2t}dy\Big|\\
&=\Big|\int_{-v(1)}^{v(1)}\eta(v^{-1}(z))\widehat{\om}_0(\al,v^{-1}(z))e^{-i\al z^2t}(v^{-1}(z))'dz\Big|\\
&=\Big|\int_{|z|^2\leq \f{\pi}{2\al t}}\eta(v^{-1}(z))\widehat{\om}_0(\al,v^{-1}(z))e^{-i\al z^2t}(v^{-1}(z))'dz\Big|\\
&\quad+\Big|\int_{v(1)^2\geq |z|^2>\f{\pi}{2\al t}}\eta(v^{-1}(z))\widehat{\om}_0(\al,v^{-1}(z))e^{-i\al z^2t}(v^{-1}(z))'dz\Big|.
\end{align*}
For the first term, we have
\begin{align*}
&\Big|\int_{|z|^2\leq \f{\pi}{2\al t}}\eta(v^{-1}(z))\widehat{\om}_0(\al,v^{-1}(z))e^{-i\al z^2t}(v^{-1}(z))'dz\Big|\\
&\leq C\|\eta\|_{L^{\infty}}\|\widehat{\om}_0(\al,\cdot)\|_{L^{\infty}_y}\Big(\int_{|z|^2\leq \f{\pi}{2\al t}}\cos (z^2\al t)dz
+\int_{|z|^2\leq \f{\pi}{2\al t}}\sin (z^2\al t)dz\Big)\\
&\leq C(\al t)^{-\f12}\|\eta\|_{L^{\infty}}\|\widehat{\om}_0(\al,\cdot)\|_{L^{\infty}_y}\Big(\int_{-\f{\pi}{2}}^{\f{\pi}{2}}\f{\cos y}{\sqrt{y}}dy
+\int_{-\f{\pi}{2}}^{\f{\pi}{2}}\f{|\sin y|}{\sqrt{y}}dy\Big)\\
&\leq C(\al t)^{-\f12}\|\eta\|_{H^1}\|\widehat{\om}_0(\al,\cdot)\|_{H^1_y}.
\end{align*}
For the second term, by integration by parts and Lemma \ref{lem:u}, we get
\begin{align*}
&\Big|\f{2}{\al t}\int_{v(1)^2\geq |z|^2>\f{\pi}{2\al t}}\eta(v^{-1}(z))\widehat{\om}_0(\al,v^{-1}(z))\f{(v^{-1}(z))'}{z}de^{-i\al z^2t}\Big|\\
&\leq \f{C}{\sqrt{\al t}}\|\eta\|_{L^{\infty}}\|\widehat{\om}_0(\al,\cdot)\|_{L^{\infty}_y}
+\f{2}{\al t}\int_{v(1)^2\geq |z|^2>\f{\pi}{2\al t}}\Big|\pa_z\Big(\eta(v^{-1}(z))\widehat{\om}_0(\al,v^{-1}(z))z^{-1}(v^{-1}(z))'\Big)\Big|dz\\
&\leq \f{C}{\sqrt{\al t}}\Big(\|\eta\|_{{L^{\infty}}}\|\widehat{\om}_0(\al,\cdot)\|_{L^{\infty}_y}
+\|\eta\|_{H^1}\|\widehat{\om}_0(\al\cdot)\|_{L^{2}_y}
+\|\eta\|_{L^2}\|\widehat{\om}_0(\al,\cdot)\|_{H^1_y}\Big).
\end{align*}
This shows that
\beno
\Big|\int_{-1}^1\widehat{\om}(t,\al,y)\eta(y)dy\Big|\le  C(\al t)^{-\f12}\|\eta\|_{H^1}\|\widehat{\om}_0(\al,\cdot)\|_{H^1_y},
\eeno
which implies \eqref{eq:om-decay}.\smallskip

\section{Sketch and key ideas of the proof}

In terms of the stream function $\psi$, the linearized Euler equations \eqref{eq:Euler-L} take
\beq\label{eq:Euler-L-stream}
\partial_t\Delta{\psi}+u(y)\partial_x\Delta{\psi}-u''(y)\partial_x{\psi}=0.
\eeq
Taking the Fourier transform in $x$, we get
\beno
(\partial_y^2-\al^2)\partial_t\widehat{\psi}=i\al\big(u''(y)-u(y)(\partial_y^2-\al^2)\big)\widehat{\psi}.
\eeno
Inverting the operator $(\partial_y^2-\al^2)$, we find
\ben\label{eq:Euler-L-operator}
-\frac{1}{i\al}\partial_t\widehat{\psi}=\mathcal{R}_\al\widehat{\psi},
\een
where
\ben\label{def:ray ope}
\mathcal{R}_\al\widehat{\psi}=-(\partial_y^2-\al^2)^{-1}\big(u''(y)-u(\partial_y^2-\al^2)\big)\widehat{\psi}.
\een
Let $\Omega$ be a simple connected domain including the spectrum $\sigma(\mathcal{R}_\al)$ of $\mathcal{R}_\al$.
We have the following representation formula of the solution to (\ref{eq:Euler-L-operator}):
\ben\label{eq:stream formula}
\widehat{\psi}(t,\al,y)=\frac{1}{2\pi i}\int_{\partial\Omega}
e^{-i\al tc}(c-\mathcal{R}_\al)^{-1}\widehat{\psi}(0,\al,y)dc.
\een
Then the large time behaviour of the solution $\widehat{\psi}(t,\al,y)$
is reduced to the study of the resolvent $(c-\mathcal{R}_\al)^{-1}$.

Let $\Phi(\al,y,c)$ be the solution of the inhomogeneous Rayleigh equation
with $f(\al,y,c)=\f {\widehat{\om}_0(\al,y)} {i\al(u-c)}$ and $c\in \Om$:
\beq\label{eq:Rayleigh-Ihom-S1}
\left\{
\begin{aligned}
&\Phi''-\al^2\Phi-\frac{u''}{u-c}\Phi=f,\\
&\Phi(-1)=\Phi(1)=0.
\end{aligned}
\right.
\eeq
Then we find that
\beno
(c-\mathcal{R}_\al)^{-1}\widehat{\psi}(0,\al,y)=i\al\Phi(\al,y,c).
\eeno

To study the decay of the continuous spectrum, one of key ingredients is to establish the limiting absorption principle:
\beno
\Phi(\al,y,c\pm i\varepsilon)\to \Phi_\pm(\al,y,c)\quad \text{for}\,\, c\in \text{Ran}\,u,
\eeno
as $\varepsilon\to 0$. If $u$ has critical points, then $\f 1 {u(y)-c}$ is more singular, thus the problem
becomes highly nontrivial. Indeed, the authors (in \cite{BM}, P. 952) pointed out that

{\it These results (the limit and its properties) are the difficult aspects of the discussion, from a mathematical point of view.}\\
The arguments  in \cite{BM} are based on the hypothesis of the limiting absorption principle, which was verified by using  numerical computations.

Using blow-up analysis and compactness argument, we establish the uniform $H^1$ estimate of the solution $\Phi(\al,y,c\pm i\varepsilon)$ for $c\in \text{Ran}\,u$ and $\varepsilon>0$ and the limiting absorption principle: there exists $\Phi_\pm(\al,y,c)\in H^1_y$ so that as $\varepsilon\to 0$,
\beno
\Phi(\al,y,c\pm i\varepsilon)\to \Phi_\pm(\al,y,c)\quad \text{in}\,\,C([-1,1]).
\eeno
With these information, we are in a position to prove the linear damping and the vorticity depletion phenomena for a class of shear flows in $\mathcal{K}$, i.e., Theorem \ref{thm:general} and Theorem \ref{thm:vorticity}.
See section 6 and section 7.

To obtain the explicit decay rate of the velocity, we need to know more precise behaviour of the limit function $\Phi_\pm(\al,y,c)$. To this end, we consider a class of symmetric flows. The main advantage is that we can decompose the solution of \eqref{eq:Rayleigh-Ihom-S1} into the odd part and even part due to the symmetry of $u(y)$:
\beno
\left\{
\begin{aligned}
&\Phi_{o}''-\al^2\Phi_{o}-\frac{u''}{u-c}\Phi_{o}=f_{o},\\
&\Phi_{o}(0)=\Phi_{o}(1)=0,
\end{aligned}
\right.
\eeno
and
\beno
\left\{
\begin{aligned}
&\Phi_{e}''-\al^2\Phi_{e}-\frac{u''}{u-c}\Phi_{e}=f_e,\\
&\Phi_{e}'(0)=\Phi_{e}(1)=0,
\end{aligned}
\right.
\eeno
where $f_o$ and $f_e$ are the odd part and even part of $f$ respectively. These two equations can be dealt as in the monotonic case.
Main difference is that $\f 1 {u(y)-c}$ is more singular for $c=u(0)$. The price to pay is that we can only establish the uniform estimates
of the solution of the homogeneous Rayleigh equation in the weighted space(see section 4), which also lead to some essential difficulties.

To solve the inhomogeneous Rayleigh equation, a key point is to prove that the Wronskian of the solution does not vanish for the homogeneous Rayleigh equations
in the odd and even case. In section 5, we show that this key fact is equivalent to our spectral assumption: the linearized operator has no embedding eigenvalues.

For the symmetric flow, we give the precise formula of $\Phi_\pm(\al, y,c)$  in section 8. Thus, we obtain the representation formula of the stream function
\begin{align*}
\widehat{\psi}(t,\al,y)
=&\frac{1}{2\pi}\int_{u(0)}^{u(1)}\al\widetilde{\Phi}_o(y,c)e^{-i\al ct}dc+\frac{1}{2\pi}\int_{u(0)}^{u(1)}\al\widetilde{\Phi}_e(y,c)e^{-i\al ct}dc\nonumber\\
=&\widehat{\psi}_o(t,\al,y)+\widehat{\psi}_e(t,\al,y).
\end{align*}
In the monotonic case, we first derive the formula of the vorticity from the stream function. Then we prove that the vorticity is bounded in Sobolev space. Finally, the decay estimate of the velocity is proved by using a dual argument.
In the present case, we find that it is difficult to follow the same procedure. New idea is that we directly derive the decay estimates of the velocity based on the formulation of stream function and the dual argument. More precisely, in section 9, we will derive the following two important formulas: for $f=g''-\al^2g$ with $g\in H^2(0,1)\cap H_0^1(0,1)$,
\begin{align*}
&\int_0^1\widehat{\psi}_o(t,\al,y)f(y)dy
=-\int_{u(0)}^{u(1)}K_o(c,\al)e^{-i\al ct}dc,
\end{align*}
and for $f=g''-\al^2g$ with $g\in H^2(0,1)$ and $g'(0)=g(1)=0$,
\begin{align*}
&\int_0^1\widehat{\psi}_e(t,\al,y)f(y)dy=-\int_{u(0)}^{u(1)}K_e(c,\al)e^{-i\al ct}dc,
\end{align*}
where
\beno
&&K_o(c,\al)=\f{\Lambda_1(\widehat{\om}_{o})(c)\Lambda_2(g)(c)}{(\mathrm{A}(c)^2+\mathrm{B}(c)^2){u'(y_c)}},\\
&&K_e(c,\al)=\f{\Lambda_3(\widehat{\omega}_{e})(c)\Lambda_4(g)(c)}{u'(y_c)(\rmA_2^2+\rmB_2^2)(c)}.
\eeno
If $K_0$ and $K_e$ are bounded in $c$ in the Sobolev space $W^{2,1}$, then it is easy to derive the decay estimate of the velocity.
However, the fact that the factor $u'(y_c)$ vanishes for $y_c=0(c=u(0))$ leads to the essential difficulties. To cancel this singularity,
we have to make very precise estimates(especially using the cancellation properties) for some quantities and singular integral operators, which are main tasks of section 10-section 12 and the appendix, where we introduce many tricks to deal with singular integral operators.


\section{The Homogeneous Rayleigh Equation for symmetric flow}
In this section, we solve the homogeneous Rayleigh equation on $[-1,1]$:
\ben\label{eq:Rayleigh-H}
(u-c)(\phi''-\al^2\phi)-u''\phi=0,
\een
where the complex constant $c\in $ will be taken in four kinds of domains:
\beno
&&D_0=\big\{c\in (u(0),u(1))\big\},\\
&&D_{\epsilon_0}=\big\{c=c_r+i\epsilon,~c_r\in (u(0),u(1)), 0<|\epsilon|<\epsilon_0\big\},\\
&&B_{\epsilon_0}^l=\big\{c=u(0)+\epsilon e^{i\theta},~0<\epsilon<\epsilon_0,~\frac{\pi}{2}\leq \theta\leq\frac{3\pi}{2}\big\},\\
&&B_{\epsilon_0}^r=\big\{c=u(1)-\epsilon e^{i\theta},~0<\epsilon<\epsilon_0,~\frac{\pi}{2}\leq \theta\leq\frac{3\pi}{2}\big\},
\eeno
for some $\epsilon_0\in (0,1)$. We denote
\ben\label{eq:Omega}
\Om_{\epsilon_0}=D_0\cup D_{\epsilon_0}\cup B_{\epsilon_0}^l\cup B^r_{\epsilon_0}.
\een
For $c\in \Om_{\epsilon_0}$, we denote
\ben\label{def:cr}
c_r=\textrm{Re}\, c \quad \textrm{for }c\in D_{\epsilon_0},\quad c_r=u(0)\quad \textrm{for }c\in B_{\epsilon_0}^l,
\quad c_r=u(1)\quad \textrm{for }c\in B_{\epsilon_0}^r.
\een
Then we define $y_c\in [0,1]$ so that $u(y_c)=c_r$.\smallskip

In this section, we always assume that $u\in C^4([-1,1])$ satisfies $(S)$.

\subsection{Basic properties of symmetric flow}

We introduce
\ben\label{def:v}
v(y)=\sqrt{u(y)-u(0)}\quad \text{for}\,\, y\in [0,1]\quad \text{and}\quad v(y)=-v(-y)\quad \text{for}\,\,y\in [-1,0].
\een

\begin{lemma}\label{lem:u}
It holds that
\begin{itemize}
\item[1.]  there exists $c_0>0$ such that
\beno
u'(y)\geq c_0y\quad \text{for}\,\, |y|\in [0,1];
\eeno

\item[2.] the function $v\in C^4([-1,0)\cup(0,1])$ satisfies
\beno
v\in C^{3}([-1,1]),\quad yv(y)\in C^4([-1,1]), \quad\tc v^{-1}(\tc)\in C^4([-v(1),v(1)]),
\eeno
and there exists $c_1>0$ so that $v'(y)\ge c_1$;

\item[3.] there exists $C>0$ so that for any $1\geq y\geq y'\geq 0$,
\beno
&C^{-1}(y-y')(y+y')\leq u(y)-u(y')\leq C(y-y')(y+y')\quad\text{or} \\
&C^{-1}\leq\frac{(u'(y)+u'(y'))(y-y')}{u(y)-u(y')}\leq C.
\eeno
\end{itemize}
\end{lemma}

\begin{proof}
The first property is obvious. Thanks to $u'(0)=0$, we get
\begin{align*}
v(y)^2=u(y)-u(0)
=\int_0^yu'(y')\,dy'=\int_0^y\int_0^{y'}u''(y'')\,dy''\,dy'
=y^2m(y),
\end{align*}
where $m(y)=\int_0^1(1-t)u''(ty)\,dt$. Due to $u''(0)>0$, there exits $\d_0>0$ and $c_0>0$ so that for $y\in [0,\d_0]$,
\beno
 m(y)\geq \int_0^1c_0(1-t)\,dt\geq \f12c_0.
\eeno
For $y\in[\d_0,1]$,
\beno
m(y)\geq y^2m(y)=u(y)-u(0)\geq \int_0^{\d_0}c_0y\,dy\geq \f{c_0\d_0^2}{2}.
\eeno
Thus, $m(y)\geq \f{c_0\d_0^2}{2}$ and $v(y)=ym(y)^{\f12}\in C^{2}([-1,1])$.  So,
\beno
v'(y)=m(y)^{\f12}+\f{ym'(y)}{2m(y)^{\f12}},
\eeno
where
\begin{align*}
ym'(y)=\f{u'(y)-2ym(y)}{y}=\int_0^{1}u''(ty)\,dt-2m(y)\in C^{2}([-1,1])
\end{align*}
with the bound $\|ym'\|_{L^{\infty}}\leq 5\|u\|_{C^2}\leq C$.
Thus, we get $v'\in C^{2}([-1,1])$, so $v\in C^{3}([-1,1])$. Moreover,
for $y\in \big[0,\min\big\{\d_0,\frac{c_0\d_0^2}{C_0}\big\}\big]$,
\beno
v'(y)\geq \sqrt{\f{c_0\d_0^2}{2}}-\f{C_0y}{2\sqrt{\f{c_0\d_0^2}{2}}}\geq \f12\sqrt{\f{c_0\d_0^2}{2}}.
\eeno
For $y\in \big[\min\big\{\d_0,\frac{c_0\d_0^2}{C_0}\big\},1\big]$, we have $u'(y)\geq c_0\min\big\{\d_0,\frac{c_0\d_0^2}{C_0}\big\}$, thus,\begin{align*}
v'(y)=\frac{u'(y)}{2ym(y)^{\f12}}\geq c_1>0
\end{align*}
for some $c_1>0$.

Because of $6(v'')^2+8v'v'''+2vv''''=u''''$, we have
\beno
2yv''''(y)=\f{y}{v(y)}(u''''-6(v'')^2-8v'v''')(y)\in C([-1,1]),
\eeno
which shows that $yv(y)\in C^4([-1,1])$ and $\tc v^{-1}(\tc)\in C^{4}([-v(1),v(1)])$.

Finally, for $1\geq y\geq y'\geq 0$, we have
\begin{align*}
u(y)-u(y')&=(v(y)-v(y'))(v(y)+v(y'))\\
&=(y-y')(y+y')\Big(\int_0^1v'(y'+t(y-y'))dt\Big)
\Big(\int_0^1v'(-y'+t(y+y'))dt\Big),
\end{align*}
which implies the third property.
\end{proof}

\subsection{Rayleigh integral operator}
Given $|\al|\ge 1$, let $A$ be a constant larger than $C|\al|$ with $C\ge 1$ only depending on $c_0$ and $\|u\|_{C^4}$. Unlike the monotonic case, we need to introduce the weighted functional spaces.

\begin{definition}\label{def:X-space}
For a function $f(y,c)$ defined on $[0,1]\times \Om_{\epsilon_0}$ and $k\ge 0$, we define
\beno
&&\|f\|_{X^k_0}\eqdef\sup_{(y,c)\in [0,1]\times D_0}\bigg|\frac{u'(y_c)^kf(y,c)}{\cosh(A(y-y_c))}\bigg|,\\
&&\|f\|_{{X}^k}\eqdef\sup_{(y,c)\in [0,1]\times D_{\epsilon_0}}\bigg|\frac{u'(y_c)^kf(y,c)}{\cosh(A(y-y_c))}\bigg|,\\
&&\|f\|_{X_l^k}\eqdef\sup_{(y,c)\in [0,1]\times B_{\epsilon_0}^l}\bigg|\frac{u'(y)^kf(y,c)}{\cosh(Ay)}\bigg|,\\
&&\|f\|_{X_r}\eqdef\sup_{(y,c)\in [0,1]\times B_{\epsilon_0}^r}\bigg|\frac{f(y,c)}{\cosh(A(y-1))}\bigg|.
\eeno
\end{definition}

\begin{definition}\label{def:C^1/2-space}
For a function $f(y,c)$ defined on $[0,1]\times \Om_{\ep_0}\setminus D_0$, we define
\beno
&&\|f\|_{C_{D_{\ep_0}}^{1/2}}\eqdef\sup_{(y,y_c)\in [0,1]\times [0,1]}\sup_{\ep_1\neq \ep_2,\ {\ep_1\ep_2\geq 0}}\bigg|\frac{f(y,u(y_c)+i\ep_1)-f(y,u(y_c)+i\ep_2)}{|\ep_1-\ep_2|^{\f12}\cosh(A(y-y_c))}\bigg|,\\
&&\|f\|_{C_l^{1/2}}\eqdef\sup_{(y,\th)\in [0,1]\times [\pi/2,3\pi/2]}\sup_{\ep_1\neq \ep_2,\ {\ep_1\ep_2\geq 0}}\bigg|\frac{f(y,u(0)+\ep_1e^{i\th})-f(y,u(0)+\ep_2e^{i\th})}{|\ep_1-\ep_2|^{\f12}\cosh(Ay)}\bigg|.
\eeno
\end{definition}

\begin{definition}
For a function $f(y,c)$ defined on $[0,1]\times \Om_{\epsilon_0}$, we define
\beno
&&\|f\|_{Y_0}\eqdef \sum_{k=0}^2\sum_{\b+\gamma=k}A^{-k}\|\partial_y^{\b}\partial_c^{\gamma}f\|_{X^{\gamma}_0}+{A^{-3}\|\pa_c^2\pa_yf\|_{X^3_1},}\\
&&\|f\|_{Y}\eqdef\|f\|_{{X}^0}+\frac{1}{A}\big(\|\partial_yf\|_{{X}^1}+\|\partial_{c}f\|_{{X}^1}+\|f\|_{C_{D_{\ep_0}}^{1/2}}\big),\\
&&\|f\|_{Y_l}\eqdef\|f\|_{X_l^0}+\frac{1}{A}\big(\|\partial_yf\|_{X_l^1}+\|\partial_{c}f\|_{X_l^1}+\|f\|_{C_l^{1/2}}\big),\\
&&\|f\|_{Y_r}\eqdef\|f\|_{X_r^0}+\frac{1}{A}\big(\|\partial_yf\|_{X_r}+\|\partial_{c}f\|_{X_r}\big),
\eeno
where
\beno
\|f\|_{X^3_1}\eqdef\sup_{(y,c)\in [0,1]\times D_0}\bigg|\frac{Au'(y_c)^2f(y,c)}{\cosh(A(y-y_c))(A+y_c/(y+y_c)^2)}\bigg|.
\eeno
\end{definition}

Now we introduce the Rayleigh integral operator, which will be used to give the solution formula
of the homogeneous Rayleigh equation.

\begin{definition}
The Rayleigh integral operator $T$ is defined by
\beno
T\triangleq T_0\circ T_{2,2}\eqdef\int_{y_c}^y\f 1 {(u(y')-c)^2}\int_{y_c}^{y'}f(z,c)(u(z)-c)^2dzdy',
\eeno
where
\beno
&&T_0f(y,c)\eqdef\int_{y_c}^yf(z,c)dz,\\
&&T_{k,j}f(y,c)\eqdef\frac{1}{(u(y)-c)^j}\int_{y_c}^yf(z,c)(u(z)-c)^kdz.
\eeno
\end{definition}

\begin{lemma}\label{lem:T-bound}
There exists a constant $C$ independent of $A$ so that
\beno
&&\|Tf\|_{Y_0}\leq \frac{C}{A^2}\|f\|_{Y_0},\quad\|Tf\|_{Y}\leq \frac{C}{A^2}\|f\|_Y,\\
&&\|Tf\|_{Y_l}\leq \frac{C}{A^2}\|f\|_{Y_l},\quad\|Tf\|_{Y_r}\leq \frac{C}{A^2}\|f\|_{Y_r}.
\eeno
\end{lemma}
\begin{proof}
Let us prove the first inequality.  A direct calculation shows that for $k\ge 0$,
\begin{align}
\|T_0f\|_{X^k_0}
=&\sup_{(y,c)\in[0,1]\times D_0}\bigg|
\frac{1}{\cosh A(y-y_c)}\int_{y_c}^y
\frac{u'(y_c)^kf(z,c)}{\cosh A(z-y_c)}\cosh A(z-y_c)dz\bigg|\nonumber\\
\leq &\sup_{(y,c)\in[0,1]\times D_0}\bigg|
\frac{1}{\cosh A(y-y_c)}\int_{y_c}^y
\cosh A(z-y_c)dz\bigg|\|f\|_{X^k_0}\nonumber\\
\leq &\frac{1}{A}\|f\|_{X^k_0}.\label{eq:T0-est}
\end{align}
Using the fact that for any $0<y<y'<y_c$ or $y_c<y'<y<1$
\ben\label{eq:u-mono}
|u(y')-c|\leq |u(y)-c|,
\een
we deduce that
\begin{align}
\|T_{2,2}f\|_{X^k_0}
\leq &\sup_{(y,c)\in[0,1]\times D_0}
\bigg|\frac{y-y_c}{\cosh A(y-y_c)}
\int_0^1
\cosh tA(y-y_c)dt
\bigg|\|f\|_{X^k_0}\nonumber\\
\leq& \frac{C}{A}\|f\|_{X^k_0},\label{eq:T22-est}
\end{align}
from which and \eqref{eq:T0-est}, we deduce that for $k\ge 0$,
\ben\label{eq:T-bound-1}
\|Tf\|_{X^k_0}\le \frac C {A^2}\|f\|_{X^k_0}.
\een

By Lemma \ref{lem:u}, we have
\beno
0\leq\frac{(u'(y)+u'(y_c))(y-y_c)}{u(y)-u(y_c)}\leq C,
\eeno
which along with \eqref{eq:u-mono} implies that
\begin{align}
&\|T_{k,k+1}f\|_{X^{j+1}_0}+\|u'(y)T_{k,k+1}f\|_{X^j_0}\nonumber\\
&\leq C\sup_{(y,c)\in[0,1]\times D_0}
\bigg|\frac{1}{\cosh A(y-y_c)}\int_0^1
\cosh tA(y-y_c)dt
\bigg|\|f\|_{X^j_0}\nonumber\\
&\leq  C\|f\|_{X^j_0}.\label{eq:Tk-est}
\end{align}
It is easy to see that
\beno
&&\partial_yTf(y,c)=T_{2,2}f(y,c),\\
&&\partial_cTf(y,c)=2T_0\circ T_{2,3}f(y,c)-2T_0\circ T_{1,2}f(y,c)+T\partial_cf(y,c),\\
&&\partial_y^2Tf(y,c)=-2u'(y)T_{2,3}f(y,c)+f(y,c).
\eeno
Thus, it follows from \eqref{eq:T0-est}, \eqref{eq:T22-est} and  \eqref{eq:Tk-est} that
\ben\label{eq:T-bound-2}
&\|\partial_{y}Tf\|_{X_0^0}+\frac{1}{A}\|\partial_{y}^2Tf\|_{X_0^0}+\|\partial_{c}Tf\|_{X_0^1}\leq \frac{C}{A}\|f\|_{X_0^0}
+\frac C {A^2}\|\partial_{c}f\|_{X_0^1}.
\een

For $c\in D_0,$ let
\beno
u_1(y,c)=\f {u(y)-c} {y-y_c}=\int_0^1u'(y_c+t(y-y_c))dt.
\eeno
Then we have
\begin{align*}
T_{k,k+1}f(y,c)
&=\int_0^1f(y_c+t(y-y_c),c)\frac{t^{k}u_1(y_c+t(y-y_c),c)^{k}}{u_1(y,c)^{k+1}}dt.
\end{align*}
By Lemma \ref{lem:u}, $u_1(y,c)\sim y+y_c$. Then we can deduce that
\begin{align}
\left|\frac{u'(y_c)^j\pa_cT_{k,k+1}f(y,c)}{\cosh A(y-y_c)}\right|
&\leq C\left(\frac{\|f\|_{X^{j-1}_0}}{u_1(y,c)^2}+\frac{\|\pa_{y}f\|_{X^{j-1}_0}+\|\pa_{c}f\|_{X^{j}_0}}{u_1(y,c)}\right).\label{eq:Tk-est2}
\end{align}
Direct calculation gives
\begin{align*}
\partial_y\partial_cTf(y,c)
&=2T_{2,3}f(y,c)-2T_{1,2}f(y,c)+T_{2,2}\partial_cf(y,c),\\
\partial_y\partial_c^2Tf(y,c)
&=2\partial_cT_{2,3}f(y,c)-2\partial_cT_{1,2}f(y,c)\\
&\quad+2T_{2,3}\partial_cf(y,c)
-2T_{1,2}\partial_cf(y,c)+T_{2,2}\partial_c^2f(y,c),
\end{align*}
and
\begin{align*}
\partial_c^2Tf(y,c)=&2\pa_cT_0T_{2,3}f(y,c)-2\pa_cT_0T_{1,2}f(y,c)+\partial_cT\partial_cf(y,c)\\
=&2T_0\partial_cT_{2,3}f(y,c)-2T_0\partial_cT_{1,2}f(y,c)+\frac{f(y_c,c)}{3u'(y_c)^2}\\
&+2T_0T_{2,3}\partial_cf(y,c)-2T_0T_{1,2}\partial_cf(y,c)+T\partial_c^2f(y,c).
\end{align*}
Then it follows from \eqref{eq:Tk-est} and \eqref{eq:T22-est} that
\ben
\|\partial_y\partial_cTf\|_{X_{0}^1}\leq C\|f\|_{X_0^0}+\frac C {A}\|\partial_{c}f\|_{X_1}.\label{eq:T-est3}
\een
By \eqref{eq:Tk-est2}, we have
\begin{align*}
\left|\frac{u'(y_c)^2T_0\pa_cT_{k,k+1}f(y,c)}{\cosh A(y-y_c)}\right|
\le& \f {u'(y_c)} {\cosh A(y-y_c)}\int_{y_c}^y\cosh A(y'-y_c)\left|\frac{u'(y_c)\pa_cT_{k,k+1}f(y',c)}{\cosh A(y'-y_c)}\right|dy'\\
\le& Cu'(y_c)\int_{y_c}^y\f 1 {u_1(y',c)^2}dy'\|f\|_{X^0_0}+\f C A\big(\|\partial_{y}f\|_{X_0^0}+\|\partial_{c}f\|_{X_0^1}\big)\\
\le& C\|f\|_{X^0_0}+\f C A\big(\|\partial_{y}f\|_{X_0^0}+\|\partial_{c}f\|_{X_0^1}\big),
\end{align*}
which along with \eqref{eq:T0-est} and \eqref{eq:Tk-est} gives
\ben\label{eq:T-est5}
\|\partial_c^2Tf\|_{X^{2}_0}\leq C\|f\|_{X^0_0}+\frac C {A}\|\partial_{y}f\|_{X_0^0}+\frac C {A}\|\partial_{c}f\|_{X^1_0}+\frac C {A^2}\|\partial_{c}^2f\|_{X^2_0}.
\een
We infer from \eqref{eq:Tk-est} and \eqref{eq:Tk-est2} that
\beno
\left|\frac{u'(y_c)\partial_y\partial_c^2Tf(y,c)}{\cosh A(y-y_c)}\right|\leq \frac {C \|f\|_{X^0_0}}{u_1(y,c)^2}+ \frac{C \|\partial_{y}f\|_{X^0_0}}{u'(y_c)}+ \frac{C \|\partial_{c}f\|_{X^1_0}}{u'(y_c)}+\frac {C \|\partial_{c}^2f\|_{X^2_0}}{Au'(y_c)},
\eeno
which implies that
\ben \label{eq:T-est5}
\|\pa_y\pa_c^2Tf\|_{X^{3}_1}\le CA\|f\|_{X^0_0}+ C\|\partial_{y}f\|_{X^0_0}+ C \|\partial_{c}f\|_{X^1_0}+\frac C {A}\|\partial_{c}^2f\|_{X^2_0}.
\een

Putting (\ref{eq:T-bound-1}), \eqref{eq:T-bound-2} and (\ref{eq:T-est3})-\eqref{eq:T-est5} together, we conclude the first inequality.
\smallskip

The other inequalities of the lemma  can be deduced in a similar derivation as in \eqref{eq:T-bound-2} except the H\"{o}lder estimates $\|Tf\|_{C^{1/2}_{D_{\ep_0}}}$ and $\|Tf\|_{C^{1/2}_l}$. For $c_k=u(y_c)+i\ep_k$ with $0\leq \ep_2<\ep_1$ and $k=1,2$,  we have
\begin{align*}
|Tf(y,c_1)-Tf(y,c_2)|
&=\int_{y_c}^y\int_{y_c}^{y'}\Big(\f{(u(z)-c_1)^2}{(u(y')-c_1)^2}
-\f{(u(z)-c_2)^2}{(u(y')-c_2)^2}\Big)f(z,c_1)dzdy'\\
&\quad+\int_{y_c}^y\int_{y_c}^{y'}\f{(u(z)-c_2)^2}{(u(y')-c_2)^2}\big(f(z,c_1)-f(z,c_2)\big)dzdy'.
\end{align*}
Using the facts that for $y_c\le z\le y'$ or $y'\le z\le y_c$,
\begin{align*}
\Big|\f{(u(z)-c_1)^2}{(u(y')-c_1)^2}
-\f{(u(z)-c_2)^2}{(u(y')-c_2)^2}\Big|
&\leq C\Big|\f{(u(z)-c_1)}{(u(y')-c_1)}-\f{(u(z)-c_2)}{(u(y')-c_2)}\Big|\\
&\leq C\f{|\ep_1-\ep_2||u(z)-u(y')|}{|u(y')-u(y_c)+i\ep_1||u(y')-u(y_c)+i\ep_2|}\\
&\leq C\f{\ep_1-\ep_2}{|u_1(y',c)(y'-y_c)|+|\ep_1|},
\end{align*}
and $u_1(y',c)\geq C^{-1}(y'+y_c)$ and $\ep_1-\ep_2\leq \ep_1$, we deduce that
\begin{align*}
&\f{1}{\cosh A(y-y_c)}\left|\int_{y_c}^y\int_{y_c}^{y'}\Big(\f{(u(z)-c_1)^2}{(u(y')-c_1)^2}
-\f{(u(z)-c_2)^2}{(u(y')-c_2)^2}\Big)f(z,c_1)dzdy_1\right|\\
&\leq C\f{1}{\cosh A(y-y_c)}\left|\int_{y_c}^y\f{(\ep_1-\ep_2)|y'-y_c|}{|u_1(y',c)(y'-y_c)|+|\ep_1|}\cosh A(y'-y_c)dy'\right|\|f\|_{X^0}\\
&\leq \f{|\ep_1-\ep_2|^{\f12}}{\cosh A(y-y_c)}\left|\int_{y_c}^y\cosh A(y'-y_c)dy'\right|\|f\|_{X^0}\\
&\leq \f{C}{A}|\ep_1-\ep_2|^{\f12}\|f\|_{X^0}.
\end{align*}
For the second term, we have
\begin{align*}
&\left|\int_{y_c}^y\int_{y_c}^{y'}\f{(u(z)-c_2)^2}{(u(y')-c_2)^2}\big(f(z,c_1)-f(z,c_2)\big)dzdy'\right|\\
&\leq C|\ep_1-\ep_2|^{\f12}\left|\int_{y_c}^y\int_{y_c}^{y'}\f{(u(z)-c_2)^2}{(u(y')-c_2)^2}\cosh A(z-y_c)dzdy'\right|\|f\|_{C_{D_{\ep_0}}^{1/2}}\\
&\leq \f{C\cosh A(y-y_c)}{A^2}|\ep_1-\ep_2|^{\f12}\|f\|_{C_{D_{\ep_0}}^{1/2}}.
\end{align*}
This shows that
\beno
\|Tf\|_{C^{1/2}_{D_{\ep_0}}}\le \f{C}{A}\|f\|_{X^0}+\f{C}{A^2}\|f\|_{C^{1/2}_{D_{\ep_0}}}\le \f C {A}\|f\|_{Y}.
\eeno

Similarly, for $c_k=u(0)+\ep_ke^{i\th}$ with $k=1, 2$, we have
\begin{align*}
|Tf(y,c_1)-Tf(y,c_2)|
&=\int_{0}^y\int_{0}^{y'}\Big(\f{(u(z)-c_1)^2}{(u(y')-c_1)^2}
-\f{(u(z)-c_2)^2}{(u(y')-c_2)^2}\Big)f(z,c_1)dzdy'\\
&\quad+\int_{0}^y\int_{0}^{y'}\f{(u(z)-c_2)^2}{(u(y')-c_2)^2}\big(f(z,c_1)-f(z,c_2)\big)dzdy'.
\end{align*}
Using the fact that for $0\leq z\leq y'$,
\begin{align*}
\Big|\f{(u(z)-c_1)^2}{(u(y')-c_1)^2}
-\f{(u(z)-c_2)^2}{(u(y')-c_2)^2}\Big|
&\leq C\Big|\f{(u(z)-c_1)}{(u(y')-c_1)}-\f{(u(z)-c_2)}{(u(y')-c_2)}\Big|\\
&\leq C\f{|c_1-c_2||u(z)-u(y')|}{|u(y')-u(0)+\ep_1e^{i\th}||u(y')-u(0)+\ep_2e^{i\th}|}\\
&\leq C\f{\ep_1-\ep_2}{y'^2+|\ep_1|},
\end{align*}
we deduce that
\begin{align*}
&\f{1}{\cos Ay}\left|\int_{0}^y\int_{0}^{y'}\Big(\f{(u(z)-c_1)^2}{(u(y')-c_1)^2}
-\f{(u(z)-c_2)^2}{(u(y')-c_2)^2}\Big)f(z,c_1)dzdy'\right|\\
&\leq C\f{1}{\cosh Ay}\left|\int_{0}^y\f{(\ep_1-\ep_2)y'}{y'^2+|\ep_1|}\cosh Ay'dy'\right|\|f\|_{X^0_l}\\
&\leq \f{C}{A}|\ep_1-\ep_2|^{\f12}\|f\|_{X^0_l},
\end{align*}
and
\begin{align*}
&\left|\int_{0}^y\int_{0}^{y'}\f{(u(z)-c_2)^2}{(u(y')-c_2)^2}\big(f(z,c_1)-f(z,c_2)\big)dzdy'\right|\\
&\leq C|\ep_1-\ep_2|^{\f12}\left|\int_{0}^y\int_{0}^{y'}\f{(u(z)-c_2)^2}{(u(y')-c_2)^2}\cosh Azdzdy'\right|\|f\|_{C_{l}^{1/2}}\\
&\leq \f{C\cosh Ay}{A^2}|\ep_1-\ep_2|^{\f12}\|f\|_{C_{l}^{1/2}}.
\end{align*}
This shows that
\beno
\|Tf\|_{C^{1/2}_{l}}\le \f C {A}\|f\|_{Y_l}.
\eeno

The proof of the lemma is completed.
\end{proof}

\subsection{Existence of the solution}

In this subsection, the constant $C$ may depend on $\al$.

\begin{proposition}\label{prop:Rayleigh-Hom}
There exists a solution $\phi(y,c)\in  C^1\big([0,1]\times \Omega_{\epsilon_0}\setminus D_0\big)\cap C\big([0,1]\times \Omega_{\epsilon_0}\big)$
of the Rayleigh equation \eqref{eq:Rayleigh-H}. Let
\ben\label{def:phi1}
\phi_1(y,c)=\frac {\phi(y,c)} {u(y)-c}.
\een
There exists $\epsilon_1>0$ such that for any $\epsilon_0\in[0,\epsilon_1)$ and $(y,c)\in [0,1]\times \Om_{\epsilon_0}$,
\beno
&&|\phi_{1}(y,c)|\ge \f12,\quad |\phi_{1}(y,c)-1|\leq C |y-y_c|^2.
\eeno
\end{proposition}

\begin{lemma}\label{lem:Ray-D}
Let $c\in D_{\epsilon_0}$.
Then there exists a solution $\phi(x,c)\in Y$ to the Rayleigh equation
\beno
\left\{
\begin{array}{l}
\phi''-\al^2\phi-\frac{u''}{u-c}\phi=0,\\
\frac{\phi(y_c,c)}{u(y_c)-c}=1,\quad \Big(\frac{\phi(y,c)}{u(y)-c}\Big)'\big|_{y=y_c}=0.
\end{array}\right.
\eeno
Moreover, there holds
\begin{align*}
\|\phi_1\|_{Y}+\|\phi\|_{Y}\leq C.
\end{align*}
\end{lemma}

\begin{lemma}\label{Lem:Ray-Bl}
Let $c\in B^l_{\epsilon_0}$.
Then there exists a solution $\phi(x,c)\in Y_l$ to the Rayleigh equation
\beno
\left\{
\begin{array}{l}
\phi''-\al^2\phi-\frac{u''}{u-c}\phi=0,\\
\frac{\phi(0,c)}{u(0)-c}=1,\quad \Big(\frac{\phi(y,c)}{u(y)-c}\Big)'\big|_{y=0}=0.
\end{array}\right.
\eeno
Moreover, there holds
\begin{align*}
\|\phi_1\|_{Y_l}+\|\phi\|_{Y_l}\leq C.
\end{align*}
\end{lemma}

\begin{lemma}\label{Lem:Ray-Br}
Let $c\in B^r_{\epsilon_0}$.
Then there exists a solution $\phi(x,c)\in Y_r$ to the Rayleigh equation
\beno
\left\{
\begin{array}{l}
\phi''-\al^2\phi-\frac{u''}{u-c}\phi=0,\\
\frac{\phi(1,c)}{u(1)-c}=1,\quad \Big(\frac{\phi(y,c)}{u(y)-c}\Big)'\big|_{y=1}=0.
\end{array}\right.
\eeno
Moreover, there holds
\begin{align*}
\|\phi_1\|_{Y_r}+\|\phi\|_{Y_r}\leq C.
\end{align*}
\end{lemma}

\begin{lemma}\label{lem:Ray-D0}
Let $c\in D_0$.
Then there exists a solution $\phi(y,c)\in Y_0$ to the Rayleigh equation
\beno
\left\{
\begin{array}{l}
\phi''-\al^2\phi-\frac{u''}{u-c}\phi=0,\\
\phi(y_c,c)=0,\quad \phi'(y_c,c)=u'(y_c).
\end{array}\right.
\eeno
Moreover, there holds
\begin{align*}
\|\phi_1\|_{Y_0}+\|\phi\|_{Y_0}\leq C.
\end{align*}
\end{lemma}

With Lemma \ref{lem:T-bound} and the formula
\begin{align}
\phi_1(y,c)=&1+\int_{y_c}^y\frac{\al^2}{(u(y')-c)^2}\int_{y_c}^{y'}\phi_{1}(z,c)(u(z)-c)^2dzdy'\nonumber\\
=&1+\al^2T\phi_1(y,c),\label{eq:phi1}
\end{align}
the above lemmas can be proved as in Lemma 4.6-Lemma 4.9 in \cite{WZZ1}. So, we omit the details.
\smallskip

As in the proof of Proposition 4.5 in \cite{WZZ1}, we define
\beno
\phi(y,c)\eqdef\left\{
\begin{array}{l}
\phi^0(y,c)\quad \textrm{for } c\in D_0,\\
\phi^\pm(y,c)\quad \textrm{for } c\in D_{\epsilon_0},\\
\phi^l(y,c)\quad \textrm{for } c\in B_{\epsilon_0}^l,\\
\phi^r(y,c)\quad \textrm{for } c\in B_{\epsilon_0}^r,
\end{array}
\right.
\eeno
where $\phi^\pm$, $\phi^l, \phi^r, \phi^0$ are given by Lemma \ref{lem:Ray-D}, Lemma \ref{Lem:Ray-Bl}, Lemma \ref{Lem:Ray-Br} and Lemma \ref{lem:Ray-D0} respectively. Then $\phi(y,c)$ is our desired solution.
\medskip

We need the following further properties of $\phi_1(y,c)$.

\begin{lemma}\label{lem:phi1}
For $c\in D_0$, it holds that
\begin{itemize}
\item[1.] $\phi_1(y,c)\geq 1$;
\item[2.] $\pa_y\phi_1(y,c)>0$ for $y\in (y_c,1]$ and $\pa_y\phi_1(y,c)<0$ for $y\in [0,y_c)$;
\item[3.] for any given $M_0> 0$, there exists a constant $C$ only depending on $M_0$ such that for
$\b+\g\leq 2$, $|y-y_c|\leq M_0/|\al|$
\beno
&&|u'(y_c)^\gamma\partial_y^{\b}\partial_c^{\gamma}\phi_1(y,c)|\leq C|\al|^{\beta+\gamma},\\
&&|u'(y_c)^2\pa_y\pa_c^2\phi_1(y,c)|\leq C(|\al|^3+\al^2y_c/(y+y_c)^2).
\eeno
\end{itemize}
\end{lemma}

\begin{proof}
The first two properties are a direct consequence of \eqref{eq:phi1}. The third property follows from Lemma \ref{lem:Ray-D0} and the definition of $Y_0$ norm.
\end{proof}

\begin{lemma}\label{lem:phi1-con}
For $c_{\ep}\in D_{\ep_0}\cup D_0$ and $c_{\ep}\in B_{\ep}^r$, $\phi_1'(0,c_{\ep})=\pa_y\phi_1(0,c_\ep)$ is continuous in $\ep$. For $c_{\ep}\in B_{\ep}^l$, $\phi_1'(0,c_{\ep})=0$.
\end{lemma}

\begin{proof}
Let $c_{\ep_1}=u(y_c)+i\ep_1\in D_{\ep_0}\cup D_0$ and $c_{\ep_2}=u(y_c)+i\ep_2\in D_{\ep_0}\cup D_0$ with $0\leq \ep_1<\ep_2\leq \ep_0$.
By \eqref{eq:phi1}, we have
\beno
\phi_1'(0,c_{\ep_1})=\f{\al^2}{(u(0)-c_{\ep_1})^2}\int_{y_c}^0(u(y')-c_{\ep_1})^2\phi_1(y',c_{\ep_1})dy'.
\eeno
Therefore,
\begin{align*}
|\phi_1'(0,c_{\ep_1})-\phi_1'(0,c_{\ep_2})|
&\leq C\Big|\int_0^{y_c}\f{(u(y')-c_{\ep_1})^2}{(u(0)-c_{\ep_1})^2}(\phi_1(y',c_{\ep_1})-\phi_1(y',c_{\ep_2}))dy'\Big|\\
&\quad+C\Big|\int_0^{y_c}\Big(\f{(u(y')-c_{\ep_1})^2}{(u(0)-c_{\ep_1})^2}-\f{(u(y')-c_{\ep_1})^2}{(u(0)-c_{\ep_1})^2}\Big)\phi_1(y',c_{\ep_2})dy'\Big|.
\end{align*}
Using the fact that $|\pa_{\ep}\phi_1|\leq \f{C}{u'(y_c)}$ and $\big|\f{(u(y')-c_{\ep_1})^2}{(u(0)-c_{\ep_1})^2}\big|\leq C$ for $y'\le y_c$,  the first term is bounded by $C|\ep_1-\ep_2|$.  For the second term, we have
\begin{align*}
&\Big|\int_0^{y_c}\Big(\f{(u(y')-c_{\ep_1})^2}{(u(0)-c_{\ep_1})^2}-\f{(u(y')-c_{\ep_1})^2}{(u(0)-c_{\ep_1})^2}\Big)\phi_1(y',c_{\ep_2})dy'\Big|\\
&\leq C\int_0^{y_c}\Big|\f{(u(y')-c_{\ep_1})}{(u(0)-c_{\ep_1})}-\f{(u(y')-c_{\ep_1})}{(u(0)-c_{\ep_1})}\Big|dy'\\
&\leq C\int_0^{y_c}\Big|\f{(u(z)-u(0))(c_{\ep_1}-c_{\ep_2})}{(u(0)-c_{\ep_1})(u(0)-c_{\ep_2})}\Big|dz\\
&\leq C\f{y_c(\ep_1-\ep_2)}{y_c^2+\ep_1}\leq C\sqrt{\ep_1-\ep_2}.
\end{align*}
This shows that
\beno
 |\phi_1'(0,c_{\ep_1})-\phi_1'(0,c_{\ep_2})|\leq C\sqrt{\ep_1-\ep_2}.
\eeno

Similarly, we can prove that for $c_{\ep_1}, c_{\ep_2}\in B^r_{\ep}$,
\beno
|\phi_1'(0,c_{\ep_1})-\phi_1'(0,c_{\ep_2})|\leq C|\ep_1-\ep_2|.
\eeno

For $c_\ep\in B_\ep^l$, we have $y_c=0$, thus $\phi_1'(0,c_\ep)=0$.
\end{proof}

\subsection{Uniform estimates of the solution}

In this subsection, we establish some uniform estimates in wave number $\al$ for the solution $\phi(y,c)$ for $c\in D_0$
of the Rayleigh equation given by Lemma \ref{lem:Ray-D0}.
Without loss of generality, we always assume $\al\ge 1$ in the sequel.\smallskip

Recall that $\phi_1(y,c)=\frac {\phi(y,c)} {u(y)-c}$ and $c=u(y_c)$ for $c\in D_0$.
We introduce
\ben
&&\cF(y,c)=\f{\pa_y\phi_1(y,c)}{\phi_1(y,c)},\quad \cG(y,c)=\f{\pa_c\phi_1(y,c)}{\phi_1(y,c)},\label{def:FG}\\
&&\cG_1(y,c)=\f{\cF(y,c)}{u'(y_c)}+\cG(y,c)=\f{1}{\phi_1}\Big(\f{\pa_y}{u'(y_c)}+\pa_c\Big)\phi_1(y,c).
\een
It is easy to see that
\beno
\cF(y_c,c)=\cG(y_c,c)=\cG_1(y_c,c)=0.
\eeno

\begin{proposition}\label{prop:phi1}
There exists a constant $C$ independent of $\al$ such that
\beno
&&\phi_1(y,c)-1\le C\min\{\al^2|y-y_c|^2,1\}\phi_1(y,c),\\
&&C^{-1}e^{C^{-1}\al |y-y_c|}\leq \phi_1(y,c)\leq e^{C\al|y-y_c|},\\
&&C^{-1}\al\min\{\al|y-y_c|,1\}\leq |\cF(y,c)|\leq C\al\min\{\al|y-y_c|,1\},
\eeno
and for $\beta+\gamma\le 2$,
\beno
&&|u'(y_c)^\gamma\partial_y^{\b}\partial_c^{\gamma}\phi_1(y,c)|\leq C \al^{\b+\g}\phi_1(y,c),\\
&&|\pa_c\phi_1(y,c)|\leq C\Big(\f{\al^2|y-y_c|}{u'(y_c)}
+\frac{\al^3|y-y_c|^2}{u'(y_c)}\Big)\phi_1(y,c).
\eeno
Moreover, there holds
\begin{align*}
&\Big|\Big(\frac{\partial_y}{u'(y_c)}+\partial_c\Big)\phi_1(y,c)\Big|\leq C\frac{\min\{1,\al^2|y-y_c|^2\}\phi_1}{u'(y_c)^2},\\
&\Big|\Big(\frac{\partial_y}{u'(y_c)}+\partial_c\Big)^2\phi_1(y,c)\Big|\leq C\frac{\min\{1,\al^2|y-y_c|^2\}\phi_1}{u'(y_c)^4}.
\end{align*}
\end{proposition}
\begin{proof}
{\bf Step 1.} Estimates of $\cF$ and $\pa_y^k\phi_1$ for $k=0,1,2$.
\smallskip

Recall that $\phi_1=1+\al^2T(\phi_1)$, which implies that
\ben\label{eq:phi-est1}
\phi_1(y,c)-1\leq C\min\{\al^2|y-y_c|^2,1\}\phi_1(y,c).
\een

It is easy to check that $\cF$ satisfies
\ben
\cF'+\cF^2+\f{2u'}{u-c}\cF=\al^2,\quad \cF(y_c,c)=0.\label{eq:cF}
\een
Notice that
\beno
\lim_{y\to y_c}\f{u'\cF(y,c)}{u-c}=\pa_y\cF(y_c,c).
\eeno
So, $\pa_y\cF(y_c,c)=\f {\al^2} 3>0$. Thanks to $\cF(y,c)\geq 0$ for $y\geq y_c$ and $\cF(y,c)\leq 0$ for $y\leq y_c$, we also have  $\f{2u'}{u-c}\cF\ge 0$.
Then \eqref{eq:cF} implies that
$|\cF(y,c)|\leq \al$, which in particular gives
\ben
e^{-\al|y-y'|}\leq \f{\phi_1(y',c)}{\phi_1(y,c)}\leq e^{\al|y-y'|}.\label{eq:cF-2}
\een

Using the equation
\begin{align*}
\pa_y\phi_1(y,c)=\f{\al^2}{(u(y)-c)^2}\int_{y_c}^y\phi_1(y',c)(u(y')-c)^2\,dy',
\end{align*}
and the fact that
$(u(y')-c)^2\leq (u(y)-c)^2$ and $\phi_1(y',c)\leq \phi_1(y,c)$  for $y_c\leq y'\leq y$ or $y\leq y'\leq y_c$(by Lemma \ref{lem:phi1}), we infer that
\beno
|\cF(y,c)|\leq C\al^2|y-y_c|.
\eeno
So, we get
\ben\label{eq:cF-up}
|\cF(y,c)|\leq C\al\min\{\al|y-y_c|,1\}.
\een

Using Lemma \ref{lem:u} and \eqref{eq:cF-2}, we deduce that for $0\leq y_c\leq y\leq \f{1}{\al}+y_c\leq 1$,
\begin{align*}
\cF(y,c)&=\f{\al^2}{(u(y)-c)^2}\int_{y_c}^{y}\f{\phi_1(y',c)}{\phi_1(y,c)}(u(y')-c)^2\,dy'\\
&\geq \f{\al^2}{(u(y)-c)^2}\int_{y_c}^{y}e^{-\al|y-y'|}(u(y')-c)^2\,dy'\\
&\geq C^{-1} \f{\al^2}{(y+y_c)^2(y-y_c)^2}\int_{\f{y_c+y}{2}}^y(y'+y_c)^2(y'-y_c)^2\,dy'\\
&\geq C^{-1}\al^2|y-y_c|,
\end{align*}
and for $0\leq y_c-\f{1}{\al}\leq y\leq y_c\leq 1$,
\begin{align*}
-\cF(y,c)&\geq \f{\al^2}{(u(y)-c)^2}\int_{y}^{y_c}e^{-\al|y-y'|}(u(y')-c)^2\,dy'\\
&\geq C^{-1}\f{\al^2}{(y+y_c)^2(y-y_c)^2}\int_{y}^{\f{y_c+y}{2}}(y'+y_c)^2(y'-y_c)^2\,dy'\\
&\geq C^{-1}\al^2|y-y_c|.
\end{align*}
For $0\leq \f{1}{\al}+y_c\leq y\leq 1$, we have
\begin{align*}
\cF(y,c)&=\f{\al^2}{(u(y)-c)^2}\int_{y_c}^{y}\f{\phi_1(y',c)}{\phi_1(y,c)}(u(y')-c)^2\,dy'\\
&\geq \f{\al^2}{(u(y)-c)^2}\int_{y_c}^{y}e^{-\al|y-y'|}(u(y')-c)^2\,dy'\\
&\geq C^{-1}\f{\al^2}{(y+y_c)^2(y-y_c)^2}\int_{\max\big\{\f{y_c+y}{2},y-\f{1}{\al}\big\}}^ye^{-\al|y-y'|}(y'+y_c)^2(y'-y_c)^2\,dy'\\
&\geq C^{-1}\al.
\end{align*}
For $0\leq y\leq y_c-\f{1}{\al}\leq 1$, we have
\begin{align*}
-\cF(y,c)&=\f{\al^2}{(u(y)-c)^2}\int_{y}^{y_c}\f{\phi_1(y',c)}{\phi_1(y,c)}(u(y')-c)^2\,dy'\\
&\geq C^{-1} \f{\al^2}{(u(y)-c)^2}\int_{y}^{y_c}e^{-\al|y-y'|}(u(y')-c)^2\,dy'\\
&\geq C^{-1}\f{\al^2}{(y+y_c)^2(y-y_c)^2}\int^{\max\big\{\f{y_c+y}{2},y+\f{1}{\al}\big\}}_ye^{-\al|y-y'|}(y'+y_c)^2(y'-y_c)^2\,dy'\\
&\geq C^{-1}\al.
\end{align*}
This shows that
\ben
C^{-1}\al\min\{\al|y-y_c|,1\}\leq |\cF(y,c)|\leq C\al\min\{\al|y-y_c|,1\},\label{eq:F-est1}
\een
which along with $\phi_1(y_c,c)=1$ implies that
\ben\label{eq: phi_1_up_low}
&&C^{-1}e^{C^{-1}\al |y-y_c|}\leq \phi_1(y,c)\leq e^{C\al|y-y_c|},\\
&&|\pa_y\phi_1(y,c)|\leq C\al\min\{\al|y-y_c|,1\}\phi_1(y,c).
\een

Using $\phi_1''+\f{2u'}{u-c}\phi_1'=\al^2\phi_1$ and $\big(u'(y)+u'(y_c)\big)|y-y_c|\leq C|u-c|$, we obtain
\beno
|\pa_y^2\phi_1(y,c)|\leq C\al^2\phi_1(y,c).
\eeno

{\bf Step 2}. Estimates of $\pa_c \phi_1, \pa_y\pa_c\phi_1$ and $\Big(\frac{\partial_y}{u'(y_c)}+\partial_c\Big)\phi_1$.

It is easy to check that
\beno
&&\partial_c\cF=\partial_y\cG,\quad \Big(\frac{\partial_y}{u'(y_c)}+\partial_c\Big)\cF=\partial_y\cG_1,\\
&&(\phi^2\partial_c\cF)'+2\phi_1^2u'\cF=0,\quad (\phi^2\partial_y\cF)'+2\phi_1^2(u''(u-c)-u'^2)\cF=0.
\eeno
Then we deduce that
\ben
&&\partial_c\cF(y,c)=\frac{-2}{\phi(y,c)^2}\int_{y_c}^y\phi_1(y',c)^2u'(y')\cF(y',c)dy'=\frac{-2}{\phi(y,c)^2}T_0(\phi_1^2u'\cF),\label{eq:pcF}\\
&&\partial_y\cF(y,c)=\frac{-2}{\phi(y,c)^2}T_0\big(\phi_1^2(u''(u-c)-u'^2)\cF\big),\label{eq:pyF}\\
&&\Big(\frac{\partial_y}{u'(y_c)}+\partial_c\Big)\cF=\partial_y\cG_1=\frac{-2}{\phi(y,c)^2}T_0\Big(\phi_1^2\frac{a_1}{u'(y_c)}\cF\Big),\label{eq:pycF}
\een
where $a_1(y,c)=u'(y)u'(y_c)+u''(y)(u(y)-c)-u'(y)^2$.

It is easy to see that
\beno
&&\partial_ya_1(y,c)=u''(y)(u'(y_c)-u'(y))+u'''(y)(u(y)-c),\\
&&\partial_ca_1(y,c)=u'(y)\frac{u''(y_c)}{u'(y_c)}-u''(y),\quad a_1(y_c,c)=0,
\eeno
which imply that
\ben
|\partial_ya_1(y,c)|\leq C|y-y_c|,\quad |a_1(y,c)|\leq C|y-y_c|^2.\label{eq:a1}
\een
Let
$$a_2(y,c)=u'(y)(u''(y_c)-u''(y))+u'''(y)(u(y)-c).$$
Then we find that
\beno
&&\Big(\frac{\partial_y}{u'(y_c)}+\partial_c\Big)a_1(y,c)=\frac{a_2(y,c)}{u'(y_c)},\\
&&\partial_ya_2(y,c)=u''(y)(u''(y_c)-u''(y))+u''''(y)(u(y)-c),
\eeno
and $a_2(y_c,c)=0$, therefore,
\ben
|\partial_ya_2(y,c)|\leq C|y-y_c|,\quad |a_2(y,c)|\leq C|y-y_c|^2. \label{eq:a2}
\een

By Lemma \ref{lem:phi1},  $\partial_y\phi_1(y,c)$ has the same sign as $y-y_c$. For $y_c\leq y'\leq y$ or $y\leq y'\leq y_c$, $u'(y')\leq C(u'(y)+u'(y_c))$. So, we get by \eqref{eq:pcF} and \eqref{eq:phi-est1} that
\begin{align}\label{eq: pa_cf}
|\partial_y\cG(y,c)|&=|\partial_c\cF(y,c)|\nonumber\\
&\leq \frac{C(u'(y)+u'(y_c))}{\phi(y,c)^2}\int_{y_c}^y\phi_1(y',c)\partial_{y'}\phi_1(y',c)dy'\nonumber\\
&\leq C\frac{\phi_1(y,c)^2-1}{\phi(y,c)^2}(u'(y)+u'(y_c))\nonumber\\
&\leq C\frac{\phi_1(y,c)-1}{\phi_1(y,c)}\f{1}{|y-y_c|^2(u'(y)+u'(y_c))}\nonumber\\
&\leq \frac{C\min\{1,\al^2|y-y_c|^2\}}{|y-y_c|^2(u'(y)+u'(y_c))},
\end{align}
and by \eqref{eq:pycF} and \eqref{eq:a1},
\begin{align}\label{eq: good_1F}
\Big|\Big(\frac{\partial_y}{u'(y_c)}+\partial_c\Big)\cF(y,c)\Big|
&=|\partial_y\cG_1(y,c)|\nonumber\\
&\leq \frac{C|y-y_c|^2}{u'(y_c)\phi(y,c)^2}\int_{y_c}^y\phi_1(y',c)\partial_{y'}\phi_1(y',c)dy'\nonumber\\
&\leq C\frac{\phi_1(y,c)^2-1}{u'(y_c)\phi(y,c)^2}|y-y_c|^2\nonumber\\
&\leq \frac{C\min\{1,\al^2|y-y_c|^2\}}{u'(y_c)}\frac{|y-y_c|^2}{(u(y)-c)^2},
\end{align}
from which and $\cG(y_c,c)=\cG_1(y_c,c)=0$,  it follows that
\beq\label{eq:g}
|\cG(y,c)|\leq C\int_{y_c}^y\frac{\min\{1,\al^2|y'-y_c|^2\}}{|y'-y_c|^2u'(y_c)}dy'\leq \frac{C\al{\min(1,\al|y-y_c|)}}{u'(y_c)},
\eeq
and
\beno
|\cG_1(y,c)|\leq C\frac{\min\{1,\al^2|y-y_c|^2\}}{u'(y_c)}\left|\int_{y_c}^y\frac{|y'-y_c|^2}{(u(y')-c)^2}dy'\right|\leq \frac{C\min\{1,\al^2|y-y_c|^2\}}{u'(y_c)^2}.
\eeno
Thus, we deduce that
\ben
|\partial_c\phi_1(y,c)|\leq \f{C\al \phi_1(y,c)}{u'(y_c)},\quad |\partial_y\partial_c\phi_1(y,c)|\leq \f{C\al^2 \phi_1(y,c)}{u'(y_c)},\label{eq:phi-pc}
\een
and
\ben\label{eq: good_1phi_1}
\Big|\Big(\frac{\partial_y}{u'(y_c)}+\partial_c\Big)\phi_1(y,c)\Big|\leq C\frac{\min\{1,\al^2|y-y_c|^2\}\phi_1}{u'(y_c)^2}.
\een
On the other hand, we have
\beno
\pa_c\phi_1=2\al^2T_0\circ T_{2,3}(\phi_1)-2\al^2T_0\circ T_{1,2}(\phi_1)+\al^2 T(\pa_c\phi_1),
\eeno
from which and \eqref{eq:phi-pc}, we infer that
\begin{align*}
|\pa_c\phi_1(y,c)|&\leq
C\al^2\Big|\int_{y_c}^y\f{\phi_1(y',c)}{u'(y')+u'(y_c)}\,dy'\Big|
+C\al^3\Big|\int_{y_c}^y|y'-y_c|\frac{\phi_1(y',c)}{u'(y_c)}\,dy'\Big|\\
&\leq C\phi_1(y,c)\Big(\f{\al^2|y-y_c|}{u'(y_c)}
+\frac{\al^3|y-y_c|^2}{u'(y_c)}\Big).
\end{align*}
Here we again used the fact that $\phi_1(y',c)\leq \phi_1(y,c)$  for $y_c\leq y'\leq y$ or $y\leq y'\leq y_c$.

{\bf Step 3.} Estimates of $\pa_c^2\phi_1$ and $\Big(\frac{\partial_y}{u'(y_c)}+\partial_c\Big)^2\phi_1$.

By Lemma \ref{lem:phi1}, we get for $|y-y_c|\leq \f{M_0}{\al}$
\begin{align*}
|\pa_c^2\cF(y,c)|&\leq \f{C|\pa_c^2\pa_y\phi_1|}{\phi_1}+\f{C|\pa_c\pa_y\phi_1\pa_c\phi_1|}{\phi_1^2}+\f{C|\pa_c^2\phi_1\pa_y\phi_1|}{\phi_1^2}+\f{C|\pa_c\phi_1|^2|\pa_y\phi_1|}{\phi_1^3}\\
&\leq \f{C\al^3}{u'(y_c)^2}+\f{C\al^2}{y_c(y+y_c)^2}.
\end{align*}
It follows from \eqref{eq:pcF} that
\begin{align*}
\partial_c^2\cF(y,c)&=\partial_y\partial_c\cG(y,c)\\
&=-2\int_{y_c}^y\partial_c\Big(\frac{\phi_1(z,c)^2}{\phi(y,c)^2}\cF(z,c)\Big)u'(z)dz\\
&=-2\int_{y_c}^y2\f{u'(z)\phi_1(z,c)}{\phi(y,c)}\Big(\f{\pa_c\phi_1(z,c)}{\phi(y,c)}-\f{\phi_1(z,c)\pa_c\phi(y,c)}{\phi(y,c)^2}\Big)\f{\pa_z\phi_1(z,c)}{\phi_1(z,c)}dz\\
&\quad-2\int_{y_c}^y\frac{\phi_1(z,c)^2}{\phi(y,c)^2}\partial_c\cF(z,c)u'(z)dz,
\end{align*}
which along with the fact $\int_{y_c}^y2\phi_1(z,c)\pa_z\phi_1(z,c)dz=\phi_1(y,c)^2-1$, \eqref{eq:phi-est1}, \eqref{eq:F-est1} and \eqref{eq: pa_cf}  implies that for $|y-y_c|\geq \f{M_0}{\al}$,
\begin{align*}
|\partial_c^2\cF(y,c)|\leq \frac{C\al}{|y-y_c|^2u'(y_c)^2}.
\end{align*}

Thanks to $\pa_c\cG(y_c,c)=\pa_c^2\phi_1(y_c,c)=\al^2\pa_c^2T\phi_1=\f{\al^2}{3u'(y_c)^2}$, we get
\beno
|\partial_c\cG(y,c)|\leq  \frac{C\al^2}{u'(y_c)^2},
\quad\text{thus,}\,\, |\partial_c^2\phi_1|\leq \frac{C\al^2\phi_1}{u'(y_c)^2}.
\eeno
Using the formula
$$\Big(\frac{\partial_y}{u'(y_c)}+\partial_c\Big)\big(T_0(a)\big)=\frac{a(y,c)}{u'(y_c)}-\frac{a(y_c,c)}{u'(y_c)}+T_0(\partial_ca)
=T_0\Big(\Big(\frac{\partial_y}{u'(y_c)}+\partial_c\Big)a\Big),$$
we obtain
\begin{align}
\Big(\frac{\partial_y}{u'(y_c)}+\partial_c\Big)\partial_y\cG_1(y,c)=&\Big(\Big(\frac{\partial_y}{u'(y_c)}+\partial_c\Big)\frac{-2}{\phi(y,c)^2}\Big)T_0
\Big(\phi_1^2\frac{a_1}{u'(y_c)}\cF\Big)\nonumber\\
&+\frac{-2}{\phi(y,c)^2}T_0\Big(\Big(\frac{\partial_y}{u'(y_c)}+\partial_c\Big)\Big(\phi_1^2\frac{a_1}{u'(y_c)}\cF\Big)\Big).\label{eq:G-pcy}
\end{align}

Using \eqref{eq: good_1phi_1} and the fact that
$$
\Big|\Big(\frac{\partial_y}{u'(y_c)}+\partial_c\Big)(u(y)-c)\Big|=\Big|\frac{u'(y)}{u'(y_c)}-1\Big|\leq C\frac{|y-y_c|}{u'(y_c)}\leq C\frac{|u(y)-c|}{u'(y_c)^2},
$$
we infer that
\ben\label{eq:phi-pcy-1}
\left|\Big(\frac{\partial_y}{u'(y_c)}+\partial_c\Big)\frac{-2}{\phi(y,c)^2}\right|\leq \frac{C}{\phi(y,c)^2u'(y_c)^2}.
\een
Using the formula
\begin{align*}
&\Big(\frac{\partial_y}{u'(y_c)}+\partial_c\Big)\Big(\phi_1^2\frac{a_1(y,c)}{u'(y_c)}\cF\Big)\\
&=2\phi_1\frac{a_1(y,c)}{u'(y_c)}\cF\Big(\frac{\partial_y}{u'(y_c)}+\partial_c\Big)\phi_1+\phi_1^2\frac{a_2(y,c)}{u'(y_c)^2}\cF\\
&\quad-\phi_1^2\frac{a_1(y,c)}{u'(y_c)^3}u''(y_c)\cF+\phi_1^2\frac{a_1(y,c)}{u'(y_c)}\Big(\frac{\partial_y}{u'(y_c)}+\partial_c\Big)\cF,
\end{align*}
we get by \eqref{eq:a1}, \eqref{eq:a2}, \eqref{eq: good_1phi_1} and \eqref{eq: good_1F}  that
\begin{align*}
&\left|\Big(\frac{\partial_y}{u'(y_c)}+\partial_c\Big)\Big(\phi_1^2\frac{a_1(y,c)}{u'(y_c)}\cF\Big)\right|\\
&\leq
C\phi_1^2\frac{|y-y_c|^2}{u'(y_c)^3}|\cF|+C\phi_1^2\frac{|y-y_c|^2}{u'(y_c)}
\frac{\phi_1(y,c)^2-1}{u'(y_c)\phi_1(y,c)^2}\frac{|y-y_c|^2}{(u(y)-c)^2}\\
&\leq C\phi_1^2\frac{|y-y_c|^2}{u'(y_c)^3}|\cF|+C|y-y_c|
\frac{\phi_1(y,c)^2-1}{u'(y_c)^3},
\end{align*}
which gives
\begin{align}
&\left|T_0\left(\Big(\frac{\partial_y}{u'(y_c)}+\partial_c\Big)\Big(\phi_1^2\frac{a_1}{u'(y_c)}\cF\Big)\right)\right|\nonumber\\
&\leq
C \frac{|y-y_c|^2}{u'(y_c)^3}|T_0(\phi_1^2|\cF|)|+C\Big|T_0\Big(|y-y_c|
\frac{\phi_1(y,c)^2-1}{u'(y_c)^3}\Big)\Big|\nonumber\\
&\leq C|y-y_c|^2
\frac{\phi_1(y,c)^2-1}{u'(y_c)^3}.\label{eq:T-F-pc}
\end{align}
Here we used
\beno
|T_0(\phi_1^2|\cF|)|\le \int_{y_c}^y\phi_1(y',c)\pa_y\phi_1(y',c)dy'\le \phi_1(y,c)^2-1.
\eeno

It follows from \eqref{eq:G-pcy}, \eqref{eq:phi-pcy-1}, \eqref{eq:pycF}, \eqref{eq: good_1F} and \eqref{eq:T-F-pc} that
\begin{align*}
&\Big|\Big(\frac{\partial_y}{u'(y_c)}+\partial_c\Big)\partial_y\cG_1(y,c)\Big|\\
&\leq C \frac{|\partial_y\cG_1(y,c)|}{u'(y_c)^2}+C\frac{|y-y_c|^2}{\phi(y,c)^2}
\frac{\phi_1(y,c)^2-1}{u'(y_c)^3}\\
&\leq C\frac{\phi_1(y,c)^2-1}{u'(y_c)^3\phi_1(y,c)^2}\frac{|y-y_c|^2}{(u(y)-c)^2}.
\end{align*}
This together with $\Big(\frac{\partial_y}{u'(y_c)}+\partial_c\Big)\cG_1(y_c,c)=\pa_c\cG_1(y_c,c)=0$ gives
\beno
\Big|\Big(\frac{\partial_y}{u'(y_c)}+\partial_c\Big)\cG_1(y,c)\Big|\leq C\frac{\phi_1(y,c)^2-1}{u'(y_c)^3\phi_1(y,c)^2}\Big|\int_{y_c}^y\frac{|y'-y_c|^2}{(u(y')-c)^2}dy'\Big|\leq \frac{C(\phi_1(y,c)^2-1)}{u'(y_c)^4\phi_1(y,c)^2},
\eeno
from which and \eqref{eq: good_1phi_1}, we infer that
\beno
\Big|\Big(\frac{\partial_y}{u'(y_c)}+\partial_c\Big)^2\phi_1(y,c)\Big|\leq C\frac{\min\{1,\al^2|y-y_c|^2\}\phi_1}{u'(y_c)^4}.
\eeno

This completes the proof of the proposition.
\end{proof}

It is easy to see from the proof of Proposition \ref{prop:phi1} that

\begin{lemma}\label{lem:FG}
It holds that
\begin{align*}
&C^{-1}\al\min\{\al y_c,1\}\leq |\cF(0,c)|\leq C\al\min\{\al y_c,1\},\\
&|\pa_c\cF(0,c)|\leq \f{C|\cF(0,c)|^2}{\al^2y_c^3},\quad
|\pa_c^2\cF(0,c)|\leq \f{C\al\min\{\al y_c,1\}}{y_c^4},\\
&|\cG(0,c)|\leq \f{C\al\min\{1,\al y_c\}}{y_c},\quad
|\pa_c\cG(0,c)|\leq \f{C\al^2}{y_c^2}.
\end{align*}
\end{lemma}

\section{Spectral analysis of the linearized operator}

In this section, we study the spectrum of the linearized operator $\cR_\al$
defined by
\beno
\mathcal{R}_\al\widehat{\psi}=-(\partial_y^2-\al^2)^{-1}\big(u''(y)-u(\partial_y^2-\al^2)\big)\widehat{\psi}.
\eeno

\subsection{Spectrum and embedding eigenvalues}
We consider the shear flows in $\mathcal{K}$.
Let us first recall some classical facts about the spectrum $\sigma(\mathcal{R}_\al)$ of the operator $\mathcal{R}_\al$ in $L^2(-1,1)$(see \cite{Ros, Ste, WZZ1} for more details):

\begin{itemize}

\item[1.] The spectrum $\sigma(\mathcal{R}_\al)$ is compact;

\item[2.] $\sigma_d(\mathcal{L})=\bigcup_{\al}\sigma_d(i\al\mathcal{R}_\al)$;

\item[3.] The continuous spectrum $\sigma_c(\mathcal{R}_\al)$ is contained in the range $ \textrm{Ran}\,u$ of $u(y)$;

\item[4.] The eigenvalues of $\mathcal{R}_\al$ can not cluster except possibly along on $\textrm{Ran}\,u$;


\end{itemize}

Usually, $c\in \text{Ran}\, u$ is called an embedding eigenvalue of $\mathcal{R}_\al$ if there exists
$0\neq \psi\in H_0^1(-1,1)$ so that
\ben\label{eq:R-eigen}
\mathcal{R}_\al\psi=c\psi\quad \text{or}\quad (u-c)(-\pa_y^2+\al^2)\psi+{u''\psi}=0.
\een
However, this definition is too general so that one can construct a nontrivial solution $\psi\in H_0^1(-1,1)$ to \eqref{eq:R-eigen}
for all $c\in \text{Ran}\, u$. Thus, we introduce the following definition.

\begin{definition}\label{embed}
We say that $c\in \mathrm{Ran}\, u$ is an embedding eigenvalue of $\mathcal{R}_\al$ if there exists $0\neq \psi\in H_0^1(-1,1)$ such that  for all $\varphi\in H_0^1(-1,1)$,
\begin{align}\label{u1}
\int_{-1}^1(\psi'\varphi'+\al^2\psi\varphi)dy+p.v.\int_{-1}^1\frac{u''\psi\varphi}{u-c}dy+i\pi\sum_{y\in u^{-1}\{c\},u'(y)\neq 0}\frac{(u''\psi\varphi)(y)}{|u'(y)|}=0.
\end{align}
\end{definition}

Let us first explain that the integration $p.v.\int_{-1}^1\frac{u''\psi\varphi}{u-c}dy $ is well defined.  It suffices to consider the point where $u'(y_0)=0$ and  $u(y_0)=c$. Without loss of generality, we assume $u'(y)>0$ in $(y_0,y_0+\d)$. Then in $(y_0,y_0+\d)$, $\psi(y)$ satisfies  \eqref{eq:R-eigen}, thus has the form
\beno
\psi(y)=C_1\phi(y)+C_2\phi(y)\int_{y_0+\d}^y\f{1}{\phi(y')^2}dy'
\eeno
for two constants $C_1, C_2$, where $\phi(y)=(u(y)-u(y_0))\phi_1(y)$ with $\phi_1(y)$ solving $\phi_1(y)=1+\al^2T(\phi_1)$ as in section 4.3.  Thus,
\beno
\phi(y)\int_{y_0+\d}^y\f{1}{\phi(y')^2}dy'\sim \f{1}{y-y_0},
\eeno
which is not in $L^2_{loc}$, thus $C_2=0$. Then $ \frac{\psi}{u-c}$ is bounded and the integral is well-defined.

It is easy to see from \eqref{u1} that
\begin{align}\label{u2}
-\psi''+\al^2\psi+\frac{u''\psi}{u-c}=0\quad \text{in}\quad [-1,1]\setminus u^{-1}\{c\}.
\end{align}
Moreover, for $u(y)$ satisfying the condition $(F)$: $u''(y_1)u''(y_2)\geq0$ if $u(y_1)=u(y_2)$, we can deduce from \eqref{u1} that $\psi\in H^2(-1,1)$ is a classical solution to \eqref{u2} see \cite{Lin1}. Obviously, the condition $(S)$ implies $(F)$.

\begin{lemma}
Let $u(y)$ satisfy the condition (F). If $c$ is an embedding eigenvalue of $\mathcal{R}_\al$,  then $u''(y)= 0 $ for some $y\in u^{-1}\{c\}\setminus\{\pm1\}$.
\end{lemma}

\begin{proof}
If $u''(y)\neq 0 $ for all $y\in u^{-1}\{c\}\setminus\{\pm1\}$, {taking $ \varphi=\overline{\psi}$ in} \eqref{u1}, we deduce that $\psi(y)=0$ for $u(y)=c$ and
$$\int_{-1}^1\left(|\psi'|^2+\al^2|\psi|^2+\frac{u''|\psi|^2}{u-c}\right)dy=0.$$
Thus, we can use integration by parts to obtain
\beno
\int_{-1}^1\Big|\psi'-u'\frac{\psi}{u-c}\Big|^2dy+\al^2\int_0^1|\psi|^2dy=0,
\eeno
which implies $\psi\equiv 0$.
\end{proof}

\subsection{A criterion on embedding eigenvalues}
Here we give an equivalent characterization of embedding eigenvalues for symmetric flow $u\in C^4([-1,1])$ satisfying $(S)$, which is key to solve
the inhomogeneous Rayleigh equation.\smallskip

For $c\in D_0$, let $\phi(y,c)$ be the solution of the homogeneous Rayleigh equation given by Lemma \ref{lem:Ray-D0}  and $\phi_1(y,c)=\frac {\phi(y,c)} {u(y)-c}$. We introduce
\ben
&&\mathrm{A}(c)=\mathrm{A}_1(c)+u'(y_c)\rho(c)\mathrm{II}_3(c),\quad\mathrm{B}(c)=\pi\rho(c)\frac{u''(y_c)}{u'(y_c)^2},\label{def:A}\\
&&\rmA_2(c)=(u(0)-c)\rmA(c)+J(c),\quad \rmB_2(c)=(u(0)-c)\rmB(c),\label{def:B}
\een
where
\ben
&&\rho(c)=(c-u(0))(u(1)-c),\quad J(c)=\f{u'(y_c)(u(1)-c)}{\phi_1(0,c)\phi_1'(0,c)},\label{def:rho-J}\\
&&\rmA_1(c)=\rho u'(y_c)\pa_c\Big(p.v.\int_0^1\f{dy}{u(y)-c}\Big),\label{def:A1}\\
&&\mathrm{II}_3(c)=\int_0^1\frac{1}{(u(y)-c)^2}\Big(\frac{1}{\phi_{1}(y,c)^2}-1\Big)dy.\label{eq:Pi-3}
\een

By Lemma \ref{lem:phi1} and Proposition \ref{prop:Rayleigh-Hom}, $J(c)$ and $\mathrm{II}_3(c)$ are well-defined.

\begin{proposition}\label{prop:spectral}
$c\in D_0$ is an embedding eigenvalue of $\mathcal{R}_\al$ if and only if
\beno
\mathrm{A}(c)^2+\mathrm{B}(c)^2=0\quad\textrm{or}\quad \rmA_2(c)^2+\rmB_2(c)^2=0.
\eeno
\end{proposition}

\begin{proof}
Let us first assume that $c\in D_0$ is an embedding eigenvalue of $\mathcal{R}_{\al}$. Thus, $u''(y_c)=0$(thus, $\rmB=\rmB_2=0$) and there exists a nontrivial solution $\varphi$ of
\beq\label{eq:embedding eigenvalue_S}
\left\{
\begin{array}{l}
\varphi''-\al^2\varphi-\frac{u''}{u-c}\varphi=0,\\
\varphi(-1)=\varphi(1)=0.
\end{array}\right.
\eeq
Let
\beno
\phi_o(y,c)=\f{\varphi(y,c)-\varphi(-y,c)}{2},\quad  \phi_e(y,c)=\f{\varphi(y,c)+\varphi(-y,c)}{2}.
\eeno
Then $\phi_0$ and $\phi_e$ also satisfies the homogeneous Rayleigh equation in [0,1] with
\beno
\phi_o(0,c)=\phi_o(1,c)=0,\quad \pa_y\phi_e(0,c)=\phi_e(1,c)=0.
\eeno

If $\phi_o$ is a nontrivial solution, then $\pa_y\phi_o(0,c)\neq 0$. Otherwise, it is trivial by the fact that  $u'(y)>0$ for $y\in (0,1)$ and $u''(y_c)=0$. Without loss of generality, we assume $\pa_y\phi_{{o}}(0,c)=1$. Thus,
$\phi_o(y,c)$ has the following representation formula
\beno
\phi_o(y,c)=\phi(0,c)\overline{\varphi}(y,c),
\eeno
where $\overline{\varphi}(y,c)$ is given by
\begin{align*}
\overline{\varphi}(y,c)=&\frac{\phi_1(y,c)}{u'(y_c)}(u(y)-u(0))\\
&+\frac{\phi_1(y,c)}{u'(y_c)}(u(y)-c)(u(0)-c)\int_0^y\frac{u'(y_c)-u'(z)}{(u(z)-c)^2}dz\\
&+\phi_1(y,c)(u(y)-c)(u(0)-c)\int_0^y\frac{1}{(u(z)-c)^2}\Big(\frac{1}{\phi_1(z,c)^2}-1\Big)dz.
\end{align*}
In fact, the solution $\phi_{{o}}(y,c)$ can be formally written as
\beno
\phi_{{o}}(y,c)={\phi(0,c)}\phi(y,c)\int_0^y\f 1 {\phi(y',c)^2}dy'.
\eeno
Then we find that $\phi_o(1,c)=0$ is equivalent to
\begin{align*}
&0=\phi_o(1,c)={\phi(0,c)\overline{\varphi}(1,c)}\\
=&\frac{\phi_1(0,c)\phi_1(1,c)}{u'(y_c)}(u(1)-u(0))
-\frac{\phi_1(0,c)\phi_1(1,c)}{u'(y_c)}\rho(c)\int_0^1\frac{u'(y_c)-u'(y)}{(u(y)-c)^2}dy\\
&-\phi_1(0,c)\phi_1(1,c)\rho(c)\int_0^1\frac{1}{(u(y)-c)^2}\Big(\frac{1}{\phi_1(y,c)^2}-1\Big)dy\\
=&-\frac{\phi_1(0,c)\phi_1(1,c)}{u'(y_c)}\mathrm{A}(c).
\end{align*}
Thus, $\rmA(c)=0$. Here we used
\ben\label{def:A1-Pi2}
\rmA_1(c)=u(0)-u(1)-\rho(c)\mathrm{II}_2(c),\quad \mathrm{II}_2(c)=p.v.\int_0^1\f {u'(y)-u'(y_c)} {(u(y)-c)^2}dy.
\een

If $\phi_e$ is a nontrivial solution, we may assume that
$\pa_y\phi_e(1,c)=1$. Thus, ${\phi}_e(y,c)$ has the following representation formula
\beno
{\phi}_e(y,c)={\phi(1,c)}{\widetilde{\varphi}}(y,c),
\eeno
where ${\widetilde{\varphi}}(y,c)$ is given by
\begin{align*}
{\widetilde{\varphi}}(y,c)=&\frac{\phi_1(y,c)}{u'(y_c)}(u(y)-u(1))\\
&+\frac{\phi_1(y,c)}{u'(y_c)}(u(y)-c)(u(1)-c)\int_1^y\frac{u'(y_c)-u'(z)}{(u(z)-c)^2}dz\\
&+\phi_1(y,c)(u(y)-c)(u(1)-c)\int_1^y\frac{1}{(u(z)-c)^2}\Big(\frac{1}{\phi_1(z,c)^2}-1\Big)dz.
\end{align*}
Then we find that
\begin{align*}
0
=&\pa_y{\widetilde{\varphi}}(0,c)=\frac{\pa_y\phi_1(0,c)}{u'(y_c)}(u(0)-u(1))\\
&+\frac{\pa_y\phi_1(0,c)}{u'(y_c)}(u(0)-c)(u(1)-c)\int_1^0\frac{u'(y_c)-u'(z)}{(u(z)-c)^2}dz\\
&+(u(1)-c)\frac{1}{u(0)-c}\frac{1}{\phi_1(0,c)}\\
&+\pa_y\phi_1(0,c)(u(0)-c)(u(1)-c)\int_1^0\frac{1}{(u(z)-c)^2}\Big(\frac{1}{\phi_1(z,c)^2}-1\Big)dz\\
=&\f {\pa_y\phi_1(0,c)} {(u(0)-c)u'(y_c)}\rmA_2(c).
\end{align*}
Thus, $\rmA_2(c)=0$.

It remains to show that if $\mathrm{A}(c)^2+\mathrm{B}(c)^2=0$ or $\rmA_2(c)^2+\rmB_2(c)^2=0$, then $c\in D_0$ must be an embedding eigenvalue of $\mathcal{R}_\al$. In this case, we have $u''(y_c)=0$. Thus, we can construct a nontrivial solution of the homogeneous Rayleigh equation with $\varphi(0,c)=0, \pa_y\varphi(0,c)=1$(or $\varphi(1,c)=0, \pa_y\varphi(1,c)=1$).
The above argument shows that $\varphi(1,c)=0$(or $\varphi(0,c)=0$). This means that $c$ is an embedding eigenvalue of
$\mathcal{R}_\al$ with the eigenfunction $\varphi(y,c)$.
\end{proof}

\section{The limiting absorption principle}

In this section, we establish the limiting absorption principle for the inhomogeneous Rayleigh equation when $u\in \mathcal{K}$:
\begin{align}\label{u3}
(u-c)(\Phi''-\al^2\Phi)-{u''\Phi}={\omega},\quad \Phi(-1)=\Phi(1)=0,
\end{align}
where $c\in \Om_{\epsilon_0}\setminus D_0$.

The key ingredient is the following uniform estimate of the solution.

\begin{proposition}\label{prop:ub}
If $\mathcal{R}_\al$ has no embedding eigenvalues, then there exists $\epsilon_0$ such that for $c\in \Om_{\epsilon_0}\setminus D_0$, the solution to \eqref{u3} has the the following uniform bound
\beno
\|\Phi\|_{H^1(-1,1)}\leq C\|\omega\|_{H^1(-1,1)}.
\eeno
Here $C$ is a constant independent of $\epsilon_0$.
\end{proposition}

The proof is based on blow-up analysis and compactness argument. 

\begin{lemma}\label{lem:ub1}
Let $\psi_n, \ \omega_n\in H^1(a,b),\ u_n\in H^3(a,b)$ be a sequence, which satisfies
\beno
&&\psi_n\rightharpoonup\psi,\ \omega_n\to\omega\quad \text{in}\quad H^1(a,b),\\
&&u_n\to u_0\quad \text{in}\quad H^3(a,b),\\
&&u_n(\psi_n''-\al_n^2\psi_n)-{u_n''\psi_n}={\omega_n},
\eeno
and $\al_n\to\al,\ \mathrm{Im}\ u_n<0$. Moreover, $\mathrm{Im}\, u_0=0 $ in $[a,b]$, $u(y_0)=0,\ y_0\in[a,b], \ u'(y)u'(y_0)>0$ in $[a,b]$.
Then we have $\psi_n\to\psi$ in $H^1(a,b)$, and for all $\varphi\in H_0^1(a,b)$,
\begin{align}\label{u4}
\int_{a}^b(\psi'\varphi'+\al^2\psi\varphi)dy+p.v.\int_{a}^b\frac{(u_0''\psi+\omega)\varphi}{u_0}dy
+i\pi\frac{((u_0''\psi+\omega)\varphi)(y_0)}{|u'(y_0)|}=0.
\end{align}
\end{lemma}

\begin{proof}
Without loss of generality, we may assume that $u_0'(y)>0$ for $y\in [a,b]$. Otherwise, we can consider $\overline{\psi}_n, \overline{\omega}_n,-\overline{u}_n$.

As $\psi_n\rightharpoonup\psi,\ \omega_n\to\omega$ in $H^1(a,b)$ and $u_n\to u_0$ in $H^3(a,b)$,  $\psi_n, \omega_n, u_n'' $ are uniformly bounded in $H^1(a,b)$. Let $ g_n={u_n''\psi_n}+{\omega_n}$. Then $g_n$ is uniformly bounded in $H^1(a,b)$ and $\psi_n''-\al_n^2\psi_n=g_n/u_n.$ Choose $N$ large enough so that $\mathrm{Re}\,u_n'(y)\geq c_0>0$ for $n\geq N$ and $y\in [a,b]$.

It is easy to see that
\beno
\left(\psi_n'-\dfrac{g_n}{u_n'}\ln u_n\right)'=\al_n^2\psi_n-\left(\dfrac{g_n}{u_n'}\right)'\ln u_n.
\eeno
Here $\ln (re^{i\theta})=\ln r+i\theta$ for $r>0,\theta\in[-\pi,\pi]$. We know that $\left(\dfrac{g_n}{u_n'}\right)'$ is uniformly bounded in $L^2(a,b)$, and $\dfrac{g_n}{u_n'} $ is uniformly bounded in $L^{\infty}(a,b).$ Due to
$\mathrm{Re}\,u_n'(y)\geq c_0>0$ for $y\in [a,b]$,  $\ln u_n$ is uniformly bounded in $L^p(a,b), p>2$.
Therefore,  $\psi_n'-\dfrac{g_n}{u_n'}\ln u_n $ is uniformly bounded in $L^2(a,b)$ and $ \dot{W}^{1,1}(a,b)$, thus bounded in $L^{\infty}(a,b).$ Thus, we conclude
\begin{align}\label{u5}
\lim_{\d\to0}\sup_n\|\psi_n'\|_{L^2((y_0-\d,y_0+\d)\cap (a,b))}=0.
\end{align}

On the other hand,  $\psi_n'' $ is bounded in $L^{\infty}((a,b)\setminus(y_0-\d,y_0+\d))$. So, $\psi_n\to\psi$ in $H^1((a,b)\setminus(y_0-\d,y_0+\d))$, which together with \eqref{u5} shows that $\psi_n\to\psi$ in $H^1(a,b)$ and $g_n\to g_0={u_0''\psi}+{\omega}$ in $H^1(a,b)$.

Now we prove \eqref{u4}. As $\psi_n''-\al_n^2\psi_n=g_n/u_n,$ we have
$$\int_{a}^b(\psi_n'\varphi'+\al_n^2\psi_n\varphi)dy+\int_{a}^b\frac{g_n\varphi}{u_n}dy=0$$
for any $\varphi\in H_0^1(a,b)$. Obviously,
\beno
\lim\limits_{n\to\infty}\int_{a}^b(\psi_n'\varphi'+\al_n^2\psi_n\varphi)dy=\int_{a}^b(\psi'\varphi'+\al^2\psi\varphi)dy.
\eeno
It remains to show that
\begin{align*}
\lim_{n\to\infty}\int_{a}^b\frac{g_n\varphi}{u_n}dy=p.v.\int_{a}^b\frac{g_0\varphi}{u_0}dy
+i\pi\frac{(g_0\varphi)(y_0)}{|u_0'(y_0)|}.
\end{align*}
Using the facts that
\beno
 &&\int_{a}^b\frac{g_n\varphi}{u_n}dy=-\int_{a}^b\big(\frac{g_n\varphi}{u_n'}\big)'\ln u_ndy,\quad \big(\frac{g_n\varphi}{u_n'}\big)'\to\big(\frac{g_0\varphi}{u_0'}\big)',\\
 &&\ln u_n\to\ln |u_0|-i\pi \chi_{[a,y_0)}\quad \text{in}\quad L^2(a,b),
 \eeno
we infer that
\begin{align*}
\int_{a}^b\frac{g_n\varphi}{u_n}dy&\to-\int_{a}^b\big(\frac{g_0\varphi}{u_0'}\big)'\ln |u_0|dy+i\pi\int_{a}^{y_0}\big(\frac{g_0\varphi}{u_0'}\big)'dy\\
&\quad=p.v.\int_{a}^b\frac{g_0\varphi}{u_0}dy
+i\pi\frac{(g_0\varphi)(y_0)}{|u_0'(y_0)|}.
\end{align*}
This proves \eqref{u4}.
\end{proof}

\begin{lemma}\label{lem:ub3}
Assume that $\psi, \ \omega\in H^1(a,b),\ u\in H^3(a,b)$ satsify
\beno
(u-c)(\psi''-\al^2\psi)-{u''\psi}={\omega},
\eeno
where $c\not\in \R$, and $u$ is real-valued and $u'(y_0)=0,\ y_0=\f{a+b}{2}\in(a,b),\ u''(y)u''(y_0)>0$ in $[a,b]$.
Then we have
\begin{align}
|(\psi''-\al^2\psi)(y_0)|\leq C\min(|u(y_0)-c|^{-\frac{3}{4}},|u(y_0)-c|^{-1})\big(\|\psi\|_{H^1(a,b)}+\|\omega\|_{H^1(a,b)}\big),\nonumber
\end{align}
where the constant $C$ depends only on $u, \al$ and $\d=y_0-a.$
 \end{lemma}
\begin{proof}
We may assume that $|u''|>c_0>0$ in $[a,b]$. Then we have
\beno
&&|u'(y)|\leq C|y-y_0|,\quad |u(y)-u(y_0)|\leq C|y-y_0|^2,\\
&&\ |u'(y)-u'(y_1)|\geq c_0|y-y_1|\quad \text{for}\,\, y,y_1\in [a,b].
\eeno
It is easy to see that
 $$|(\psi''-\al^2\psi)(y_0)|=\left|\frac{{u''\psi}-{\omega}}{u-c}\right|(y_0)\leq C\frac{\|\psi\|_{H^1(a,b)}+\|\omega\|_{H^1(a,b)}}{|u(y_0)-c|}.$$
Thus, it suffices to  consider the case $|u(y_0)-c|<\d^2$. Let $\d_1=|u(y_0)-c|^{\f12}$. So, $\d_1<\d.$

First of all, we consider the case of $\omega(y_0)=0$. We normalize $\|\psi\|_{H^1(a,b)}+\|\omega\|_{H^1(a,b)}=1. $ Then for $|y-y_0|<\d_1,$
\beno
|\omega(y)|\leq C\d_1^{\f12},\quad
|u(y)-c|\leq |u(y)-u(y_0)|+|u(y_0)-c|\leq C\d_1^2
\eeno
Let $g=((u-c)\psi'-u'\psi)'=\al^2\psi(u-c)+\omega$. Then for $|y-y_0|<\d_1$,
\beno
|g(y)|\leq C|u(y)-c|+|\omega(y)|\leq C\d_1^{\f12},
\eeno
which gives
\beno
\left|((u-c)\psi'-u'\psi)\big|_{y_1}^{y_2}\right|=\left|\int_{y_1}^{y_2}g(y)dy\right|\leq C\d_1^{\f32}.
\eeno

Choose $\d_2\in(\d_1/2,\d_1),\ y_1=y_0-\d_2,\ y_2=y_0+\d_2$ so that
 \beno
|\psi'(y_1)|^2+|\psi'(y_2)|^2\leq \d_1^{-1}\|\psi'\|_{L^2(a,b)}^2\leq \d_1^{-1},
\eeno
which gives
\beno
\left|(u-c)\psi'\big|_{y_1}^{y_2}\right|\leq(|\psi'(y_1)|+|\psi'(y_2)|)\|u-c\|_{L^\infty([y_1,y_2])}\leq C\d_1^{-\f12}\d_1^2=C\d_1^{\f32}.
\eeno
This shows that $\left|u'\psi\big|_{y_1}^{y_2}\right|\leq C\d_1^{\f32}. $ Notice that
\beno
u'\psi\big|_{y_1}^{y_2}=\psi(y_0)u'\big|_{y_1}^{y_2}+u'(y_1)\psi\big|_{y_1}^{y_0}+u'(y_2)\psi\big|_{y_0}^{y_2}.
\eeno
Thus, we have
\beno
c_0\d_1|\psi(y_0)|\leq \left|\psi u'\big|_{y_1}^{y_2}\right|+C\|u'\|_{L^\infty([y_1,y_2])}\d_2^{\f12}\|\psi'\|_{L^2(a,b)}\leq C\d_1^{\f32}+C\d_1\d_1^{\f12}=C\d_1^{\f32},
\eeno
which gives
\ben
|\psi(y_0)|\leq C\d_1^{\f12}\leq C|u(y_0)-c|^{\f12},\label{eq:psi-y0}
\een
and thus,
$$
|(\psi''-\al^2\psi)(y_0)|=\left|\frac{{u''\psi}(y_0)}{u(y_0)-c}\right|\leq C\d_1^{-2}|\psi(y_0)|\leq C\d_1^{-\f32}.
$$

For the general case, we introduce
\beno
\psi_*(y)=\psi(y)-\frac{\omega(y_0)}{u''(y_0)}\cosh\al(y-y_0),\quad
\omega_*(y)=\omega(y)-u''(y)\frac{\omega(y_0)}{u''(y_0)}\cosh\al(y-y_0).
\eeno
We find that $\psi_*, \omega_*\in H^1(a,b),\, \omega_*(y_0)=0$, and
\beno
(u-c)(\psi_*''-\al^2\psi_*)-{u''\psi_*}={\omega_*}.
\eeno
Then the above argument shows that
\begin{align*}
|(\psi''-\al^2\psi)(y_0)|&=|(\psi_*''-\al^2\psi_*)(y_0)|\\
&\leq C\d_1^{-\f32}(\|\psi_*\|_{H^1(a,b)}+\|\omega_*\|_{H^1(a,b)})\\
&\leq C\d_1^{-\f32}(\|\psi\|_{H^1(a,b)}+\|\omega\|_{H^1(a,b)}+|\omega(y_0)|)\\
&\leq C\d_1^{-\f32}(\|\psi\|_{H^1(a,b)}+\|\omega\|_{H^1(a,b)}).
\end{align*}
This gives our result by recalling $\d_1=|u(y_0)-c|^{\f12}$.
\end{proof}

\begin{lemma}\label{lem:ub2}
Let $\psi_n, \ \omega_n\in H^1(a,b),\ u\in H^3(a,b)$ be a sequence,
which satisfies
\beno
&&\psi_n\rightharpoonup0,\,\, \omega_n\to0\quad \text{in}\,\, H^1(a,b),\\
&&(u-c_n)(\psi_n''-\al^2\psi_n)-{u''\psi_n}={\omega_n},
\eeno
and $\mathrm{Im}\, c_n>0,\ c_n\to u(y_0)$, $u'(y_0)=0,\ y_0=\f{a+b}{2}\in(a,b),\  \d=y_0-a\in(0,1), \ u''(y)u''(y_0)>0$ in $[a,b]$. Then we have $\psi_n\to0$ in $H^1(a,b)$.
\end{lemma}

\begin{proof}
Without loss of generality, we may assume that $y_0=0,u''(y_0)=2,u(y_0)=0$. Then $[a,b]=[-\d,\d]$. Let $c_n=r_n^2e^{2i\theta_n},\ r_n>0,\ 0<\theta_n<\pi/2$. So, $r_n\to0.$

First of  all, we consider the case of $\omega_n(0)=0.$ By Lemma \ref{lem:ub3} and \eqref{eq:psi-y0}, we have
\beno
|(\psi_n''-\al^2\psi_n)(0)|\leq Cr_n^{-\f32},\quad |\psi_n(0)|\leq Cr_n^{\f12},\quad |(u''\psi_n+\omega_n)(0)|\leq Cr_n^{\f12}.
\eeno
We introduce
\beno
\widetilde{\psi}_n(y)=r_n^{-\f12}{\psi_n}(r_ny),\quad \widetilde{\omega}_n(y)=r_n^{-\f12}{\omega_n}(r_ny),\quad u_n(y)=r_n^{-2}(u(r_ny)-u(0)).
\eeno
It holds that
\beno
(u_n-e^{2i\theta_n})(\widetilde{\psi}_n''-(\al r_n)^2\widetilde{\psi}_n)-{u_n''\widetilde{\psi}_n}={\widetilde{\omega}_n},
\eeno
and
\beno
&&|\widetilde{\psi}_n(0)|=|r_n^{-\f12}{\psi}_n(0)|\leq C,\quad
\|\widetilde{\psi}_n'\|_{L^2(a/r_n,b/r_n)}=\|{\psi}_n'\|_{L^2(a,b)}\leq C,\\
&&|\widetilde{\omega}_n(0)|=0,\quad \|\widetilde{\omega}_n'\|_{L^2(a/r_n,b/r_n)}=\|{\omega}_n'\|_{L^2(a,b)}\to 0.
\eeno
Thus, $\widetilde{\psi}_n $ is bounded in $H_{loc}^1(\R),$ and $\widetilde{\omega}_n\to0 $ in $H_{loc}^1(\R).$ Up to a subsequence, we may assume that $\widetilde{\psi}_n\rightharpoonup \widetilde{\psi}_0$ in $H_{loc}^1(\R),$ $\theta_n\to \theta_0 $ with $\widetilde{\psi}_0'\in L^2(\R),\ \theta_0\in[0,\pi/2]. $

Using the facts that
\beno
u_n(y)=y^2\int_0^1\int_0^1tu''(r_nyts)dsdt,\quad
u_n'(y)=y\int_0^1u''(r_nyt)dt,
\eeno
and $u_n''(y)=u''(r_n y), \ u_n'''(y)=r_nu'''(r_n y)$, we can deduce that $u_n\to y^2$ in $H_{loc}^3(\R)$.

Now, if $\theta_0\neq 0 $,  then $\widetilde{\psi}_n'' $ is bounded in $L_{loc}^{\infty}(\R)$ and $\widetilde{\psi}_n\to \widetilde{\psi}_0$ in $C_{loc}^1(\R),$ and
\beno
(y^2-e^{2i\theta_0})\widetilde{\psi}_0''={2\widetilde{\psi}_0}.
\eeno
If $\theta_0= 0$, then $\widetilde{\psi}_n'' $ is bounded in $L_{loc}^{\infty}(\R\setminus\{\pm1\})$ and $\widetilde{\psi}_n\to \widetilde{\psi}_0$ in $C_{loc}^1(\R\setminus\{\pm1\}),$ and
\beno
 (y^2-1)\widetilde{\psi}_0''={2\widetilde{\psi}_0}\quad\text{in}\,\,\R\setminus\{\pm1\}.
 \eeno
Lemma \ref{lem:ub1}  also ensures that  $\widetilde{\psi}_n\to \widetilde{\psi}_0 $ in $H^1(1-\d,1+\d)\cap H^1(-1-\d,-1+\d)$, thus,
\beno
\widetilde{\psi}_n\to \widetilde{\psi}_0 \quad  \textrm{in}\  H_{loc}^1(\R)\cap C_{loc}^1(\R\setminus\{\pm1\}).
\eeno

Thanks to $(y^2-e^{2i\theta_0})\widetilde{\psi}_0''={2\widetilde{\psi}_0},$ we deduce that for $y>1$,
\beno
\widetilde{\psi}_0(y)=(y^2-e^{2i\theta_0})\left(C_1+C_2\int_y^{+\infty}\frac{dz}{(z^2-e^{2i\theta_0})^2}\right),
\eeno
where $C_1,C_2$ are constants. As $\widetilde{\psi}_0'\in L^2(\R) $,  we have $C_1=0$, thus $\widetilde{\psi}_0\in L^2(2,+\infty). $ Similarly, $\widetilde{\psi}_0\in L^2(-\infty,-2)$. Hence,
$\widetilde{\psi}_0\in L^2(\R),\ \widetilde{\psi}_0\in H^1(\R).$

If $\theta_0\neq 0$, then
$$\int_{\R}|\widetilde{\psi}_0'|^2dy=-\int_{\R}\widetilde{\psi}_0''\overline{\widetilde{\psi}_0}dy
=-\int_{\R}\frac{2|\widetilde{\psi}_0|^2}{y^2-e^{2i\theta_0}}dy.$$
Multiplying $e^{i\theta_0}$ on both sides, and then taking the imaginary part,
we obtain
$$\sin \theta_0\int_{\R}|\widetilde{\psi}_0'|^2dy
=-\sin\theta_0\int_{\R}\frac{2(y^2+1)|\widetilde{\psi}_0|^2}{|y^2-e^{2i\theta_0}|^2}dy,$$
which implies $\widetilde{\psi}_0=0. $

If $\theta_0=0 $,  we first claim that for any $\varphi\in H^1(\R)$,
\begin{align}\label{u10}
\int_{\R}\widetilde{\psi}_0'\varphi'dy+p.v.\int_{\R}\frac{2\widetilde{\psi}_0\varphi}{y^2-1}dy
+i\pi\sum_{y=\pm1}\widetilde{\psi}_0\varphi=0.
\end{align}
Indeed, it holds that for $y\in \R\setminus\{\pm1\}$,  $(y^2-1)\widetilde{\psi}_0''={2\widetilde{\psi}_0}$, thus \eqref{u10} holds for $\varphi\in H^1(\R),\ \text{supp}\ \varphi\subset\R\setminus\{\pm1\}$. Lemma \ref{lem:ub1} ensures that \eqref{u10} holds for $\varphi\in H^1(\R),\ \text{supp}\ \varphi\subset[1-\d,1+\d]$ or $[-1-\d,-1+\d].$ Therefore, \eqref{u10} holds for any $\varphi\in H^1(\R). $

Taking $\varphi=\overline{\widetilde{\psi}_0}$ in \eqref{u10} and taking the imaginary part, we deduce that $\widetilde{\psi}_0(\pm1)=0, $ which implies $ \widetilde{\psi}_0\in H^2(\R).$  As $((y^2-1)\widetilde{\psi}_0'-{2y\widetilde{\psi}_0})'=0 $ and $((y^2-1)\widetilde{\psi}_0'-{2y\widetilde{\psi}_0})|_{\pm1}=0$, we infer that
\beno
(y^2-1)\widetilde{\psi}_0'-{2y\widetilde{\psi}_0}=0,\quad((y^2-1)^{-1}\widetilde{\psi}_0)'=0\quad \text{in}\,\, \R\setminus\{\pm1\},
\eeno
which mean that  there exist constants $C_3,C_4,C_5$ so that $\widetilde{\psi}_0=C_3(y^2-1)$ in $(-\infty,-1)$, $\widetilde{\psi}_0=C_4(y^2-1)$ in $(-1,1)$, $\widetilde{\psi}_0=C_5(y^2-1)$ in $(1,+\infty)$. However, $\widetilde{\psi}_0\in H^2(\R)$, thus  $C_3=C_4=C_5=0$.
Then $\widetilde{\psi}_0=0.$

Up to now,  we proved that $\widetilde{\psi}_n\to 0$ in $H_{loc}^1(\R)\cap C_{loc}^1(\R\setminus\{\pm1\}),$ thus $\|{\psi}_n'\|_{L^2(-2r_n,2r_n)}=\|\widetilde{\psi}_n'\|_{L^2(-2,2)}\to0. $ Since ${\psi}_n'' $ is bounded in $L_{loc}^{\infty}((a,b)\setminus\{0\})$ and $\psi_n\rightharpoonup0$ in $H^1(a,b)$, we have ${\psi}_n\to0 $ in $C_{loc}^{1}((a,b)\setminus\{0\}).$ Integration by parts gives
\begin{align}\label{u11}
\int_{2r_n}^b\left(|{\psi}_n'|^2+\al^2|{\psi}_n|^2+\frac{u''|{\psi}_n|^2}{u-c_n}\right)dy
=-\int_{2r_n}^b\frac{u''{\omega}_n\overline{{\psi}_n}}{u-c_n}dy+{\psi}_n'\overline{{\psi}_n}\big|_{2r_n}^b.
\end{align}
Notice that for $n$ sufficiently large and $y\in[2r_n,b]$, we have $u(y)\geq u(2r_n)>2r_n^2=2|c_n|$. Then we have
\beno
&&\text{Re}\dfrac{1}{u(y)-c_n}\geq \dfrac{1}{2(u(y)+|c_n|)}\geq \dfrac{1}{4u(y)}\geq \dfrac{1}{Cy^2},\\
&&\ u''(y)\geq C^{-1},\quad \left|\dfrac{1}{u(y)-c_n}\right|\leq \dfrac{1}{u(y)-|c_n|}\leq \dfrac{2}{u(y)}\leq \dfrac{C}{y^2}.
\eeno
 Taking the real part of \eqref{u11},  we get
\beno
\int_{2r_n}^b(|{\psi}_n'|^2+\al^2|{\psi}_n|^2)dy+C^{-1}\left\|\frac{{\psi}_n}{y}\right\|_{L^2(2r_n,b)}^2
\leq C\left\|\frac{{\psi}_n}{y}\right\|_{L^2(2r_n,b)}\left\|\frac{{\omega}_n}{y}\right\|_{L^2(2r_n,b)}
+\left|{\psi}_n'\overline{{\psi}_n}\big|_{2r_n}^b\right|,
\eeno
which implies that
\beno
\|{\psi}_n'\|_{L^2(2r_n,b)}^2
\leq C\left\|{{\omega}_n}/{y}\right\|_{L^2(2r_n,b)}^2
+\left|{\psi}_n'\overline{{\psi}_n}\big|_{2r_n}^b\right|.
\eeno
Notice that ${\psi}_n'\overline{{\psi}_n}\big|_{2r_n}^b={\psi}_n'\overline{{\psi}_n}(b)-\widetilde{\psi}_n'\overline{\widetilde{\psi}_n}(2)\to0$, and by Hardy's inequality, $\left\|{{\omega}_n}/{y}\right\|_{L^2(2r_n,b)}\leq C\left\|{{\omega}_n}\right\|_{H^1(a,b)}\to0 $. Thus, $\|{\psi}_n'\|_{L^2(2r_n,b)}\to0. $ Similarly, we have $\|{\psi}_n'\|_{L^2(a,-2r_n)}\to0$. This shows that $\|{\psi}_n'\|_{L^2(a,b)}\to0. $ As $\psi_n\rightharpoonup0$ in $H^1(a,b)$, we have $\|{\psi}_n\|_{L^2(a,b)}\to0 $ and then $\|{\psi}_n\|_{H^1(a,b)}\to0. $

For general case, we consider
\beno
&&\psi_{n*}(y)=\psi_n(y)-\frac{\omega_n(y_0)}{u''(y_0)}\cosh\al(y-y_0),\\\quad
&&\omega_{n*}(y)=\omega_n(y)-u''(y)\frac{\omega_n(y_0)}{u''(y_0)}\cosh\al(y-y_0).
\eeno
Then $\psi_{n*},\,\omega_{n*}\in H^1(a,b)$, $\omega_{n*}(y_0)=0$ and
\beno
 (u-c_n)(\psi_{n*}''-\al^2\psi_{n*})-{u''\psi_{n*}}={\omega_{n*}}.
\eeno
Notice that
\beno
&&\|\omega_{n*}\|_{H^1(a,b)}\leq\|\omega_{n}\|_{H^1(a,b)}+C|\omega_n(y_0)|\leq C\|\omega_{n}\|_{H^1(a,b)}\to 0,\\
&&\|\psi_{n}-\psi_{n*}\|_{H^1(a,b)}\leq C|\omega_n(y_0)|\leq C\|\omega_{n}\|_{H^1(a,b)}\to 0,
\eeno
thus, $\psi_{n*}\rightharpoonup0$ in $H^1(a,b)$. This is reduced to
 the case of $\om(y_0)=0$. Thus, $\|\psi_{n*}\|_{H^1(a,b)}\to0 $, hence $\|\psi_{n}\|_{H^1(a,b)}\to0$.
 \end{proof}

Now we are in a position to prove Proposition \ref{prop:ub}.
\begin{proof}
We only consider $\mathrm{Im}\, c>0$, and the case of  $\mathrm{Im}\, c<0$ can be proved by taking conjugation. We use the contradiction argument.

Assume  that there exists $\psi_n\in H_0^1(-1,1)$, $\omega_n \in H^1(-1,1),$ and $c_n$ with $\mathrm{Im}\, c_n>0$, such that $\|\psi_n\|_{H^1(-1,1)}=1,\ \|\omega_n\|_{H^1(-1,1)}\to 0,\ c_n\to c \in \text{Ran}\, u,$ and $$(u-c_n)(\psi_n''-\al^2\psi_n)-{u''\psi_n}={\omega_n}.$$ Then there exists a subsequence of $\{\psi_n\}$ (still denoted by $\{\psi_n\}$) and $\psi\in H_0^1(-1,1) $ such that $\psi_n\rightharpoonup\psi $ weakly in $H^1(-1,1)$. Since $\psi_n'' $ is bounded in $L_{loc}^{\infty}((-1,1)\setminus u^{-1}\{c\}),$ we have $\psi_n\to\psi $ in $H^1_{loc}((-1,1)\setminus u^{-1}\{c\})$ and $\psi $ satisfies \eqref{u2} in $[-1,1]\setminus u^{-1}\{c\} $ and \eqref{u1} holds for $\varphi\in H_0^1(-1,1),\ \text{supp}\ \varphi\subset [-1,1]\setminus u^{-1}\{c\}$.

For $y_0\in u^{-1}\{c\}$,  if $u'(y_0)=0$, then $y_0\neq\pm1,\ u''(y_0)\neq 0,$ and there exists $0<\d<1-|y_0|$ so that $u''(y)u''(y_0)>0$ for $|y-y_0|\leq \d$. By the argument in section 5.1, we know that $ \psi/(u-c)$ is bounded in $[y_0-\d,y_0+\d],$ and then \eqref{u1} holds for $\varphi\in H_0^1(-1,1),\ \text{supp}\ \varphi\subset(y_0-\d,y_0+\d)$. If $u'(y_0)\neq0,$ there exists $0<\d<1$ so that $u'(y)u'(y_0)>0$ for $|y-y_0(1-\d)|\leq\d$. Then Lemma \ref{lem:ub1} ensures that \eqref{u1} holds for $\varphi\in H_0^1(-1,1),\ \text{supp}\ \varphi\subseteq[y_0(1-\d)-\d,y_0(1-\d)+\d]$. Thus, \eqref{u1} holds for all $\varphi\in H_0^1(-1,1)$. Since $\mathcal{R}_\al$ has no embedding eigenvalues, $\psi=0 $. Thus, $\psi_n\to0 $ in $H^1_{loc}((-1,1)\setminus u^{-1}\{c\}).$  Furthermore,
if $u'(y_0)=0$, Lemma \ref{lem:ub2}  gives $\psi_n\to0 $ in $H^1(y_0-\d,y_0+\d)$; if $u'(y_0)\neq0$, Lemma \ref{lem:ub1} gives $\psi_n\to0 $ in $H^1(y_0(1-\d)-\d,y_0(1-\d)+\d)$. Thus, $\psi_n\to0 $ in $H^1(-1,1)$, which leads to a contradiction.
\end{proof}

Now we establish the limiting absorption principle.

\begin{proposition}\label{prop:limit}
Assume that $\mathcal{R}_\al$ has no embedding eigenvalues. Then there exists $\Phi_{\pm}(\al,y,c)\in H_0^1(-1,1) $ for $c\in \mathrm{Ran}\, u,$ such that $\Phi(\al,\cdot,c\pm i\epsilon)\to \Phi_{\pm}(\al,\cdot,c) $ in $C([-1,1])$ as $\epsilon\to0+ $ and
\beno
\|\Phi_{\pm}(\al,\cdot,c)\|_{H^1(-1,1)}\leq C\|\om\|_{H^1(-1,1)}.
\eeno
\end{proposition}
\begin{proof}
We view $\Phi$ as the map
\beno
\Phi(\al,\cdot,\cdot): c\in \Om_{\epsilon_0}\setminus\R&\to& \Phi(\al,\cdot,c)\in C([-1,1]).
\eeno
Let us claim that the map $\Phi(\al,\cdot,\cdot)$ is uniformly continuous in $\Omega_+=\{c+i\epsilon|c\in \text{Ran}\, u,0<\epsilon\leq \epsilon_0/2\}$ and $\Omega_-=\{c-i\epsilon|c\in \text{Ran}\, u,$ $0<\epsilon\leq \epsilon_0/2\}.$ We first consider the domain $\Omega_+$.

If it is not true, then there exists $c_n^1,\, c_n^2\in \Omega_+$ and $\d>0$ so that
\beno
\|\Phi(\al,\cdot,c_n^1)-\Phi(\al,\cdot,c_n^1)\|_{C([-1,1])}>\d,\quad |c_n^1-c_n^2|\to 0.
\eeno
 By Proposition \ref{prop:ub}, $\Phi(\al,\cdot,c_n^j)\ (j=1,2,n>0) $ is bounded in $ H^1(-1,1)$. Up to a subsequence,  we can assume $\Phi(\al,\cdot,c_n^j)\rightharpoonup\Phi_j $ weakly in $ H^1(-1,1)$ for some $\Phi_j\in H^1(-1,1)$, and $c_n^j\to c_0$ for $j=1,2$. Moreover, $\Phi(\al,\cdot,c_n^j)\to\Phi_j,\ j=1,2 $ in $ C([-1,1])$. Thus, $\|\Phi_1-\Phi_2\|_{C([-1,1])}\geq\d$.

If $c_0\not\in \R$, then $\Phi_1,\, \Phi_2$ are solutions of  \eqref{u3}
with $c=c_0$. Thus, $0\neq\Phi_1-\Phi_2$ is a solution of  \eqref{u3} with $\om(y)=0$,  which is impossible by Proposition \ref{prop:ub}.

If $c_0\in \R$, then $c_0\in \text{Ran}\, u$. Following the arguments in Proposition \ref{prop:ub}, we have
\beno
\int_{-1}^1(\Phi_j'\varphi'+\al^2\Phi_j\varphi)dy+p.v.\int_{-1}^1\frac{(u''\Phi_j+{\om})\varphi}{u-c_0}dy
+i\pi\sum_{y\in u^{-1}\{c_0\},u'(y)\neq 0}\frac{((u''\Phi_j+\omega)\varphi)(y)}{|u'(y)|}=0
\eeno
 for $j=1,2, \varphi\in H_0^1(-1,1),\ \text{supp}\ \varphi\subset[-1,1]\setminus\{y: u(y)=c,\,u'(y)=0\}.$ Moreover, $\psi=\Phi_1-\Phi_2 $ satisfies \eqref{u1} with $c=c_0$ for any $\varphi\in H_0^1(-1,1)$, since $\psi/(u-c)$ is bounded near $y_0$ where $u(y_0)=c, u'(y_0)=0$(see section 5.1).
Hence, $c_0$ is an embedding eigenvalue of $\mathcal{R}_\al$, which leads to a contradiction.

By taking conjugation, the map $\Phi(\al,\cdot,\cdot)$ is uniformly continuous in $\Omega_-$.

Then the result follows easily from Proposition \ref{prop:ub}.
\end{proof}

\section{Linear damping and vorticity depletion}
In this section, we prove the linear damping and the vorticity depletion phenomena for shear flows in $\mathcal{K}$.

In terms of the stream function $\widehat{\psi}(t,\al,y)$,
the linearized Euler equations take as follows
\beno
\frac 1 {i\al}\partial_t\widehat{\psi}=(\partial_y^2-\al^2)^{-1}\big(u''(y)-u(y)(\partial_y^2-\al^2)\big)\widehat{\psi}=-\mathcal{R}_\al\widehat{\psi}
\eeno
with the initial data
\beno
\psi(0,\al, y)=(\al^2-\pa_y^2)^{-1}\widehat{\om}_0(\al,y).
\eeno
Let $\Omega$ be a simply connected domain including the spectrum $\sigma(\mathcal{R}_\al)$ of $\mathcal{R}_\al$.
Then we have
\beno
\widehat{\psi}(t,\al,y)=\frac{1}{2\pi i}\int_{\partial\Omega}
e^{-i\al tc}(c-\mathcal{R}_\al)^{-1}\widehat{\psi}(0,\al,y)dc.
\eeno
Due to $P_{\mathcal{R}_\al}\widehat{\psi}(0,\al,y)=0$, we have
\beno
\widehat{\psi}(t,\al,y)=\frac{1}{2\pi i}\int_{\partial\Omega_{\epsilon}}e^{-i\al tc}(c-\mathcal{R}_\al)^{-1}\widehat{\psi}(0,\al,y)dc,
\eeno
where $\Om_\epsilon$ is defined by (\ref{eq:Omega}) with $\epsilon$ sufficiently small.
Let $\Phi(\al,y,c)$ be the solution of  (\ref{u3})
with $\om=\widehat{\om}_0(\al,y)/(i\al)$. It is easy to see that
\beno
(c-\mathcal{R}_\al)^{-1}\widehat{\psi}(0,\al,y)=i\al\Phi(\al,y,c).
\eeno
This shows that
\beq\label{eq:solution formu}
\widehat{\psi}(t,\al,y)=\frac{1}{2\pi}\int_{\partial\Omega_\epsilon}\al\Phi(\al,y,c)e^{-i\al ct}dc.
\eeq

Now we prove the linear damping.

\begin{proof}[Proof of Theorem \ref{thm:general}]
By Proposition \ref{prop:ub}, $|\Phi(\al,y,c)|\leq C\|\widehat{\om}_0(\al,\cdot){/\al}\|_{H^1}{\leq C\|\widehat{\om}_0(\al,\cdot)\|_{H^1} }$, which along with Proposition \ref{prop:limit} ensures that
\begin{align*}
\widehat{\psi}(t,\al,y)
&=\lim_{\epsilon\to0+}\frac{1}{2\pi i}\int_{\partial\Omega_{\epsilon}}e^{-i\al tc}i\al\Phi(\al,y,c)dc\\
&=\frac{1}{2\pi i}\int_{\text{Ran}\,u}e^{-i\al tc}i\al\big(\Phi_{-}(\al,y,c)-\Phi_{+}(\al,y,c)\big)dc\\
&=\frac{1}{2\pi}\int_{\text{Ran}\, u}e^{-i\al tc}\widetilde{\Phi}(\al,y,c)dc,
\end{align*}
where $\widetilde{\Phi}(\al,y,c)=\al\big(\Phi_{-}(\al,y,c)-\Phi_{+}(\al,y,c)\big)$.

For $k, l=0,1$, we have
\beno
\partial_t^k\partial_y^l\widehat{\psi}(t,\al,y)=\frac{1}{2\pi}\int_{\text{Ran}\,u}(-i\al c)^ke^{-i\al tc}\partial_y^l\widetilde{\Phi}(\al,y,c)dc,
\eeno
from which and Plancherel's formula, we infer that
\begin{align*}
&\|\widehat{V}(t,\al,y)\|_{H_t^1L_y^2}^2
=\int_{\R}\big(\|\widehat{V}(t,\al,\cdot)\|_{L_y^2}^2
+\|\partial_t\widehat{V}(t,\al,\cdot)\|_{L_y^2}^2\big)dt\\
&=\int_{-1}^1\int_{\R}\big(\al^2|\widehat{\psi}(t,\al,y)|^2+|\partial_y\widehat{\psi}(t,\al,\cdot)|^2
+\al^2|\partial_t\widehat{\psi}(t,\al,y)|^2
+|\partial_t\partial_y\widehat{\psi}(t,\al,y)|^2\big)dtdy\\
&=\frac{1}{2\pi|\al|}\int_{-1}^1\int_{\text{Ran}\,u}(1+|\al c|^2)\big(\al^2|\widetilde{\Phi}(\al,y,c)|^2+|\partial_y\widetilde{\Phi}(\al,y,c)|^2\big)dcdy\\
&\leq C\int_{\text{Ran}\,u}\|\widetilde{\Phi}(\al,\cdot,c)\|_{H^1_y}^2dc
\leq C\int_{\text{Ran}\,u}\|\widehat{\om}_0(\al,\cdot)\|_{H^1_y}^2dc=C\|\widehat{\om}_0(\al,\cdot)\|_{H^1_y}^2.
\end{align*}
By Sobolev embedding, we deduce that as $t\to\infty,$
\beno
\|\widehat{V}(t,\al,\cdot)\|_{L_y^2}^2\leq C\int_{t-1}^{t+1}\big(\|\widehat{V}(s,\al,\cdot)\|_{L_y^2}^2
+\|\partial_s\widehat{V}(s,\al,\cdot)\|_{L_y^2}^2\big)ds\to0.
\eeno
This shows the linear damping.
\end{proof}

Next, we prove the vorticity depletion phenomena.

\begin{proof}[Proof of Theorem \ref{thm:vorticity}]
Let $W(\al,y,c)=i\al(\al^2-\pa_y^2)\Phi(\al,y,c)$. It follows from \eqref{eq:solution formu} that
\begin{align*}
\widehat{\om}(t,\al,y)=\frac{1}{2\pi i}\int_{\partial\Omega_{\epsilon}}e^{-i\al tc}W(\al,y,c)dc.
\end{align*}
By Lemma \ref{lem:ub3} and Proposition \ref{prop:ub}, we have for $c\in \Om_{\epsilon_0}\setminus\R$,
\begin{align*}
|W(\al,y_0,c)|
&\leq C|u(y_0)-c|^{-\f34}(\|\Phi(\al,\cdot,c)\|_{H^1_y}+\|\widehat{\om}_0(\al,\cdot)\|_{H^1_y})\\
&\leq C|u(y_0)-c|^{-\f34}\|\widehat{\om}_0(\al,\cdot)\|_{H^1_y},
\end{align*}
which gives
\beno
&&|W(\al,y_0,c)|\leq C\epsilon^{-\f34}\|\widehat{\om}_0(\al,\cdot)\|_{H^1_y}\quad\textrm{for}\  c\in \partial\Om_{\epsilon},\\
&&|W(\al,y_0,c\pm i\epsilon)|\leq C|u(y_0)-c|^{-\f34}\|\widehat{\om}_0(\al,\cdot)\|_{H^1_y}\quad  \textrm{for} \ c\in \text{Ran}\,u.
\eeno
Thus, $W(\al,y_0,\cdot\pm i\epsilon)$ is bounded in $L_c^p(\text{Ran}\,u)(1<p<4/3)$ for $0<\epsilon<\epsilon_0.$ There exists a subsequence $\epsilon_n\to0+$ and $W_{\pm}(c) \in L_c^p(\text{Ran}\,u)(1<p<4/3)$ so that $W(\al,y_0,\cdot\pm i\epsilon_n)\rightharpoonup W_{\pm}$ weakly in $L_c^p(\text{Ran}\,u)$. Then we deduce that
\begin{align*}
\widehat{\om}(t,\al,y_0)
&=\lim_{\epsilon\to0+}\frac{1}{2\pi i}\int_{\partial\Omega_{\epsilon}}e^{-i\al tc}W(\al,y_0,c)dc\\
&=\lim_{n\to\infty}\frac{1}{2\pi i}\int_{\text{Ran}\,u}\Big(e^{-i\al t(c-i\epsilon_n)}W(\al,y_0,c-i\epsilon_n)-e^{-i\al t(c+i\epsilon_n)}W(\al,y_0,c+i\epsilon_n)\Big)dc\\
&=\frac{1}{2\pi i}\int_{\text{Ran}\,u}e^{-i\al tc}\big(W_-(c)-W_+(c)\big)dc\to 0\quad \textrm{as}\quad  t\to\infty,
\end{align*}
where in the last step we used Riemann-Lebesgue lemma.
\end{proof}

\section{The Inhomogeneous Rayleigh Equation for symmetric flow}

In this section, we solve the inhomogeneous Rayleigh equation when $u\in C^4([-1,1])$ satisfies $(S)$:
\beq\label{eq:Rayleigh-Ihom}
\left\{
\begin{aligned}
&\Phi''-\al^2\Phi-\frac{u''}{u-c}\Phi=f,\\
&\Phi(-1)=\Phi(1)=0,
\end{aligned}
\right.
\eeq
where $c\in \Om_{\epsilon_0}\setminus D_0$. In what follows, we will suppress the variable $\al$ for simplicity.

\subsection{Representation formula of the solution}
We decompose $f$ into the odd part and even part, i.e.,
\beno
f=f_{{o}}+f_e,
\eeno
where
\beno
f_{{o}}(y,c)=\frac{f(y,c)-f(-y,c)}{2},\quad f_e(y,c)=\frac{f(y,c)+f(-y,c)}{2}.
\eeno
Let $\Phi_0$ and $\Phi_e$ be the solution of the following Rayleigh equations
\beno
\left\{
\begin{aligned}
&\Phi_{o}''-\al^2\Phi_{o}-\frac{u''}{u-c}\Phi_{o}=f_{o},\\
&\Phi_{o}(-1)=\Phi_{o}(1)=0,
\end{aligned}
\right.
\eeno
and
\beno
\left\{
\begin{aligned}
&\Phi_{e}''-\al^2\Phi_{e}-\frac{u''}{u-c}\Phi_{e}=f_e,\\
&\Phi_{e}(-1)=\Phi_{e}(1)=0.
\end{aligned}
\right.
\eeno
Thus, $\Phi=\Phi_o+\Phi_e$ is the solution of the inhomogeneous Rayleigh equation \eqref{eq:Rayleigh-Ihom} with $\Phi_o$ being an odd function and $\Phi_e$ being an even function. So, it suffices to solve the following Rayleigh equation in $[0,1]$:
\beq\label{eq:Rayleigh-Ihom-odd-half}
\left\{
\begin{aligned}
&\Phi_{o}''-\al^2\Phi_{o}-\frac{u''}{u-c}\Phi_{o}=f_{o},\\
&\Phi_{o}(0)=\Phi_{o}(1)=0,
\end{aligned}
\right.
\eeq
and
\beq\label{eq:Rayleigh-Ihom-even-half}
\left\{
\begin{aligned}
&\Phi_{e}''-\al^2\Phi_{e}-\frac{u''}{u-c}\Phi_{e}=f_e,\\
&\Phi_{e}'(0)=\Phi_{e}(1)=0.
\end{aligned}
\right.
\eeq
Let $\varphi(y)$ be a solution of the homogenous Rayleigh equation
\beno
\varphi''-\al^2\varphi-\frac{u''}{u-c}\varphi=0.
\eeno
Then the inhomogeneous Rayleigh equations \eqref{eq:Rayleigh-Ihom-odd-half} and \eqref{eq:Rayleigh-Ihom-even-half} are equivalent to
\ben\label{eq:Ray-ihom-o}
\left\{
\begin{aligned}
&\Big(\varphi^2\big(\frac{\Phi_o}{\varphi}\big)'\Big)'=f_o\varphi,\\
&\Phi_o(0)=\Phi_o(1)=0,
\end{aligned}
\right.
\een
and
\ben\label{eq:Ray-ihom-e}
\left\{
\begin{aligned}
&\Big(\varphi^2\big(\frac{\Phi_e}{\varphi}\big)'\Big)'=f_e\varphi,\\
&\Phi_e'(0)=\Phi_e(1)=0.
\end{aligned}
\right.
\een
In particular, let $\varphi=\phi$ be the solution of the homogeneous Rayleigh equation constructed in Proposition \ref{prop:Rayleigh-Hom}.
We deduce by integrating \eqref{eq:Ray-ihom-o} and \eqref{eq:Ray-ihom-e} twice and matching the boundary conditions that for $y\in [0,1]$,
\begin{align}\label{eq: Phi_o}
\Phi_{o}(y,c)=&\phi(y,c)\int_{0}^y\frac{1}{\phi(y',c)^2}
\int_{y_c}^{y'}f_o\phi(y'',c)dy''dy'
+\mu^{o}(c)\phi(y,c)\int_{0}^y\frac{1}{\phi(y',c)^2}dy'\nonumber\\
=&\phi(y,c)\int_{1}^y\frac{1}{\phi(y',c)^2}
\int_{y_c}^{y'}f_o\phi(y'',c)dy''dy'
+\mu^{o}(c)\phi(y,c)\int_{1}^y\frac{1}{\phi(y',c)^2}dy',
\end{align}
where
\beno
\mu^o(c)=\frac{-\int_{0}^1\frac{1}{\phi(y',c)^2}
\int_{y_c}^{y'}f_o\phi(y'',c)dy''dy'}{\int_{0}^{1}\frac{1}{\phi(y',c)^2}dy'},
\eeno
and for $y\in [0,1]$,
\begin{align}\label{eq: Phi_e}
\Phi_{e}(y,c)=&\phi(y,c)\int_{0}^y\frac{1}{\phi(y',c)^2}
\int_{y_c}^{y'}f_e\phi(y'',c)dy''dy'\nonumber\\
&+\mu^e(c)\phi(y,c)\int_{0}^y\frac{1}{\phi(y',c)^2}dy'
+\nu^e(c)\phi(y,c)\nonumber\\=&\phi(y,c)\int_{1}^y\frac{1}{\phi(y',c)^2}
\int_{y_c}^{y'}f_e\phi(y'',c)dy''dy'
+\mu^e(c)\phi(y,c)\int_{1}^y\frac{1}{\phi(y',c)^2}dy',
\end{align}
where $\mu^e$ and $\nu^e$ are determined by solving
\begin{align*}
&\int_{0}^1\frac{1}{\phi(y',c)^2}
\int_{y_c}^{y'}f_e\phi(y'',c)dy''dy'
+\mu^e(c)\int_{0}^1\frac{1}{\phi(y',c)^2}dy'
+\nu^e(c)=0,\\
&\nu^e(c)\phi(0,c)\phi'(0,c)+\int_{y_c}^{0}f_e\phi(y'',c)dy''
+\mu^e(c)=0.
\end{align*}
That is,
\begin{align*}
\nu^e(c)&=\f{\int_0^1\f{1}{\phi(y',c)^2}\int_{y_c}^{y'}f_e\phi(y'',c)\,dy''dy'-\Big(\int_{y_c}^0f_e\phi(y'',c)dy''\Big)
\Big(\int_0^1\f{1}{\phi(y',c)^2}dy'\Big)}{\phi(0,c)\phi'(0,c)\int_0^1\f{1}{\phi(y',c)^2}\,dy'-1},\\
\mu^e(c)&=\f{-\phi(0,c)\phi'(0,c)\int_{0}^1\frac{1}{\phi(y',c)^2}
\int_{y_c}^{y'}f_e\phi(y'',c)dy''dy'+\int_{y_c}^0f_e\phi(y'',c)dy''}
{\phi(0,c)\phi'(0,c)\int_0^1\f{1}{\phi(y',c)^2}\,dy'-1}.
\end{align*}

Thus, the solution $\Phi(y,c)$ of \eqref{eq:Rayleigh-Ihom}  can be written as
\ben\label{eq:Phi}
\Phi(y,c)=\Phi_o(y,c)+\Phi_e(y,c)
\een
with $\Phi_{{o}}$ and $\Phi_e$ given by  \eqref{eq: Phi_o} and \eqref{eq: Phi_e}.

\subsection{Representation formula of the limiting solution}
In this subsection, we give a precise representation formula of the limiting solution $\Phi_\pm(y,c)$ obtained in Proposition \ref{prop:limit} for the symmetric flow.

The following fact is classical \cite{Hor}.

\begin{lemma}\label{lem:limit1}
For $g\in H^2(a,b)$, the function
\beno
F(c)=\int_{a}^{b}\f{g(z)}{z-c}dz\quad \text{for}\,\,\, \mathrm{Im}c>0
\eeno
can be $C^1$ extended to the interval $(a,b) $ with
\beno
F(c)=-H(g\chi_{[a,b]})(c)
+i\pi g(c)\quad \text{for}\,\,c\in (a,b).
\eeno
Here $H(g)(c)=p.v.\int\f{g(z)}{c-z}dz$ is the Hilbert transform of $g$.
\end{lemma}

\begin{remark}\label{Rmk:limit1}
If $g\in H^2\big((-b,b)\setminus\{0\}\big)\cap C([-b,b])$, then $F(z)$
can be $C^1$ extended to the interval $(-b,b)\setminus\{0\} $ with
\beno
F(c)=-H(g\chi_{[-b,b]})(c)
+i\pi g(c)\quad \text{for}\,\,c\in (-b,b)\setminus\{0\}.
\eeno
\end{remark}

Let $v(y)$ be as in {\eqref{def:v}} and $\tc=v(y_c)$. We define
\begin{align*}
\Int(\varphi)(y)=\int_{0}^y\varphi(y')\,dy' \quad\textrm{for}\quad y\in [0,1],\quad \Int(\varphi)(y)=\Int(\varphi)(-y)\quad\textrm{for}\quad y\in [-1,0].
\end{align*}
We introduce
\begin{align}
\textrm{II}_{1,1}(\varphi)(c)
&=p.v.\int_0^1\f{\Int(\varphi)(y)-\Int(\varphi)(y_c)}{(u(y)-u(y_c))^2}\,dy\nonumber\\
&=\partial_c\Big(\int_0^1\frac{\Int(\varphi)(y)-\Int(\varphi)(y_c)}{u(y)-c}\,dy\Big)
+\frac{\varphi(y_c)}{u'(y_c)}p.v.\int_0^1\frac{1}{u(y)-c}\, dy\nonumber\\
&=\partial_c\Big(p.v.\int_0^1\frac{\Int(\varphi)(y)}{u(y)-c}dy\Big)
-\Int(\varphi)(y_c)\partial_c\Big(p.v.\int_0^1\frac{dy}{u(y)-c}\Big)\nonumber\\
&=\f{1}{2\tc}\partial_{\tc}\Big(\frac{1}{2\widetilde{c}}p.v.\int_{-v(1)}^{v(1)}\frac{\Int(\varphi)(v^{-1}(z))(v^{-1})'(z)}{z-\widetilde{c}}\,dz\Big)\nonumber\\
&\quad-\f{\Int(\varphi)\big(v^{-1}(\tc)\big)}{2\tc}\pa_{\tc}\Big(\frac{1}{2\widetilde{c}}p.v.\int_{-v(1)}^{v(1)}\f{(v^{-1})'(z)}{z-\tc}\,dz\Big).\label{eq:Pi-11}
\end{align}

Let us point out that $\Int(\varphi)(y)\notin H^2(-1,1)$ if $\varphi(0)\neq 0$. However, $\Int(\varphi){\circ v^{-1}}\in H^2\big(({-v(1),v(1)})\setminus\{0\}\big)\cap C([{-v(1),v(1)}])$.

\begin{lemma}\label{Lem: II_1 convergence}
Let $c_{\epsilon}=c+i\ep\in D_{\epsilon_0}$. Then for any $\varphi\in H^1(0,1)$,
\begin{align*}
&\lim_{\ep\to 0+}
\rho(c_{\epsilon})\int_{0}^1\f{\int_{y_{c}}^y\varphi(y')dy'}{(u(y)-c_{\epsilon})^2}dy
=\rho(c)\mathrm{II}_{1,1}(\varphi)(c)
+i\pi \f{\varphi(y_c)}{u'(y_c)^2}\rho(c),\\
&\lim_{\ep\to 0-}
\rho(c_{\epsilon})\int_{0}^1\f{\int_{y_{c}}^y\varphi(y')dy'}{(u(y)-c_{\epsilon})^2}dy
=\rho(c)\mathrm{II}_{1,1}(\varphi)(c)
-i\pi \f{\varphi(y_c)}{u'(y_c)^2}\rho(c),
\end{align*}
and
\begin{align*}
&\lim_{\ep \to 0+}
\rho(c_{\epsilon})\tc_{\epsilon}\int_{0}^1\f{dy}{(u(y)-c_{\epsilon})^2}
=\f{\rmA_1(c)-i\rmB(c)}{2v'(y_c)},\\
&\lim_{\ep\to 0-}
\rho(c_{\epsilon})\tc_{\epsilon}\int_{0}^1\f{dy}{(u(y)-c_{\epsilon})^2}
=\f{\rmA_1(c)+i\rmB(c)}{2v'(y_c)}.
\end{align*}
Here $\tc_{\epsilon}$ is a unique solution of $c_{\epsilon}-u(0)=\tc_{\epsilon}^2$ with $\mathrm{Im}\,\tc_{\epsilon}>0$.
\end{lemma}

\begin{proof}
Thanks to $u(y)-c_\epsilon=v(y)^2-\tc_{\epsilon}^2$, we have
\begin{align*}
\int_{0}^1\f{\int_{y_{c}}^y\varphi(y')dy'}{(u(y)-c_{\epsilon})^2}dy
&=p.v.\int_0^1\f{\Int(\varphi)(y)-\Int(\varphi)(y_{c})}{(u(y)-c_{\epsilon})^2}\,dy\\
&=\partial_{c}\Big(\int_0^1\frac{\Int(\varphi)(y)-\Int(\varphi)(y_{c})}{u(y)-{c_{\epsilon}}}\,dy\Big)
+\frac{\varphi(y_{c})}{u'(y_{c})}\int_0^1\frac{1}{u(y)-{c_{\epsilon}}}\, dy\\
&=\partial_{c}\Big(\int_0^1\frac{\Int(\varphi)(y)}{u(y)-{c_{\epsilon}}}dy\Big)
-\Int(\varphi)(y_{c})\partial_{c}\Big(p.v.\int_0^1\frac{dy}{u(y)-{c_{\epsilon}}}\Big)\\
&=\f{1}{2\tc_{\epsilon}}\partial_{\tc_{\epsilon}}\Big(\frac{1}{2\widetilde{c}_{\epsilon}}
p.v.\int_{-v(1)}^{v(1)}\frac{\Int(\varphi)(v^{-1}(z))(v^{-1})'(z)}{z-\widetilde{c}_{\epsilon}}\,dz\Big)\\
&\quad-\f{\Int(\varphi)\big(y_{c}\big)}{2\tc_{\epsilon}}\pa_{\tc_{\epsilon}}
\Big(\frac{1}{2\widetilde{c}_{\epsilon}}p.v.\int_{-v(1)}^{v(1)}\f{(v^{-1})'(z)}{z-\tc_{\epsilon}}\,dz\Big),
\end{align*}
from which and {Remark \ref{Rmk:limit1}}, we deduce the first limit.
Similarly, we have
\begin{align*}
&\rho(c_{\epsilon})\tc_{\epsilon}\int_{0}^1\f{dy}{(u(y)-c_{\epsilon})^2}
=\rho(c_{\epsilon})\tc_{\epsilon}\partial_{c}\Big(p.v.\int_0^1\frac{dy}{u(y)-{c_{\epsilon}}}\Big)\\
&=\f{\rho(c_{\epsilon})}{2}\partial_{\tc_{\epsilon}}
\Big(\frac{1}{2\widetilde{c}_{\epsilon}}p.v.\int_{-v(1)}^{v(1)}\f{(v^{-1})'(z)}{z-\tc_{\epsilon}}\,dz\Big)\\
&\rightarrow \f{\rho(c)}{2}\partial_{\tc}
\Big(\frac{1}{2\widetilde{c}}p.v.\int_{-v(1)}^{v(1)}\f{(v^{-1})'(z)}{z-\tc}\,dz+\f{\pi i(v^{-1})'(\tc)}{2\tc}\Big)
=\f{\rmA_1(c)-i\rmB(c)}{2v'(y_c)},
\end{align*}
as $c_{\epsilon}\to c=u(y_c)\in D_0,\  \mathrm{Im}\,c_{\epsilon}>0.$

The case of $\mathrm{Im}\, c_{\epsilon}<0$ is similar. We omit the details.
\end{proof}

We decompose
\beno
\widehat{\om}_0(y)=\widehat{\om}_o(y)+\widehat{\om}_e(y),
\eeno
where $\widehat{\om}_o(y)=\f{\widehat{\om}_{0}(y)-\widehat{\om}_{0}(-y)}{2}$ and $\widehat{\om}_{e}(y)=\f{\widehat{\om}_{0}(y)+\widehat{\om}_{0}(-y)}{2}$. We take $f$ in \eqref{eq:Rayleigh-Ihom} as follows
\beno
f=\f{\widehat{\om}_{o}(y)}{i\al (u(y)-c)}+\f{\widehat{\om}_{e}(y)}{i\al (u(y)-c)}=f_o+f_e.
\eeno
We denote
\begin{align*}
&\mathrm{C}_{e}(c)=\rho(c)\frac{\widehat{\omega}_{e}(y_c)}{u'(y_c)}\pi,\quad \mathrm{D}_{e}(c)=u'(y_c)\rho(c)\textrm{II}_1(\widehat{\om}_{e})(c),\\
&\mathrm{C}_{o}(c)=\rho(c)\frac{\widehat{\omega}_{o}(y_c)}{u'(y_c)}\pi,\quad \mathrm{D}_{o}(c)=u'(y_c)\rho(c)\textrm{II}_1(\widehat{\om}_{o})(c),\\
&\mathrm{E}_e(c)=\mathrm{E}({\widehat{\om}_e})(c)\eqdef\int_{y_c}^0\widehat{\om}_{e}\phi_1(y,c)\,dy,\\
&\mathrm{II}_1(\varphi)(c)=\mathrm{II}_{1,1}(\varphi)(c)+\mathrm{II}_{1,2}(\varphi)(c),
\end{align*}
where
\ben\label{eq:Pi-12}
\mathrm{II}_{1,2}(\varphi)(c)=\int_0^1\int_{y_c}^{z}{\varphi}(y)
\left(\f{1}{(u(z)-c)^2}\left(\f{\phi_1(y,c)\,}{\phi_1(z,c)^2}-1\right)\right)\,dydz.
\een

For $c\in D_0$, we introduce
\begin{align*}
\Phi_\pm^{o}(y,c)
\triangleq\left\{
\begin{array}{l}
\phi\int_0^y\frac{1}{\phi(z,c)^2}\int_{y_c}^z\phi f_o(y',c)dy'dz
+\mu_{\pm}^o(c)\phi\int_0^y\frac{1}{\phi(y',c)^2}dy'\quad 0\leq y\leq y_c,\\
\phi\int_1^y\frac{1}{\phi(z,c)^2}\int_{y_c}^z\phi f_o(y',c)dy'dz
+\mu_{\pm}^o(c)\phi\int_1^y\frac{1}{\phi(y',c)^2}dy'\quad y_c\leq y\leq 1,
\end{array}
\right.
\end{align*}
where
\begin{align*}
\mu_{+}^o(c)&=\frac{1}{\al}\frac{iu'(y_c)\rho(c)\mathrm{II}_1(\widehat{\om}_{o})(c)
-\rho(c)\frac{\widehat{\om}_{o}(y_c)}{u'(y_c)}\pi}
{\rmA_1(c)-i\pi\rho(c)\frac{u''(y_c)}{u'(y_c)^2}+u'(y_c)\rho(c)\mathrm{II}_3(c)}\\
&=\f{1}{\al}\f{-\rmC_o(c)+i\rmD_{o}(c)}{\rmA(c)-i\rmB(c)},\\
\mu_{-}^o(c)&=\frac{1}{\al}\frac{iu'(y_c)\rho(c)\mathrm{II}_1(\widehat{\om}_{o})(c)
+\rho(c)\frac{\widehat{\om}_{o}(y_c)}{u'(y_c)}\pi}
{\rmA_1(c)+i\pi\rho(c)\frac{u''(y_c)}{u'(y_c)^2}+u'(y_c)\rho(c)\mathrm{II}_3(c)}\\
&=\f{1}{\al}\f{\rmC_o(c)+i\rmD_{o}(c)}{\rmA(c)+i\rmB(c)}.
\end{align*}
We also introduce
\begin{align*}
\Phi_{\pm}^{e}(y,c)
\triangleq
\left\{
\begin{aligned}
&\phi\int_0^y\f{\int_{y_c}^{y'}\phi f_e(y'',c)dy''}{\phi(y',c)^2}dy'
+\mu_{\pm}^e(c)\phi\int_0^y\f{1}{\phi(y',c)^2}dy'
+\nu_{\pm}^e(c)\phi(y,c)
\quad0\leq y\leq y_c,\\
&\phi\int_1^y\f{\int_{y_c}^{y'}\phi f_e(y'',c)dy''}{\phi(y',c)^2}dy'
+\mu_{\pm}^e(c)\phi\int_1^y\f{1}{\phi(y',c)^2}dy' \quad y_c\leq y\leq 1,
\end{aligned}
\right.
\end{align*}
where
\begin{align*}
\mu_+^e(c)&=\f{1}{\al}\f{\phi(0,c)\phi'(0,c)(i\rmD_{e}-\rmC_{e})(c)-i\rmE_e(c)u'(y_c)\rho(c)}{\phi(0,c)\phi'(0,c)(\rmA-i\rmB)(c)-\rho(c)u'(y_c)},\\
\mu_-^e(c)&=\f{1}{\al}\f{\phi(0,c)\phi'(0,c)(i\rmD_{e}+\rmC_{e})(c)-i\rmE_e(c)u'(y_c)\rho(c)}{\phi(0,c)\phi'(0,c)(\rmA+i\rmB)(c)-\rho(c)u'(y_c)},\\
\nu_+^e(c)&=-\f{1}{\al}\f{i\rmD_{e}(c)-\rmC_{e}(c)-i\rmE_e(c)(\rmA-i\rmB)(c)}{\phi(0,c)\phi'(0,c)(\rmA-i\rmB)(c)-u'(y_c)\rho(c)},\\
\nu_-^e(c)&=-\f{1}{\al}\f{i\rmD_{e}(c)+\rmC_{e}(c)-i\rmE_e(c)(\rmA+i\rmB)(c)}{\phi(0,c)\phi'(0,c)(\rmA+i\rmB)(c)-u'(y_c)\rho(c)}.
\end{align*}
It is easy to see that
\beno
\phi(0,c)\phi'(0,c)(\rmA+i\rmB)(c)-u'(y_c)\rho(c)=\phi_1(0,c)\phi_1'(0,c)(u(0)-c)(\rmA_2+i\rmB_2)(c).
\eeno
Thus, if $\cR_\al$ has no embedding eigenvalues, $\Phi_\pm^{{o}}$ and $\Phi_\pm^e$ are well-defined by Proposition \ref{prop:spectral}.

For $y\in [-1,0]$,  we let
\beno
\Phi_{\pm}^{o}(y,c)=-\Phi_{\pm}^{o}(-y,c),\quad \Phi^e_{\pm}(y,c)=\Phi^e_{\pm}(-y,c),
\eeno
and for $y\in [-1,1]$,
\beno
\Phi_{\pm}(y,c)=\Phi^o_{\pm}(y,c)+\Phi^e_{\pm}(y,c).
\eeno

\begin{proposition}\label{prop:psi-conv}
Let $\Phi(y,c)$ be a solution of \eqref{eq:Rayleigh-Ihom} with $f=\f {\widehat{\om}_0(y)} {i\al(u-c)}$ for $\widehat{\om}_0\in L^2(-1,1)$ given by  the formula \eqref{eq:Phi}. If $\mathcal{R}_\al$ has no embedding eigenvalues, then $\Phi(y,c)$ is well-defined for $c\in \Om_{\epsilon_0}$ small enough. Moreover, it holds that for any $(y,c)\in[0,1]\times D_0$ and $y\neq y_c$,
\beno
\lim_{\epsilon\to 0+}\Phi(y,c_\epsilon)=\Phi_+(y,c),\quad
\lim_{\epsilon\to 0-}\Phi(y,c_\epsilon)=\Phi_-(y,c),
\eeno
where $c_\epsilon=c+i\epsilon\in D_{\epsilon_0}$.
\end{proposition}


\begin{proof}
By Proposition \ref{prop:Rayleigh-Hom},  $\phi_1(y,c)$ is continuous for $(y,c)\in [0,1]\times\Om_{\epsilon_0}$ and for $\epsilon_0$ small enough,
$$|\phi_1(y,c)|>\f12,\quad |\phi_1(y,c)-1|\le C|y-y_c|^2.$$
Thus, for $(y,c)\in [0,1]\times\Om_{\epsilon_0}$, there exists a constant $C$ so that
\beno
\left|\f{\rho(c)}{(u(y)-c)^2}\Big(\f{1}{\phi_1(y,c)^2}-1\Big)\right|\leq C,
\eeno
and for $(y,c)\in [0,1]\times\Om_{\epsilon_0}$ with  $|z-y_c|\leq |y-y_c|$,
\beno
\left|\f{\rho(c)}{(u(y)-c)^2}\Big(\f{\phi_1(z,c)}{\phi_1(y,c)^2}-1\Big)\right|\leq C.
\eeno
This implies that  as $c_{\epsilon}\to u(y_c)=c$,
\begin{align*}
&\rho(c_{\epsilon})\int_0^1\f{1}{(u(y)-c_{\epsilon})^2}\Big(\f{1}{\phi_1(y,c_{\epsilon})^2}-1\Big)dy\to \rho(c)\mathrm{II}_3(c),\\
&\rho(c_{\epsilon})\int_0^1\int_{y_c}^y\f{\varphi(z)}{(u(y)-c_{\epsilon})^2}\Big(\f{\phi_1(z,c)}{\phi_1(y,c_{\epsilon})^2}-1\Big)dzdy\to \rho(c)\mathrm{II}_{1,2}(\varphi)(c),
\end{align*}
from which and Lemma \ref{Lem: II_1 convergence}, it follows that as $\ep\to 0\pm$,
\beno
u'(y_c)\rho(c_\ep)\int_{0}^{1}\frac{1}{\phi(y',c_\ep)^2}dy'\to \rho(c)\mathrm{II}_3(c)+\rmA_1(c)\pm i\rmB(c),
\eeno
and
\beno
&&i\al u'(y_c)\rho(c_\ep)\int_{0}^1\frac{1}{\phi(y',c_\ep)^2}
\int_{y_c}^{y'}f_o\phi(y'',c_\ep)dy''dy'\\
&&\to \rho(c)u'(y_c)\mathrm{II}_{1,2}(\widehat\om)(c)+\rho(c)u'(y_c)\mathrm{II}_{1,1}(\widehat\om)(c)
\pm i\pi\f {\widehat\om(y_c)} {u'(y_c)}\rho(c).
\eeno
Thus, we get
\beno
\lim_{\ep\to 0\pm}\mu^o(c_\ep)=\mu_\pm^o(c).
\eeno
Similarly, we have
\begin{align*}
\lim_{\ep\to 0\pm}\mu^e(c_\ep)=\mu_\pm^e(c),\quad \lim_{\ep\to 0\pm}\nu^e(c_\ep)=\nu_\pm^e(c).
\end{align*}

The following convergence is obvious: for $0\le y<y_c$, as $c_{\epsilon}\to u(y_c)$,
\begin{align*}
&\phi(y,c_\epsilon)\int_{0}^y\frac{1}{\phi(y',c_\epsilon)^2}
\int_{y_{c}}^{y'}f\phi(y'',c_\epsilon)dy''dy'\to \phi(y,c)\int_{0}^y\frac{1}{\phi(y',c)^2}
\int_{y_c}^{y'}f\phi(y'',c)dy''dy',\\
&\phi(y,c_\epsilon)\int_{0}^y\frac{1}{\phi(y',c_\epsilon)^2}dy'\to \phi(y,c)\int_{0}^y\frac{1}{\phi(y',c)^2}dy',
\end{align*}
 and   for $y_c<y\le 1$, as $c_{\epsilon}\to u(y_c)$,
\begin{align*}
&\phi(y,c_\epsilon)\int_{1}^y\frac{1}{\phi(y',c_\epsilon)^2}
\int_{y_c}^{y'}f\phi(y'',c_\epsilon)dy''dy'\to \phi(y,c)\int_{1}^y\frac{1}{\phi(y',c)^2}
\int_{y_c}^{y'}f\phi(y'',c)dy''dy',\\
&\phi(y,c_\epsilon)\int_{1}^y\frac{1}{\phi(y',c_\epsilon)^2}dy'\to \phi(y,c)\int_{1}^y\frac{1}{\phi(y',c)^2}dy'.
\end{align*}

With the above information, the proposition follows easily.
\end{proof}


\section{Decay estimates of the velocity for symmetric flow}

In this section, we establish the decay estimates of the velocity for symmetric flow by using the dual method.

We introduce that for $y\in [0,1]$,
\begin{align*}
\widetilde{\Phi}(y,c)&=\widetilde{\Phi}_{o}(y,c)+\widetilde{\Phi}_{e}(y,c)\\
&=\left\{
\begin{aligned}
&(\mu_{-}^o(c)-\mu_{+}^o(c))\phi(y,c)\int_0^y\frac{1}{\phi(z,c)^2}dz\\
&(\mu_{-}^o(c)-\mu_{+}^o(c))\phi(y,c)\int_1^y\frac{1}{\phi(z,c)^2}dz
\end{aligned}
\right.\\
&\quad+\left\{
\begin{aligned}
&(\mu_-^e(c)-\mu_+^e(c))\phi(y,c)\int_0^y\frac{1}{\phi(z,c)^2}dz
+(\nu_-^e(c)-\nu_+^e(c))\phi(y,c)\quad0\leq y<y_c,\\
&(\mu_-^e(c)-\mu_+^e(c))\phi(y,c)\int_1^y\frac{1}{\phi(z,c)^2}dz\quad y_c<y\leq 1.
\end{aligned}
\right.
\end{align*}
For $y\in [-1,0]$,  we define
\beno
\widetilde{\Phi}(y,c)=-\widetilde{\Phi}_{o}(-y,c)+\widetilde{\Phi}_{e}(-y,c).
\eeno
Furthermore, we have
\begin{align*}
\mu_{-}^o(c)-\mu_{+}^o(c)&=\f{2}{\al}\f{\mathrm{A}\mathrm{C}_{o}+\mathrm{B}\mathrm{D}_{o}}{\mathrm{A}^2+\mathrm{B}^2}
\eqdef\f{2}{\al}\rho(c)\mu_1(c),\\
\nu_-^e(c)-\nu_+^e(c)
&=-\f{2}{\al}\f{\phi(0,c)\phi'(0,c)(\rmA\rmC_{e}+\rmB\rmD_{e})-u'(y_c)\rho(c)(\rmB\rmE_e+C_{e})}
{(\phi(0,c)\phi'(0,c)\rmA-u'(y_c)\rho(c))^2+\phi(0,c)^2\phi'(0,c)^2\rmB^2}\eqdef \f{2}{\al}\nu_1(c),\\
\mu_-^e(c)-\mu_+^e(c)
&=-\f{2}{\al}\phi(0,c)\phi'(0,c)\nu_1(c)\eqdef \f{2}{\al}\mu_2(c).
\end{align*}

It follows from Proposition \ref{prop:ub} and Proposition \ref{prop:psi-conv}  that
\begin{align*}
\widehat{\psi}(t,\al,y)=&\frac{1}{2\pi}\int_{u(0)}^{u(1)}\al\widetilde{\Phi}(y,c)e^{-i\al ct}dc\nonumber\\
=&\frac{1}{2\pi}\int_{u(0)}^{u(1)}\al\widetilde{\Phi}_o(y,c)e^{-i\al ct}dc+\frac{1}{2\pi}\int_{u(0)}^{u(1)}\al\widetilde{\Phi}_e(y,c)e^{-i\al ct}dc\nonumber\\
\eqdef&\widehat{\psi}_o(t,\al,y)+\widehat{\psi}_e(t,\al,y).
\end{align*}

\subsection{The odd part}
For $f=g''-\al^2g$ with $g\in H^2(0,1)\cap H^1_0(0,1),$ we have
\begin{align*}
&\int_0^1\widehat{\psi}_o(t,\al,y)f(y)dy\\
&=\dfrac{1}{\pi}\int_0^1f(y)\int_{u(0)}^{u(y)}\rho(c)\mu_1(c)\phi(y,c)
\int_1^y\dfrac{1}{\phi(z,c)^2}dze^{-i\al ct}dc dy\\
&\quad+\dfrac{1}{\pi}\int_0^1f(y)\int_{u(y)}^{u(1)}\rho(c)\mu_1(c)\phi(y,c)
\int_0^y\dfrac{1}{\phi(z,c)^2}dze^{-i\al ct}dc dy\\
&=-\f{1}{\pi}\int_{u(0)}^{u(1)}\rho(c)\mu_1(c)e^{-i\al ct}
\int_0^1\dfrac{\int_{y_c}^zf(y)\phi(y,c)dy}{\phi(z,c)^2}dz dc.
\end{align*}
First of all, we have
$$
\rho(c)\mu_1(c)=\pi\f{\rmA(c)\rho(c)\frac{\widehat{\omega}_{o}(y_c)}{u'(y_c)}+\rho(c)^2\frac{u''(y_c)}{u'(y_c)}\textrm{II}_1(\widehat{\om}_{o})(c)}{\rmA(c)^2+\rmB(c)^2}.
$$
We write
\begin{align*}
&\int_0^1\dfrac{\int_{y_c}^zf(y)\phi(y,c)dy}{\phi(z,c)^2}dz\\
&=\int_0^1\dfrac{\int_{y_c}^z(g''(y)-\al^2g(y))\phi(y,c)dy}{\phi(z,c)^2}dz\\
&=\int_0^1\dfrac{\int_{y_c}^zg(y)(\phi''-\al^2\phi)(y,c)dy+g'(z)\phi(z,c)-g(z)\phi'(z,c)+g(y_c)\phi'(y_c,c)}{\phi(z,c)^2}dz\\
&=\int_0^1\left[\dfrac{\int_{y_c}^zg(y)u''(y)\phi_1(y,c)dy}{\phi(z,c)^2}+\Big(\dfrac{g(z)}{\phi(z,c)}\Big)'
+\dfrac{g(y_c)u'(y_c)}{\phi(z,c)^2}\right]dz\\
&=\mathrm{II}_1(gu'')(c)+p.v.\int_0^1\left[\left(\dfrac{g(z)}{\phi(z,c)}\right)'
+\dfrac{g(y_c)u'(y_c)}{\phi(z,c)^2}\right]dz.
\end{align*}
For $c\in(u(0),u(1))$, we have
\beno
\f{g(z)}{\phi(z,c)}-\f{g(y_c)}{u(z)-c}=\f{g(z)(1-\phi_1(z,c))}{(u(z)-c)\phi_1(z,c)}+\f{g(z)-g(y_c)}{u(z)-c}\in C([0,1]),
\eeno
which along with $g(0)=g(1)=0$ gives
\begin{align*}
&p.v.\int_0^1\left[\left(\f{g(z)}{\phi(z,c)}\right)'
+\f{g(y_c)u'(y_c)}{\phi(z,c)^2}\right]dz\\
&=\left(\f{g(z)}{\phi(z,c)}-\f{g(y_c)}{u(z)-c}\right)\bigg|_0^1
+p.v.\int_0^1\left[\left(\f{g(y_c)}{u(z)-c}\right)'
+\f{g(y_c)u'(y_c)}{\phi(z,c)^2}\right]dz\\
&=\f{g(y_c)}{u(0)-c}-\f{g(y_c)}{u(1)-c}+
p.v.\int_0^1\left[-\f{g(y_c)(u'(z)-u'(y_c))}{(u(z)-c)^2}
+\f{g(y_c)u'(y_c)}{(u(z)-c)^2}\left(\f{1}{\phi(z,c)^2}-1\right)\right]dz,\\
&=\f{g(y_c)(u(0)-u(1))}{\rho(c)}-g(y_c)\mathrm{II}_2(c)+g(y_c)u'(y_c)\mathrm{II}_3(c)
=\f{g(y_c)}{\rho(c)}\rmA(c).
\end{align*}
Thus, we deduce that
\begin{align}\label{eq: estimate psi}
\int_0^1\widehat{\psi}_o(t,\al,y)f(y)dy
=-\int_{u(0)}^{u(1)}K_o(c,\al)e^{-i\al ct}dc,
\end{align}
where
\begin{align*}
K_o(c,\al)=\f{\big(\mathrm{A}(c){\widehat{\omega}_{o}(y_c)}
+\rho(c){u''(y_c)}\textrm{II}_1(\widehat{\om}_{o})(c)\big)\big(\mathrm{II}_1(gu'')(c)\rho(c)+{g(y_c)}\mathrm{A}(c)\big)}{(\mathrm{A}(c)^2+\mathrm{B}(c)^2){u'(y_c)}}.
\end{align*}
We introduce
\begin{align*}
&\Lambda_1(\varphi)(c)
=\Lambda_{1,1}(\varphi)(c)+\Lambda_{1,2}(\varphi)(c),\\
&\Lambda_2(\varphi)(c)
=\Lambda_{2,1}(\varphi)(c)+\Lambda_{2,2}(\varphi)(c),
\end{align*}
where
\begin{align*}
&\Lambda_{1,1}(\varphi)(c)=\rmA_1(c)\varphi(y_c)+\rho(c)u''(y_c)\mathrm{II}_{1,1}(\varphi)(c),\\
&\Lambda_{1,2}(\varphi)(c)=\rho(c)u''(y_c)\mathrm{II}_{1,2}(\varphi)+u'(y_c)\rho(c)\mathrm{II}_3(c)\varphi(y_c),\\
&\Lambda_{2,1}(\varphi)(c)=\rmA_1(c)\varphi(y_c)+\rho(c)\mathrm{II}_{1,1}(u''\varphi)(c),\\
&\Lambda_{2,2}(\varphi)(c)=\rho(c)\mathrm{II}_{1,2}(u''\varphi)(c)+u'(y_c)\rho(c)\mathrm{II}_3(c)\varphi(y_c).
\end{align*}
Then we have
\begin{align}\label{eq:Ko}
K_o(c,\al)=\f{\Lambda_1({\widehat{\om}_{o}})(c)\Lambda_2(g)(c)}{(\mathrm{A}(c)^2+\mathrm{B}(c)^2){u'(y_c)}}.
\end{align}

\subsection{The even part}
For $f=g''-\al^2g$ with $g\in H^2(0,1)$ and $g'(0)=g(1)=0$, we have
\begin{align*}
&\int_0^1\widehat{\psi}_e(t,\al,y)f(y)dy\\
&=\dfrac{1}{\pi}\int_0^1f(y)\int_{u(0)}^{u(y)}\mu_2(c)\phi(y,c)
\int_1^y\dfrac{1}{\phi(z,c)^2}dze^{-i\al ct}dc dy\\
&\quad+\dfrac{1}{\pi}\int_0^1f(y)\int_{u(y)}^{u(1)}\mu_2(c)\phi(y,c)
\int_0^y\dfrac{1}{\phi(z,c)^2}dze^{-i\al ct}dc dy\\
&\quad+\dfrac{1}{\pi}\int_0^1f(y)\int_{u(y)}^{u(1)}\nu_1(c)\phi(y,c)e^{-i\al ct}dc dy\\
&=-\f{1}{\pi}\int_{u(0)}^{u(1)}\mu_2(c)e^{-i\al ct}
\int_0^1\dfrac{\int_{y_c}^zf(y)\phi(y,c)dy}{\phi(z,c)^2}dz dc\\
&\quad+\f{1}{\pi}\int_{u(0)}^{u(1)}\nu_1(c)e^{-i\al ct}
\int_0^{y_c}f(y)\phi(y,c)dy dc\\
&=-\f{1}{\pi}\int_{u(0)}^{u(1)}\mu_2(c)e^{-i\al ct}
\Big(\int_0^1\dfrac{\int_{y_c}^zf(y)\phi(y,c)dy}{\phi(z,c)^2}dz+ \dfrac{\int_0^{y_c}f(y)\phi(y,c)dy}{\phi(0,c)\phi'(0,c)}\Big)dc.
\end{align*}
Using the facts that
\begin{align*}
&\int_0^1\dfrac{\int_{y_c}^zf(y)\phi(y,c)dy}{\phi(z,c)^2}dz
=\mathrm{II}_1(gu'')(c)+p.v.\int_0^1\left[\left(\dfrac{g(z)}{\phi(z,c)}\right)'
+\dfrac{g(y_c)u'(y_c)}{\phi(z,c)^2}\right]dz,\\
&p.v.\int_0^1\left[\left(\f{g(z)}{\phi(z,c)}\right)'
+\f{g(y_c)u'(y_c)}{\phi(z,c)^2}\right]dz
=\f{g(y_c)}{\rho(c)}\rmA(c)-\f{g(0)}{\phi(0,c)},\\
&\int_0^{y_c}f(y)\phi(y,c)dy=-\rmE(gu'')(c)-g(y_c)u'(y_c)+g(0)\phi'(0,c),
\end{align*}
we deduce that
\begin{align}
&\int_0^1\widehat{\psi}_e(t,\al,y)f(y)dy\nonumber\\
&=-\f{1}{\pi}\int_{u(0)}^{u(1)}\Big[\mu_2(c)\Big(\f{g(y_c)\rmA(c)}{\rho(c)}+\mathrm{II}_1(gu'')(c)\Big)+\nu_1(c)\big(\rmE(gu'')(c)+g(y_c)u'(y_c)\big)\Big]e^{-i\al ct}dc\nonumber\\
&=-\int_{u(0)}^{u(1)}K_e(c,\al)e^{-i\al ct}dc,\label{eq:psieven}
\end{align}
where
\beno
\pi K_e(c,\al)=\mu_2(c)\Big(\f{g(y_c)\rmA(c)}{\rho(c)}+\mathrm{II}_1(gu'')(c)\Big)+\nu_1(c)\big(\rmE(gu'')(c)+g(y_c)u'(y_c)\big).
\eeno
Recall that
\begin{align*}
\nu_1(c)&=-\f{\phi(0,c)\phi'(0,c)(\rmA\rmC_{e}+\rmB\rmD_{e})-u'(y_c)\rho(c)(\rmB\rmE_e+C_{e})}
{(\phi(0,c)\phi'(0,c)\rmA-u'(y_c)\rho(c))^2+\phi(0,c)^2\phi'(0,c)^2\rmB(c)^2}\\
&=-\f{-\rho_1(\rmA\rmC_{e}+\rmB\rmD_{e})+J(\rmB\rmE_e+C_{e})}
{\big[(-\rho_1\rmA+J)^2+\rho_1^2\rmB^2\big]\phi(0,c)\phi_1'(0,c)}\\
&=-\pi\f{-\rho_1\big(\rmA\rho\frac{\widehat{\omega}_{e}(y_c)}{u'(y_c)}+\rho\frac{u''(y_c)}{u'(y_c)^2}u'(y_c)
\rho\textrm{II}_1(\widehat{\om}_{e})\big)+J\big(\rho\frac{u''(y_c)}{u'(y_c)^2}\rmE_e+\rho\frac{\widehat{\omega}_{e}(y_c)}{u'(y_c)}\big)}
{\big[(-\rho_1\rmA+J)^2+\rho_1^2\rmB^2\big]\phi(0,c)\phi_1'(0,c)}\\
&=-\f{\pi\rho(c)}{u'(y_c)}\f{-\rho_1\big(\rmA{\widehat{\omega}_{e}(y_c)}+{u''(y_c)}
\rho\textrm{II}_1(\widehat{\om}_{e})\big)+J\big(\frac{u''(y_c)}{u'(y_c)}\rmE_e+{\widehat{\omega}_{e}(y_c)}\big)}
{\big[(-\rho_1\rmA+J)^2+\rho_1^2\rmB^2\big]\phi(0,c)\phi_1'(0,c)},\\
\mu_2(c)&=-\phi(0,c)\phi'(0,c)\nu_1(c),
\end{align*}
with $\rho_1=c-u(0)$ and
\begin{align*}
J(c)=\f{-u'(y_c)\rho(c)}{\phi(0,c)\phi_1'(0,c)}=\f{u'(y_c)(u(1)-c)}{\phi_1(0,c)\phi_1'(0,c)}.
\end{align*}
Let
\begin{align*}
\Lambda_3(\widehat{\omega}_{e})(c)
&=-\rho_1\big(\rmA(c){\widehat{\omega}_{e}(y_c)}+{u''(y_c)}
\rho\textrm{II}_1(\widehat{\om}_{e})(c)\big)+J\Big(\frac{u''(y_c)}{u'(y_c)}\rmE_e(c)+{\widehat{\omega}_{e}(y_c)}\Big)\\
&=-\rho_1(c)\Lambda_1(\widehat{\om}_e)(c)
+\Lambda_{3,1}(\widehat{\om}_e)(c),\\
\Lambda_4(g)(c)
&=-\rho_1\big(\rmA{g}(y_c)+
\rho\textrm{II}_1(gu'')\big)(c)+J\Big(\frac{\rmE(gu'')}{u'(y_c)}+g(y_c)\Big)(c)\\
&=-\rho_1(c)\Lambda_2(g)(c)
+\Lambda_{4,1}(g)(c),
\end{align*}
with
\beno
&&\Lambda_{3,1}(\widehat{\om}_e)(c)=J(c)\Big(\frac{u''(y_c)}{u'(y_c)}\rmE(\widehat{\om}_e)(c)+{\widehat{\omega}_{e}(y_c)}\Big),\\
&&\Lambda_{4,1}(g)(c)=J(c)\Big(\frac{\rmE(gu'')(c)}{u'(y_c)}+g(y_c)\Big).
\eeno
Thus, we get
\begin{align}\label{eq:Ke}
K_e(c,\al)=\f{\Lambda_3(\widehat{\omega}_{e})(c)\Lambda_4(g)(c)}
{u'(y_c)\big((-\rho_1\rmA+J)^2+\rho_1^2\rmB^2\big)(c)}
=\f{\Lambda_3(\widehat{\omega}_{e})(c)\Lambda_4(g)(c)}{u'(y_c)(\rmA_2^2+\rmB_2^2)(c)}.
\end{align}

\subsection{Decay estimates and scattering}

The decay estimates are based on the following regularity estimates of the kernel.

\begin{proposition}\label{prop:Ko}
Assume that $f=g''-\al^2g$ with $g\in H^2(0,1)\cap H^1_0(0,1)$ and
$\widehat{\om}_o(\al,y)=\f{1}{2}(\widehat{\om}_0(\al,y)-\widehat{\om}_0(\al,-y))\in H^2(0,1)$. Then it holds that
\beno
K_o(u(0),\al)=K_o(u(1),\al)=0,
\eeno
and there exists a constant $C$ independent of $\al$ so that
\begin{align*}
&\|K_o(\al,\cdot)\|_{L^1_c}
\leq C\|\widehat{\om}_o(\al,\cdot)\|_{L^2_y}\|g\|_{L^2},\\
&\|(\pa_cK_o)(\al,\cdot)\|_{L^1_c}
\leq C\|\widehat{\om}_o(\al,\cdot)\|_{H^1_y}\|g\|_{H^1},\\
&\|(\pa_c^2K_o)(\al,\cdot)\|_{L^1_c}
\leq C\al^{\f12}\|\widehat{\om}_o(\al,\cdot)\|_{H^2_y}\|f\|_{L^2}.
\end{align*}
\end{proposition}

\begin{proposition}\label{prop:Ke}
Assume that $f=g''-\al^2g$ with $g\in H^2(0,1)$ and $g'(0)=g(1)=0$, and
$\widehat{\om}_e(\al,y)=\f{1}{2}(\widehat{\om}_0(\al,y)+\widehat{\om}_0(\al,-y))\in H^2(0,1)$. Then we have
\beno
K_e(u(0),\al)=K_e(u(1),\al)=0,
\eeno
and  there exists a constant $C$ independent of $\al$ such that
\begin{align*}
&\|K_e(\al,\cdot)\|_{L^1_c}
\leq C\|\widehat{\om}_e(\al,\cdot)\|_{L^2_y}\|g\|_{L^2},\\
&\|(\pa_cK_e)(\al,\cdot)\|_{L^1_c}
\leq C\al^{\f12}\|\widehat{\om}_e(\al,\cdot)\|_{H^1_y}(\|g'\|_{L^2}+\al\|g\|_{L^2}),\\
&\|(\pa_c^2K_e)(\al,\cdot)\|_{L^1_c}
\leq C\al^{\f32}\|\widehat{\om}_e(\al,\cdot)\|_{H^2_y}\|f\|_{L^2}.
\end{align*}
\end{proposition}

With Proposition \ref{prop:Ko} and Proposition \ref{prop:Ke}, we are in a position to prove Theorem \ref{thm:main}. \smallskip

Using Proposition \ref{prop:Ko} and Proposition \ref{prop:Ke},
we get by integration by parts that
\begin{align*}
\|\widehat{\psi}_o(t,\al,\cdot)\|_{L^2_y}
&={2}\sup_{\|f\|_{L^2}=1}\left|\int_0^1\widehat{\psi}_o(t,\al,y)f(y)dy\right|\\
&={2}\sup_{\|f\|_{L^2}=1}\left|\int_{u(0)}^{u(1)}K_o(c,\al)e^{-i\al ct}dc\right|\\
&\leq C\f{1}{\al^{\f32}t^2}\|\widehat{\om}_o(\al,\cdot)\|_{H^2_y},
\end{align*}
and
\begin{align*}
\|\widehat{\psi}_e(t,\al,\cdot)\|_{L^2_y}
&={2}\sup_{\|f\|_{L^2}=1}\left|\int_0^1\widehat{\psi}_e(t,\al,y)f(y)dy\right|\\
&={2}\sup_{\|f\|_{L^2}=1}\left|\int_{u(0)}^{u(1)}K_e(c,\al)e^{-i\al ct}dc\right|\\
&\leq C\f{1}{\al^{\f12}t^2}\|\widehat{\om}_o(\al,\cdot)\|_{H^2_y}.
\end{align*}
Similarly, we have
\begin{align*}
&\al^2\|\widehat{\psi}_o(t,\al,\cdot)\|_{L^2_y}^2+\|\partial_y\widehat{\psi}_o(t,\al,\cdot)\|_{L^2_y}^2\\
&=-{2}\int_0^1\widehat{\psi}_o(t,\al,y)(\overline{\widehat{\psi}_o}''-\al^2\overline{\widehat{\psi}_o})(\al,y)dy
\leq \left\{
\begin{aligned}
&C\|\widehat{\om}_o(\al,\cdot)\|_{L^2_y}\|\widehat{\psi}_o(\al,\cdot)\|_{L^2_y},\\
&C\f{1}{\al t}\|\widehat{\om}_o(\al,\cdot)\|_{H^1_y}\|\widehat{\psi}_o(\al,\cdot)\|_{H^1_y},
\end{aligned}
\right.
\end{align*}
and
\begin{align*}
&\al^2\|\widehat{\psi}_e(t,\al,\cdot)\|_{L^2_y}^2+\|\partial_y\widehat{\psi}_e(t,\al,\cdot)\|_{L^2_y}^2\\
&=-{2}\int_0^1\widehat{\psi}_e(t,\al,y)(\overline{\widehat{\psi}_e}''-\al^2\overline{\widehat{\psi}_e})(\al,y)dy
\leq \left\{
\begin{aligned}
&C\|\widehat{\om}_e(\al,\cdot)\|_{L^2_y}\|\widehat{\psi}_e(\al,\cdot)\|_{L^2_y},\\
&C\f{1}{\al^{\f12}t}\|\widehat{\om}_e(\al,\cdot)\|_{H^1_y}\|\widehat{\psi}_e(\al,\cdot)\|_{H^1_y}.
\end{aligned}
\right.
\end{align*}
Let
\beno
V_{{o}}=\na^\perp(-\Delta)^{-1}\psi_{{o}},\quad V_e=\na^\perp(-\Delta)^{-1}\psi_e.
\eeno
Thus, we deduce that for $t\leq 1$,
\begin{align*}
\|\widehat{V}_o(t,\al,\cdot)\|_{L^2_y}
\leq& \al\|\widehat{\psi}_o(t,\al,\cdot)\|_{L^2_y}+\|\partial_y\widehat{\psi}_o(t,\al,\cdot)\|_{L^2_y}\\
\leq& C\al^{-1}\|\widehat{\om}_o(\al,\cdot)\|_{L^2_y},\\
\|\widehat{V}_e(t,\al,\cdot)\|_{L^2_y}
\leq& \al\|\widehat{\psi}_e(t,\al,\cdot)\|_{L^2_y}+\|\partial_y\widehat{\psi}_e(t,\al,\cdot)\|_{L^2_y}\\
\leq& C\al^{-1}\|\widehat{\om}_e(\al,\cdot)\|_{L^2_y},
\end{align*}
and for $t\geq 1$,
\begin{align*}
\|\widehat{V}_o(t,\al,\cdot)\|_{L^2_y}
\leq& \al\|\widehat{\psi}_o(t,\al,\cdot)\|_{L^2_y}+\|\partial_y\widehat{\psi}_o(t,\al,\cdot)\|_{L^2_y}\\
\leq& C\f{1}{\al t}\|\widehat{\om}_o(\al,\cdot)\|_{H^1_y},\\
\|\widehat{V}_e(t,\al,\cdot)\|_{L^2_y}
\leq& \al\|\widehat{\psi}_e(t,\al,\cdot)\|_{L^2_y}+\|\partial_y\widehat{\psi}_e(t,\al,\cdot)\|_{L^2_y}\\
\leq& C\f{1}{\al^{\f12}t}\|\widehat{\om}_e(\al,\cdot)\|_{H^1_y}.
\end{align*}
Thanks to $V=V_o+V_e$, we get
\beno
\|V(t)\|_{L^2_{x,y}}\leq \f{C}{\langle t\rangle}\|\om_0\|_{H^{-\f12}_xH^1_y}.
\eeno
Using $\|\widehat{V}^2\|_{L^2(I_0)}\leq C\al\|\widehat{\psi}\|_{L^2(I_0)}
\leq C\dfrac{\al^{\f12}}{t^2}\|\widehat{\om}_0\|_{H^2(I_0)}$, we get
\beno
\|V^2(t)\|_{L^2_{x,y}}\leq \f{C}{\langle t\rangle^2}\|{\om}_0\|_{H^{\f12}_xH^2_y}.
\eeno

The scattering part is the same as the monotonic case in \cite{WZZ1}. Thus, we omit the details.
\medskip

The remaining sections of this paper will be devoted to the proof of Proposition \ref{prop:Ko} and Proposition \ref{prop:Ke}.

\section{Estimates of some key quantities}

In this section, we present some estimates for some key quantities like $\rmA, \rmB$ etc. appeared in $K_o$ and $K_e$.

In the sequel, we denote by $\phi(y,c)$ the solution of the homogeneous Rayleigh equation given
by Proposition \ref{prop:Rayleigh-Hom} and $\phi_1(y,c)=\f {\phi(y,c)} {u(y)-c}$.

Let $I_v=(-v(1), v(1))$ with $v$ given by \eqref{def:v} and $\rho_0(c)=\rho(c)/u'(y_c)\sim y_c(1-y_c)$. We denote by $\|\cdot\|_{L^p_{\tc}}$ the norm of $L^p(I_v, d\tc)$, and $C$ a constant independent of $\al$, which may be different from line to line.

\subsection{Estimate of $\mathrm{II}_3$}
Recall that
\beno
\mathrm{II}_3(c)=\int_0^1\frac{1}{(u(y)-c)^2}\Big(\frac{1}{\phi_{1}(y,c)^2}-1\Big)dy\leq 0.
\eeno

We have following estimates for $\mathrm{II}_3.$
\begin{lemma}\label{lem:II3}
It holds that
\begin{align*}
&C^{-1}\min\Big\{\f{\al^2}{u'(y_c)},\f{\al}{u'(y_c)^2}\Big\}\leq -\mathrm{II}_3(c)\leq C\min\Big\{\f{\al^2}{u'(y_c)},\f{\al}{u'(y_c)^2}\Big\},\\
&\big|\rho(c)^k\pa_c^k\mathrm{II}_3(c)\big|\leq C\min\Big\{\f{\al^2}{u'(y_c)},\f{\al}{u'(y_c)^2}\Big\},\quad k=1,2.
\end{align*}
\end{lemma}
\begin{proof}
By Proposition \ref{prop:phi1}, we have
\beno
{\phi_1(y,c)\ge 1},\quad \phi_1(y,c)-1\le C\min\{\al^2|y-y_c|^2,1\}\phi_1(y,c),
\eeno
which gives
\begin{align*}
-\mathrm{II}_3(c)=&\int_0^1\frac{1}{(u(y)-c)^2}\Big(1-\frac{1}{\phi_{1}(y,c)^2}\Big)dy\\
\geq&C^{-1}\int_0^1\frac{\al^2|y-y_c|^2}{(u(y)-c)^2}\chi_{\{|y-y_c|\leq\frac{1}{\al}\}}dy\\
\geq&
\left\{
\begin{array}{ll}
C^{-1}\int_{0}^{y_c}\dfrac{\al^2|y-y_c|^2}{u'(y_c)^2|y-y_c|^2}dy\geq \dfrac{\al^2}{Cu'(y_c)},& 0< y_c\leq\frac{1}{\al},\\
C^{-1}\int_{y_c-\frac{1}{\al}}^{y_c}\dfrac{\al^2|y-y_c|^2}{u'(y_c)^2|y-y_c|^2}dy\geq \dfrac{\al}{Cu'(y_c)^2},& \frac{1}{\al}\leq y_c< 1,
\end{array}
\right.
\end{align*}
and
\begin{align*}
-\mathrm{II}_3(c)=&\int_0^1\frac{1}{(u(y)-c)^2}\Big(1-\frac{1}{\phi_{1}(y,c)^2}\Big)dy\\
\leq&C\int_0^1\frac{\al^2|y-y_c|^2}{(u(y)-c)^2}dy\le C\al^2\int_0^1\frac{1}{(y+y_c)^2}dy
\leq C\frac{\al^2}{u'(y_c)},
\end{align*}
or
\begin{align*}
-\mathrm{II}_3(c)=&\int_0^1\frac{1}{(u(y)-c)^2}\Big(1-\frac{1}{\phi_{1}(y,c)^2}\Big)dy\\
\leq&C\int_0^1\frac{\min\{\al^2|y-y_c|^2,1\}}{u'(y_c)^2|y-y_c|^2}dy
\leq C\frac{\al}{u'(y_c)^2}.
\end{align*}
This proves the first point of the lemma.

Direct calculation gives
\begin{align*}
\pa_c\mathrm{II}_3(c)=&\pa_c\int_0^1\frac{1}{(u(y)-c)^2}\Big(\frac{1}{\phi_{1}(y,c)^2}-1\Big)dy\\
=&\int_0^1\left(\frac{\partial_y}{u'(y_c)}+\partial_c\right)\Big(\frac{1}{(u(y)-c)^2}\Big(\frac{1}{\phi_{1}(y,c)^2}-1\Big)\Big)dy\\
&-\frac{1}{u'(y_c)(u(y)-c)^2}\Big(\frac{1}{\phi_{1}(y,c)^2}-1\Big)\Big|_{y=0}^1
\end{align*}
and
\begin{align*}
\pa_c^2\mathrm{II}_3(c)
=&\int_0^1\left(\frac{\partial_y}{u'(y_c)}+\partial_c\right)^2\Big(\frac{1}{(u(y)-c)^2}\Big(\frac{1}{\phi_{1}(y,c)^2}-1\Big)\Big)dy\\
&-\frac{1}{u'(y_c)}\left(\frac{\partial_y}{u'(y_c)}+\partial_c\right)
\Big(\frac{1}{(u(y)-c)^2}\Big(\frac{1}{\phi_{1}(y,c)^2}-1\Big)\Big)\Big|_{y=0}^1,\\
&-\partial_c\left(\frac{1}{u'(y_c)(u(y)-c)^2}\Big(\frac{1}{\phi_{1}(y,c)^2}-1\Big)\Big|_{y=0}^1\right).
\end{align*}
By Proposition \ref{prop:phi1} and Lemma \ref{lem:u}, we have
\beno
&&\left|\left(\frac{\partial_y}{u'(y_c)}+\partial_c\right)\Big(\frac{1}{(u(y)-c)^2}\Big(\frac{1}{\phi_{1}(y,c)^2}-1\Big)\Big)\right|\leq \frac{C}{u'(y_c)^2(u(y)-c)^2}\Big(1-\frac{1}{\phi_{1}(y,c)^2}\Big),\\
&&\left|\left(\frac{\partial_y}{u'(y_c)}+\partial_c\right)^2\Big(\frac{1}{(u(y)-c)^2}\Big(\frac{1}{\phi_{1}(y,c)^2}-1\Big)\Big)\right|\leq \frac{C}{u'(y_c)^4(u(y)-c)^2}\Big(1-\frac{1}{\phi_{1}(y,c)^2}\Big).
\eeno
On the other hand, we have for $y=0,1$,
\beno
\left|\frac{1}{u'(y_c)(u(y)-c)^2}\Big(\frac{1}{\phi_{1}(y,c)^2}-1\Big)\right|\leq C\f{\min\{\al^2|y-y_c|^2,1\}}{u'(y_c)(u(y)-c)^2}\leq C\f{\min\{\al^2,\al/u'(y_c)\}}{u'(y_c)|u(y)-c|},
\eeno
and for $y=0,1$,
\begin{align*}
&\left|\partial_c\left(\frac{1}{u'(y_c)(u(y)-c)^2}\Big(\frac{1}{\phi_{1}(y,c)^2}-1\Big)\right)\right|\\
&\leq
C\left(\frac{1}{u'(y_c)^3(u(y)-c)^2}+\frac{1}{u'(y_c)|u(y)-c|^3}\right)\Big|\frac{1}{\phi_{1}(y,c)^2}-1\Big|\\
&\quad+2\frac{1}{u'(y_c)(u(y)-c)^2}\frac{\cG(y,c)}{\phi_{1}(y,c)^2}\\
&\leq
C\left(\frac{\min\{\al^2|y-y_c|^2,1\}}{u'(y_c)^3(u(y)-c)^2}+\frac{\min\{\al^2|y-y_c|^2,1\}}{u'(y_c)|u(y)-c|^3}
+\frac{\al\min\{\al|y-y_c|,1\}}{u'(y_c)^2(u(y)-c)^2}\right)\\
&\leq
C\left(\frac{\min\{\al^2|y-y_c|^2,1\}}{u'(y_c)|u(y)-c|\rho(c)^2}+
\frac{\al\min\{\al|y-y_c|,1\}}{u'(y_c)^2(u(y)-c)^2}\right)\\
&\leq \f{C}{\rho(c)^2}\min\left\{\frac{\al^2}{u'(y_c)},\frac{\al}{u'(y_c)^2}\right\},
\end{align*}
here $\cG(y,c)=\f{\pa_c\phi_1(y,c)}{\phi_1(y,c)}$ and we used \eqref{eq:g}. Then the second inequality follows easily.
\end{proof}

\subsection{Estimate of $\rmA_1$}
Recall that
\begin{align*}
\rmA_1(c)&=\rho(c)u'(y_c)\pa_c\Big(p.v.\int_0^1\f{dy}{u(y)-c}\Big)=-\rho(c)u'(y_c)\partial_c\Big(\frac{1}{2\tc}H\big((v^{-1})'\chi_{I_v}\big)(\tc)\Big).
\end{align*}
Here we made the change of variable $\tc=v(y_c)\in I_v$.

\begin{lemma}\label{lem:A1}
It holds that for any $p\in (1,+\infty)$,
\begin{align*}
|\rmA_1(c)|\leq C\tc^2,\quad |\pa_c\rmA_1(c)|\leq C(\left|\ln(u(1)-c)\right|+1),\quad \left\|\rho\pa_c^2\rmA_1\right\|_{L_{\tc}^p}\leq C.
\end{align*}
\end{lemma}

\begin{proof}
Let $f_2(\tc)=\left((v^{-1})'(\tc)-\sqrt{\f{2}{u''(0)}}\right)\chi_{I_v}(\tc)$. Here $(v^{-1})'(0)=\f{1}{v'(0)}=\sqrt{\f{2}{u''(0)}}>0$. By Lemma \ref{lem:u},
we know that $f_2\in C^3([-v(1),0)\cup(0,v(1)])\cap C^2(I_v)$ and $zf_2(z)\in C^3(I_v)$. Then we have
\begin{align*}
\rmA_1&=-\rho(c)u'(y_c)\partial_c\Big(\frac{1}{2\tc}H(f_2)(\tc)
+\f{1}{\tc}\sqrt{\f{2}{u''(0)}}\ln \left|\f{v(1)-\tc}{v(1)+\tc}\right|\Big)\\
&=I_1+I_2,
\end{align*}
where
\begin{align*}
I_1&=-\rho(c) u'(y_c)\partial_{c}\Big(\frac{1}{2\tc}H(f_2)(\tc)\Big),\\
I_2&=-\sqrt{\f{1}{2u''(0)}}\tc(u(1)-c) u'(y_c)\partial_{\tc}\Big(\f{1}{\tc}\ln \left|\f{v(1)-\tc}{v(1)+\tc}\right|\Big).
\end{align*}

Let $\eta(\tc)\geq 0$ be a smooth even  function, such that $\eta(\tc)=1$ for $|\tc|\leq \f{v(1)}{4}$ and $\eta(\tc)=0$ for $|\tc|\geq \f{3v(1)}{4}$. Let $\eta_1(\tc)=\chi_{I_v}(\tc)-\eta(\tc)$. Using the fact that if $f(\tc)$ is odd, then $\pa_{\tc}\big(\f{f(\tc)}{\tc}\big)=\pa_{\tc}\Big(\int_{0}^{1}f'(t\tc)dt\Big)=\int_{0}^{1}tf''(t\tc)dt$ is odd, and $u'(y_c)=2\tc v'(y_c)$, we infer that
\ben\label{eq:I2-est1}
|\eta(\tc)I_2|\leq C\tc^3, \quad |\eta(\tc)\pa_cI_2|+|\tc\eta(\tc)\pa_c^2I_2|\leq C.
\een
On the other hand, we have
\begin{align*}
&\eta_1(\tc)(u(1)-c)\ln \left|\f{v(1)-\tc}{v(1)+\tc}\right|\leq C\eta_1(\tc),\\
&\eta_1(\tc)(u(1)-c)\left|\partial_{\tc}\Big(\f{1}{\tc}\ln \left|\f{v(1)-\tc}{v(1)+\tc}\right|\Big)\right|\leq C\eta_1(\tc),
\end{align*}
which along with \eqref{eq:I2-est1} give
\beno
|I_2|\leq |\eta(\tc)I_2|+|\eta_1(\tc)I_2|\leq C|\tc|^3.
\eeno
Using \eqref{eq:I2-est1} and the facts that
\begin{align*}
&\tc(u(1)-c) u'(y_c)\partial_{\tc}\Big(\f{1}{\tc}\ln \left|\f{v(1)-\tc}{v(1)+\tc}\right|\Big)\\
&=-2(u(1)-c) v'(y_c)\ln \left|\f{v(1)-\tc}{v(1)+\tc}\right|-2u'(y_c)v(1),
\end{align*}
and
\beno
&&\left|\eta_1(\tc)v'(y_c)\ln \left|\f{v(1)-\tc}{v(1)+\tc}\right|\right|\leq C|\ln (u(1)-c)|+C,\\
&&\left|\eta_1(\tc)\f{u(1)-c}{u'(y_c)}\ln \left|\f{v(1)-\tc}{v(1)+\tc}\right|\right|\leq C,
\eeno
we infer  that
\beno
|\pa_{c}I_2|\leq C(|\ln (u(1)-c)|+1),\quad
|\tc(u(1)-c)\pa_{c}^2I_2|\leq C.
\eeno

The estimates of $I_1$ are a direct consequence of Proposition \ref{prop:SIO}.
\end{proof}

\subsection{Estimate of $\rmA^2+\rmB^2$}
Recall that
\beno
\rmA(c)=\rmA_1(c)+u'(y_c)\rho(c)\mathrm{II}_3(c),\quad \rmB(c)=\pi\rho(c)\f{u''(y_c)}{u'(y_c)^2}.
\eeno

\begin{lemma}\label{lem:A2+B2}
It holds that
\begin{align*}
&C^{-1}(1+\al\rho_0(c))^2\leq \rmA(c)^2+\rmB(c)^2\leq C(1+\al\rho_0(c))^2,
\end{align*}
and
\begin{align*}
\left|\partial_c\Big(\f{1}{(\mathrm{A}^2+\mathrm{B}^2){u'(y_c)}}\Big)\right|\leq& \f{C(1+\al u'(y_c))}{(1+\al\rho_0)^3{u'(y_c)^3}}+\f{C(1+|\ln (1-y_c))|}{(1+\al\rho_0)^3{u'(y_c)}},\\
\left|\partial_c^2\Big(\f{1}{(\mathrm{A}^2+\mathrm{B}^2){u'(y_c)}}\Big)\right| \leq& \f{C(1+\al^2u'(y_c)^2)}{(1+\al\rho_0)^4{u'(y_c)^5}}
+\f{C(1+|\ln (1-y_c)|)^2}{(1+\al\rho_0)^4u'(y_c)}\\&
+\f{C(\al+|\rho\partial_c^2\mathrm{A_1}|)}{(1+\al\rho_0)^3{\rho}u'(y_c)}.
\end{align*}
\end{lemma}
\begin{proof}
By Lemma \ref{lem:II3} and Lemma \ref{lem:A1}, we get
\ben\label{eq:AB-up}
|\rmA(c)|\leq Cy_c^2+C(u(1)-c)\al y_c,\quad |\rmB(c)|\leq C(u(1)-c),
\een
which imply that
\beno
|\rmA(c)|+|\rmB(c)|\leq C(1+\al\rho_0(c)).
\eeno

Now we prove the lower bound of $\rmA^2+\rmB^2$.
Proposition \ref{prop:spectral} ensures that for any fixed $M>0$ and $\al\leq M$, $\rmA(c)^2+\rmB(c)^2\geq C^{-1}$. Thus, we may assume $\al\gg1$. If $\al y_c\le 1$,  then $|\rmB(c)|\geq C^{-1}$. Since $\rmA_1(u(1))=u(0)-u(1)<0$(by \ref{def:A1-Pi2}) and $\rmA_1(c)$ is continuous,  $\rmA_1(c)<-\f{u(0)-u(1)}{2}$ for $y_c\in [1-C/\al,1]$. Thus, by $\mathrm{II}_3<0$,  for $y_c\in [1-C/\al,1]$,
\begin{align*}
&|\rmA(c)|=-\rmA_1(c)-u'(y_c)\rho(c)\mathrm{II}_3>\f{u(0)-u(1)}{2}.
\end{align*}
For $y_c\in [0,1-C/\al]$ and $\al y_c\ge 1$, we get by Lemma \ref{lem:II3} and Lemma \ref{lem:A1} that
\begin{align*}
|\rmA(c)|\geq -u'(y_c)\rho(c)\mathrm{II}_3-|\rmA_1|\geq C_1^{-1}\al(1-y_c) y_c-C_2^{-1}y_c^2\geq \f12C_1^{-1}\al\rho_0(c).
\end{align*}
Thus, we deduce that $|\rmA(c)|+|\rmB(c)|>C^{-1}(1+\al \rho_0(c))$.

Notice that
\begin{align*}
u'(y_c)=2v(y_c)v'(y_c),\quad u''(y_c)=2v'(y_c)^2+2v(y_c)v'(y_c),
\end{align*}
which gives
\beno
\rho(c)\f{u''(y_c)}{u'(y_c)^2}=\f12\Big(u(1)-c+\f{(v(1)^2-\tc^2)\tc v''(y_c)}{v'(y_c)^2}\Big).
\eeno
Using the facts that  $\f{\pa_{\tc}}{2\tc}\Big(\f{1}{v'(y_c)^2}\Big)=\f{-v''(y_c)}{\tc v'(y_c)^3}$, $|v''(y_c)|\leq C|y_c|$, and
\begin{align*}
&\big|\pa_{\tc}\big((v(1)^2-\tc^2)\tc v''(y_c)\big)\big|\\
&\leq \big|(v(1)^2-\tc^2)v''(y_c)\big|+\big|2\tc^2 v''(y_c)\big|+\big|(v(1)^2-\tc^2)\tc v'''(y_c)/v'(y_c)\big|\leq C|\tc|,
\end{align*}
we infer that $\big|\pa_c\rmB(c)\big|\le C.$ Similarly, $|\pa_c^2\rmB|\leq \f{C}{y_c^2}$.

A direct calculation gives
\begin{align*}
&\left|\partial_c\Big(\f{1}{(\mathrm{A}^2+\mathrm{B}^2){u'(y_c)}}\Big)\right|\\
&\leq C\left|\f{\rmA\pa_c\rmA+\rmB\pa_c\rmB}{(\mathrm{A}^2+\mathrm{B}^2)^2u'(y_c)}\right|
+\f{C}{(\mathrm{A}^2+\mathrm{B}^2)u'(y_c)^3}\\
&\leq \f{C(|\pa_c\rmA_1|+|\pa_c\rmB|+|\pa_c(u'(y_c)\rho\mathrm{II}_3)|)}{(\mathrm{A}^2+\mathrm{B}^2)^{3/2}u'(y_c)}
+\f{C}{(\mathrm{A}^2+\mathrm{B}^2)u'(y_c)^3},
\end{align*}
and
\begin{align*}
&\left|\partial_c^2\Big(\f{1}{(\mathrm{A}^2+\mathrm{B}^2){u'(y_c)}}\Big)\right|\\
&\leq \f{C}{(\mathrm{A}^2+\mathrm{B}^2){u'(y_c)^5}}
+\f{C(|\mathrm{A}\partial_c\mathrm{A}|+|\mathrm{B}\partial_c\mathrm{B}|)}{(\mathrm{A}^2+\mathrm{B}^2)^2{u'(y_c)^3}}
+\f{C(|\partial_c\mathrm{A}|^2+|\partial_c\mathrm{B}|^2)}{(\mathrm{A}^2+\mathrm{B}^2)^2{u'(y_c)}}\\
&\quad+\f{C(|\mathrm{A}\partial_c^2\mathrm{A}|+|\mathrm{B}\partial_c^2\mathrm{B}|)}{(\mathrm{A}^2+\mathrm{B}^2)^2{u'(y_c)}}
+\f{C(|\mathrm{A}\partial_c\mathrm{A}|+|\mathrm{B}\partial_c\mathrm{B}|)^2}{(\mathrm{A}^2+\mathrm{B}^2)^3{u'(y_c)}}\\
&\leq \f{C}{(\mathrm{A}^2+\mathrm{B}^2){u'(y_c)^5}}+
\f{C(|\partial_c\mathrm{A}|^2+|\partial_c\mathrm{B}|^2)}{(\mathrm{A}^2+\mathrm{B}^2)^2{u'(y_c)}}
+\f{C(|\mathrm{A}\partial_c^2\mathrm{A}|+|\mathrm{B}\partial_c^2\mathrm{B}|)}{(\mathrm{A}^2+\mathrm{B}^2)^2{u'(y_c)}}\\
&\leq \f{C}{(\mathrm{A}^2+\mathrm{B}^2){u'(y_c)^5}}+
\f{C|\partial_c\mathrm{A}|^2}{(\mathrm{A}^2+\mathrm{B}^2)^2{u'(y_c)}}
+\f{C|\mathrm{A}\partial_c^2\mathrm{A}|}{(\mathrm{A}^2+\mathrm{B}^2)^2{u'(y_c)}}.
\end{align*}

Thus, by Lemma \ref{lem:II3} and Lemma \ref{lem:A1}, we obtain
\begin{align*}
&\left|\partial_c\Big(\f{1}{(\mathrm{A}^2+\mathrm{B}^2){u'(y_c)}}\Big)\right|\\
&\leq \f{C\big((\mathrm{A}^2+\mathrm{B}^2)^{1/2}+u'(y_c)+\al u'(y_c)\big)}{(\mathrm{A}^2+\mathrm{B}^2)^{3/2}u'(y_c)^3}
+\f{C|\pa_c\rmA_1|}{(\mathrm{A}^2+\mathrm{B}^2)^{3/2}u'(y_c)}\\
&\leq \f{C(1+\al u'(y_c))}{(1+\al\rho_0)^3{u'(y_c)^3}}+\f{C(1+|\ln (u(1)-c)|)}{(1+\al\rho_0)^3{u'(y_c)}}.
\end{align*}
and using the facts that
\begin{align*}
&|\pa_c\rmA|^2\leq C(|\partial_c\mathrm{A_1}|^2+|\partial_c(u'(y_c)\rho\mathrm{II}_3)|^2),\\
&|\pa_c^2\rmA|\leq C(|\partial_c^2\mathrm{A_1}|+|\partial_c^2(u'(y_c)\rho\mathrm{II}_3)|),
\end{align*}
we deduce that
\begin{align*}
\partial_c^2\Big(\f{1}{(\mathrm{A}^2+\mathrm{B}^2){u'(y_c)}}\Big)\leq& \f{C(1+\al^2u'(y_c)^2)}{(1+\al\rho_0)^4{u'(y_c)^5}}
+\f{C\al}{(1+\al\rho_0)^3{\rho}u'(y_c)}\\
&+\f{C(1+|\ln (u(1)-c)|)^2}{(1+\al\rho_0)^4u'(y_c)}
+\f{C|\rho\partial_c^2\mathrm{A_1}|}{(1+\al\rho_0)^3{\rho}u'(y_c)}.
\end{align*}
This completes the proof of the lemma.
\end{proof}

\subsection{Estimate of $\rmA_2^2+\rmB_2^2$ and $J$}
Recall that
\beno
\rmA_2(c)=(u(0)-c)\rmA(c)+J(c),\quad \rmB_2(c)=(u(0)-c)\rmB(c),\quad J(c)=\f{u'(y_c)(u(1)-c)}{\phi_1(0,c)\phi_1'(0,c)}.
\eeno

\begin{lemma}\label{lem:J}
It holds that
\beno
|J(c)|\leq \f{C(1-y_c)}{\phi_1(0,c)\al^2},\quad
|\pa_cJ(c)|\leq \f{C}{\phi_1(0,c)\al^2y_c^2}\quad \textrm{and}\quad
|\pa_c^2J(c)|\leq \f{C}{\phi_1(0,c)\al^2y_c^4}.
\eeno
\end{lemma}
\begin{proof}
We get by Lemma \ref{lem:FG} that
\beno
|J(c)|\leq \f{Cy_c(1-y_c)}{\phi_1(0,c)^2\al \min\{\al y_c,1\}}\leq \f{C(1-y_c)}{\phi_1(0,c)^2\al^2}.
\eeno
Recall that  $\cF=\frac{\partial_y\phi_1}{\phi_1}$ and $\cG=\frac{\partial_c\phi_1}{\phi_1}$ defined by \ref{def:FG}. A direct calculation gives
\begin{align*}
&\f{1}{\phi_1\pa_y\phi_1}(0,c)=\f{1}{\phi_1(0,c)^2\cF(0,c)},\\
&\pa_c\Big(\f{1}{\phi_1\pa_y\phi_1}\Big)(0,c)=-\f{\pa_c\cF}{\phi_1^2\cF^2}(0,c)-\f{2\cG}{\phi_1^2\cF}(0,c),\\
&\pa_c^2\Big(\f{1}{\phi_1\pa_y\phi_1}\Big)(0,c)
=\f{-2\partial_c\cG+4\cG^2}{\phi_1^2\cF}(0,c)
+\f{4\cG\partial_c\cF}{\phi_1^2\cF^2}(0,c)
-\f{\partial_c^2 \cF}{\phi_1^2\cF^2}(0,c)
+\f{2|\partial_c \cF|^2}{\phi_1^2\cF^3}(0,c),
\end{align*}
from which and Lemma \ref{lem:FG}, we infer that
\begin{align*}
&\left|\f{1}{\phi_1(0,c)^2\cF(0,c)}\right|\sim \f{1+\al y_c}{\al^2\phi_1(0,c)^2 y_c},\\
&\left|\partial_c\Big(\f{1}{\phi_1\phi_1'}\Big)(0,c)\right|
\leq \f{C(1+\al y_c)^2}{\al^2\phi_1(0,c)^2y_c^3},\\
&\left|\partial_c^2\Big(\f{1}{\phi_1\phi_1'}\Big)(0,c)\right|
\leq \f{C(1+\al y_c)^3}{\al^2\phi_1(0,c)^2y_c^5}.
\end{align*}
Therefore, we obtain
\begin{align*}
&|\partial_c J(c)|
\leq \f{C}{|\phi_1\phi_1'(0,c)|u'(y_c)}
+{C}{u'(y_c)}\left|\partial_c\f{1}{\phi_1\phi_1'(0,c)}\right|
\leq \f{C}{\phi_1(0,c)\al^2y_c^2},\\
&|\partial_c^2 J(c)|
\leq C\sum_{k=0}^2{u'(y_c)^{2k-3}}\left|\partial_c^k\f{1}{\phi_1\phi_1'(0,c)}\right|
\leq \f{C}{\phi_1(0,c)\al^2y_c^4}.
\end{align*}
This proves the lemma.
\end{proof}

\begin{lemma}\label{lem:A2+B2-2}
It holds that
\beno
C^{-1}\f{(1+\al\rho_0)^2(1+\al y_c)^4}{\al^4}\leq \rmA_2(c)^2+\rmB_2(c)^2\leq C\f{(1+\al\rho_0)^2(1+\al y_c)^4}{\al^4},
\eeno
and
\begin{align*}
\left|\partial_c\Big(\f{1}{(\rmA_2^2+\rmB_2^2){u'(y_c)}}\Big)\right|
\leq& \f{C\al^4}{(1+\al\rho_0)^3(1+\al y_c)^3{u'(y_c)^3}}\\
&+\f{C\al^4(1+|\ln (1-y_c)|)}{(1+\al\rho_0)^3(1+\al y_c)^4u'(y_c)},\\
\left|\partial_c^2\Big(\f{1}{(\rmA_2^2+\rmB_2^2){u'(y_c)}}\Big)\right| \leq& \f{C\al^4}{(1+\al\rho_0)^4(1+\al y_c)^2{u'(y_c)^5}}\\
&+\f{C\al^4(|\ln(1-y_c)|^2+|\rho\pa_c^2\rmA_1|+\al^2)}{(1+\al\rho_0)^3(1+\al y_c)^4{\rho}u'(y_c)}.
\end{align*}

\end{lemma}

\begin{proof}
We get by Lemma \ref{lem:J} and \eqref{eq:AB-up} that
\begin{align*}
|J(c)|\leq \f{C}{\al^2},\quad |(u(0)-c)\rmA(c)|\leq Cy_c^2(1+\al \rho_0),\quad |(u(0)-c)\rmB(c)|\leq Cy_c^2,
\end{align*}
which gives
\beno
|\rmA_2(c)|+|\rmB_2(c)|\leq \f{C}{\al^2}+Cy_c^2(1+\al \rho_0)\leq C\frac{(1+\al \rho_0)(1+\al y_c)^2}{\al^2}.
\eeno
Proposition \ref{prop:spectral} ensures that for some fixed $M$, $|\rmA_2(c)|+|\rmB_2(c)|\ge C^{-1}$ for $\al \le M$.
Thus, we may assume that $\al\gg1$. Let $\delta_1$ be a small constant.
For $\ 0\leq y_c\leq \f{\delta_1}{\al}$ , we have
\beno
|J(c)|\geq \f{y_c}{C\al\min\{\al y_c,1\}}\geq \f{1}{C\al^2},
\eeno
which gives for $0\leq y_c\leq \dfrac{\delta_1}{\al}$,
\begin{align*}
|\rmA_2(c)| &\geq |J(c)|-|(u(0)-c)\rmA(c)| \\
&\geq \f{1}{C\al^2}-Cy_c^2(y_c+\al\rho_0)\geq \f{1}{C\al^2}\geq \f{(1+\al \rho_0)(1+\al y_c)^2}{C\al^2}.
\end{align*}
For $\f{\delta_1}{\al}\leq y_c\leq \dfrac{1}{\delta_1\al}$ we have $|\rmB(c)|\geq |\rmB(0)/2|\geq 1/C,$ hence,
$$|\rmB_2(c)|\geq \f{|u(0)-c|}{C}\geq \f{y_c^2}{C}\geq\frac{(1+\al \rho_0)(1+\al y_c)^2}{C\al^2}.$$
For $\f{1}{\delta_1\al}\leq y_c\leq 1$, we have
\begin{align*}
|\rmA_2(c)|&\geq |(u(0)-c)\rmA(c)|-|J(c)|\\
&\geq \f{y_c^2( y_c+\al\rho_0)}{C}-\f{C}{\al^2}\\
&\geq \f{y_c^2( y_c+\al\rho_0)}{C}\geq\frac{(1+\al \rho_0)(1+\al y_c)^2}{C\al^2}.
\end{align*}
Here we used the fact that $|\rmA(c)|\ge C^{-1}(1+\al\rho_0(c))$ for $\al y_c\ge 1$(see the proof of Lemma \ref{lem:A2+B2}).
Summing up, we conclude the lower bound of $|\rmA_2(c)|+|\rmB_2(c)|$.

Using Lemma \ref{lem:A1}, Lemma \ref{lem:II3} and Lemma \ref{lem:J}, we deduce that
\begin{align*}
|\pa_{c}\rmA_2(c)|&\leq |\pa_{c}((u(0)-c)\mathrm{A}_1)|+|\pa_{c}((u(0)-c)u'(y_c)\rho\mathrm{II}_3)|+|\pa_{c}J|\\
&\leq Cy_c^2(1+|\ln(1-y_c)|)+C\al y_c+C(\al^2y_c^2)^{-1}\\
&\leq Cy_c^2(1+|\ln(1-y_c)|)+C(1+\al^3 y_c^3) (\al^2y_c^2)^{-1},\\
|\pa_{c}^2\rmA_2(c)|&\leq |\pa_{c}^2((u(0)-c)\rmA_1)|+|\pa_{c}^2((u(0)-c)u'(y_c)\rho\mathrm{II}_3)|+|\pa_{c}^2J|\\
&\leq C(1-y_c)^{-1}|\rho\pa_c^2\rmA_1|+|\pa_c\rmA_1|+C\al \rho_0^{-1}+C(\al^2y_c^4)^{-1}\\
&\leq C(1-y_c)^{-1}|\rho\pa_c^2\rmA_1|+C\al \rho_0^{-1}+C(\al^2y_c^4)^{-1}.
\end{align*}
We known from the proof of Lemma \ref{lem:A2+B2} that
\beno
|\mathrm{B_2}(c)|\leq Cy_c^2,\quad
|\pa_{c}\mathrm{B_2}(c)|\leq C,\quad
|\pa_{c}^2\mathrm{B_2}(c)|\leq C.
\eeno

With these estimates, we can deduce that
\begin{align*}
\left|\partial_c\Big(\f{1}{(\rmA_2^2+\rmB_2^2){u'(y_c)}}\Big)\right|
\leq& \f{C\al^4}{(1+\al\rho_0)^3(1+\al y_c)^3{u'(y_c)^3}}\\
&+\f{C\al^4(1+|\ln(1-y_c)|)}{(1+\al\rho_0)^3(1+\al y_c)^4u'(y_c)},
\end{align*}
and
\begin{align*}
&\left|\pa_c^2\Big(\f{1}{(\rmA_2^2+\rmB_2^2){u'(y_c)}}\Big)\right|\\
&\leq \f{C}{(\rmA_2^2+\rmB_2^2){u'(y_c)^5}}+
\f{C|\partial_c\rmA_2|^2}{(\rmA_2^2+\rmB_2^2)^2{u'(y_c)}}
+\f{C|\rmA_2\pa_c^2\rmA_2|}{(\rmA_2^2+\rmB_2^2)^2{u'(y_c)}}\\
&\leq \f{C\al^4}{(1+\al\rho_0)^2(1+\al y_c)^4{u'(y_c)^5}}+
\f{C\al^6(y_c^4(1+|\ln(1-y_c)|)^2+(1+\al y_c)^6/(\al y_c)^4)}{(1+\al\rho_0)^4(1+\al y_c)^8{u'(y_c)}}\\
&\quad+\f{C\al^6(|\rho\pa_c^2\rmA_1|/(1-y_c)+\al /\rho_0+1/(\al^2y_c^4))}{(1+\al\rho_0)^3(1+\al y_c)^6{u'(y_c)}}\\
&\leq \f{C\al^4}{(1+\al\rho_0)^4(1+\al y_c)^2{u'(y_c)^5}}
+\f{C\al^4((1+|\ln(1-y_c)|)^2+|\rho\pa_c^2\rmA_1|+\al^2)}{(1+\al\rho_0)^3(1+\al y_c)^4{\rho}u'(y_c)}.
\end{align*}
This proves the lemma.
\end{proof}

\section{Estimates of basic integral operators}

In this section, we present some estimates for some basic integral operators appeared in $K_{o}$ and $K_e$.

Let $I_0=(0,1)$ and $I=(-1,1)$. Let $I_v=(-v(1), v(1))$ with $v$ given by \eqref{def:v} and $\rho_0(c)=\rho(c)/u'(y_c)\sim y_c(1-y_c)$. We denote by $\|\cdot\|_{L^p_{\tc}}$ the norm of $L^p(I_v, d\tc)$, and $C$ a constant independent of $\al$, which may be different from line to line.

\subsection{Estimates of the operator $\Lambda_{j,1}$ for $j=1,2$}
Recall that
\begin{align*}
\textrm{II}_{1,1}(\varphi)(c)
&=p.v.\int_0^1\f{\Int(\varphi)(y)-\Int(\varphi)(y_c)}{(u(y)-u(y_c))^2}\,dy\\
&=\f{1}{2\tc}\partial_{\tc}\Big(\frac{1}{2\widetilde{c}}p.v.\int_{-1}^1\frac{\Int(\varphi)(y)}{v(y)-\widetilde{c}}\,dy\Big)
\\&\quad-\f{\Int(\varphi)\big(v^{-1}(\tc)\big)}{2\tc}\pa_{\tc}\Big(\frac{1}{2\widetilde{c}}p.v.\int_{-v(1)}^{v(1)}\f{(v^{-1})'(z)}{z-\tc}\,dz\Big).
\end{align*}
Here  $\tc=v(y_c)$ and
\begin{align*}
&\Int(\varphi)(y)=\int_{0}^y\varphi(y')\,dy'\quad\textrm{for}\,\, y\in [0,1], \quad \Int(\varphi)(y)=\Int(\varphi)(-y)\quad\textrm{for}\,\, y\in [-1,0].
\end{align*}
Thus, $\Int(\varphi)(y)$ is an even function with $\Int(\varphi)(0)=0$ and satisfies
\begin{align}\label{eq:Int-Lp}
&\|\Int(\varphi)\|_{W^{1,p}(I)}\leq C\|\varphi\|_{L^p(I_0)},\quad \Big\|\f{1}{y}\Int(\varphi)(y)\Big\|_{L^p(I)}\leq C\|\varphi\|_{L^p(I_0)}.
\end{align}
Let us also recall that
\begin{align*}
&\Lambda_{1,1}(\varphi)(c)=\rmA_1(c)\varphi(y_c)+\rho(c)u''(y_c)\mathrm{II}_{1,1}(\varphi)(c),\\
&\Lambda_{2,1}(\varphi)(c)=\rmA_1(c)\varphi(y_c)+\rho(c)\mathrm{II}_{1,1}(u''\varphi)(c).
\end{align*}

\begin{lemma}\label{lem:Lambda-11}
For $k=0,1,2$ and $j=1,2$, we have
\beno
&&\left\|\pa_{\tc}^k\Lambda_{j,1}(\varphi)\right\|_{L^2_{\tc}}\leq C\|\varphi\|_{H^k(I_0)},\\
&&\left\|\tc^k\pa_c^k\big(\Lambda_{j,1}(\varphi)\big)\right\|_{L^2_{\tc}}\leq C\|\varphi\|_{H^k(I_0)}.
\eeno
If $\varphi(0)=0$, then we have
\beno
&&\left\|\tc^{-1}\Lambda_{j,1}(\varphi)\right\|_{L^2_{\tc}}\leq C\|\varphi\|_{H^1(I_0)},\\
&&\left\|\tc^{-2}\Lambda_{j,1}(\varphi)\right\|_{L^2_{\tc}}+\left\|\pa_c\big(\Lambda_{j,1}(\varphi)\big)\right\|_{L^2_{\tc}}\leq C\|\varphi\|_{H^2(I_0)}.
\eeno
If $\varphi(1)=0$, then we have
\beno
&&\left\|(u(1)-c)^{-1}\Lambda_{j,1}(\varphi)\right\|_{L^2_{\tc}}\leq C\|\varphi\|_{H^1(I_0)}.
\eeno
\end{lemma}

\begin{remark}\label{rem:Lambda}
Using Sobolev embedding and Lemma \ref{lem:Lambda-11}, we can deduce that
\begin{align*}
\left\|\Lambda_{j,1}(\varphi)\right\|_{L^{\infty}_{\tc}}&\leq C\|\varphi\|_{H^1(I_0)},
\end{align*}
and
\begin{align*}
\left\|\tc\pa_c\big(\Lambda_{j,1}(\varphi)\big)\right\|_{L^{\infty}_{\tc}}&\leq C\big(\left\|\pa_{\tc}^2\big(\Lambda_{j,1}(\varphi)\big)\right\|_{L^2_{\tc}}^{\f12}
+\left\|\pa_{\tc}\big(\Lambda_{j,1}(\varphi)\big)\right\|_{L^2_{\tc}}^{\f12}\big)\left\|\pa_{\tc}\big(\Lambda_{j,1}(\varphi)\big)\right\|_{L^2_{\tc}}^{\f12}\\
&\leq C\|\varphi\|_{H^2(I_0)}^{\f12}\|\varphi\|_{H^1(I_0)}^{\f12},
\end{align*}
and if $\varphi(0)=0$,
\beno
\left\|(u'(y_c))^{-1}\Lambda_{j,1}(\varphi)\right\|_{L^{\infty}_{\tc}}\leq C\|\varphi\|_{H^2(I_0)}.
\eeno
If $\varphi(1)=0$, due to  $\Lambda_{j,1}(\varphi)(1)=0$, we have
\beno
\Big\|\f{\Lambda_{j,1}(\varphi)}{1-y_c}\Big\|_{L^{\infty}_{\tc}}\leq C\|\Lambda_{j,1}(\varphi)\|_{W^{1,\infty}}\leq C\|\Lambda_{j,1}(\varphi)\|_{H^1}^{\f12}\|\Lambda_{j,1}(\varphi)\|_{H^2}^{\f12}
\leq C\|\varphi\|_{H^1}^{\f12}\|\varphi\|_{H^2}^{\f12}.
\eeno
\end{remark}

\begin{proof}
Let us recall that for $y_c\in [0,1]$,
\begin{align*}
&\rmA_1=-\rho u'(y_c)\partial_c\Big(\frac{1}{2\tc}H\big((v^{-1})'\chi_{[v(-1),v(1)]}\big)(\tc)\Big)
=-\rho u'(y_c)\partial_c\Big(\frac{1}{2\tc}H\big(f_2\big)(\tc)\Big),\\&{\mathrm{II}_{1,1}}(\varphi)=-\pa_{c}\Big(\f{1}{2\tc}
H\big(f_1f_2\big)\Big)+
f_1\pa_{c}\Big(\f{1}{2\tc}H\big(f_2\big)\Big),
\end{align*}
and by the proof of Proposition \ref{prop:SIO}, we find that
\begin{align*}
&\Lambda_{1,1}(\varphi)(c)=\rmA_1\varphi(y_c)+\rho(c)u''(y_c)\mathrm{II}_{1,1}(\varphi)(c)\\
&=\rho(- u'(y_c)\varphi(y_c)+u''(y_c)f_1(\tc))\partial_c\Big(\frac{1}{2\tc}H\big(f_2\big)(\tc)\Big)-\rho(c)u''(y_c)\pa_{c}\Big(\f{1}{2\tc}
H\big(f_1f_2\big)\Big)\\
&=(- u'(y_c)\varphi(y_c)+u''(y_c)f_1(\tc))\left(\f{u(1)-c}{4\tc} H\left(zZ(f_2)\right)(\tc)+ \f{f_2(v(1))v(1)}{2}\right)\\
&\ \ \ -u''(y_c)\left(\f{u(1)-c}{4} H\left(Z(f_1f_2)\right)(\tc)+ \f{f_1f_2(v(1))v(1)}{2}\right)\\
&=\f{v(1)}{2v'(1)}u''(y_c)\int_{1}^{y_c}\varphi(z)dz-\f{v(1)}{2v'(1)}u'(y_c)\varphi(y_c)\\
&\quad+\f{u(1)-c}{4}\left(-u''(y_c)H(Z(f_1f_2))(\tc)+\f{u''(y_c)f_1(\tc)-\varphi(y_c)u'(y_c)}{\tc}H(zZ(f_2))\right)\\
&=\Lambda_{1,1,1}(\varphi)+\Lambda_{1,1,2}(\varphi),
\end{align*}
where $f_1(\tc)=\big(\Int(\varphi)\circ v^{-1}\big)\chi_{[-v(1),v(1)]}(\tc)$ and $f_2(\tc)=(v^{-1})'\cdot\chi_{[-v(1),v(1)]}(\tc)$. Let
\beno
&&\varphi_1(\tc)=Z(f_1f_2)(\tc),\quad \varphi_2(\tc)=\tc Z(f_2)(\tc),\\
&&\varphi_3(\tc)=\f{u''(y_c)f_1(\tc)-\varphi(y_c)u'(y_c)}{\tc},\quad\varphi_4(\tc)=-u''(y_c).
\eeno
Then $\varphi_1,\varphi_3 $ are odd in $\tc$, and  $\varphi_2,\varphi_4 $ are even and in $ H_{\tc}^2(I_v)$, thus $\Lambda_{1,1,2}(\varphi)$ is even. If $\varphi\in H^k(I_0),k=1,2, $ then by Lemma \ref{lem:Z}, $\varphi_1\in H^k(I_v)$, and $\varphi_3\in H^k(I_v\setminus\{0\})$. As $\partial_{y_c}(u''(y_c)f_1(\tc)-\varphi(y_c)u'(y_c))=u'''(y_c)f_1(\tc)-\varphi'(y_c)u'(y_c)$, we find that $\partial_{\tc}(u''(y_c)f_1(\tc)-\varphi(y_c)u'(y_c))\big|_{\tc=0}=(u''(y_c)f_1(\tc)-\varphi(y_c)u'(y_c))\big|_{\tc=0}=0 $ and $\varphi_3\big|_{\tc=0}=0$, thus $\varphi_3\in H^k(I_v).$

Let us verify that
\ben
\varphi_1(\tc)\varphi_4(\tc)+\varphi_3(\tc)\varphi_2(\tc)\big|_{\tc=\pm v(1)}=0.\label{eq:phi-relation}
\een
We only need to check the case of $\tc= v(1)$. Using the facts that $(v^{-1})'(v(1))=\f{1}{v'(1)}$, $(v^{-1})''(v(1))=\f{-v''(1)}{v'(1)^3}$, $u'(1)=2v(1)v'(1)$ and $u''(1)=2v'(1)^2+2v(1)v''(1)$, we deduce that
\begin{align*}
\varphi_1(v(1))&=(f_1f_2)'-\f{f_1f_2}{v(1)}\\
&=\f{\varphi(1)}{v'(1)^2}-f_1\Big(\f{v''(1)}{v'(1)^3}+\f{1}{v(1)v'(1)}\Big)\\
&=\f{\varphi(1)}{v'(1)^2}-f_1\f{u''(1)}{v'(1)^2u'(1)}
=-\f{v(1)\varphi_3(v(1))}{v'(1)^2u'(1)}\\
\varphi_2(v(1))&=-v(1)\f{v''(1)}{v'(1)^3}-\f{1}{v'(1)}\\
&=-\f{u''(1)}{2v'(1)^3}=\f{\varphi_4(v(1))}{2v'(1)^3}
=\f{v(1)\varphi_4(v(1))}{v'(1)^2u'(1)}.
\end{align*}
This implies \eqref{eq:phi-relation}. Then by Lemma \ref{lem:boundaryvanish} and \eqref{eq:aver}, \eqref{eq:H-p-est}, we have $(u(1)-c)H(zZ(f_2))\in L^{\infty}, $ and for $k=0,1,2,$
\begin{align*}
&\left\|\pa_{\tc}^k\Lambda_{1,1,2}(\varphi)\right\|_{L^{2}_{y_c}}\leq C\|Z(f_1f_2)\|_{H^{k}(I_v)}+C\|\varphi\|_{H^{k}(I_0)}
\leq C\|\varphi\|_{H^{k}(I_0)}.
\end{align*}
This shows that
\beno
&&\left\|\pa_{\tc}^k\big(\Lambda_{1,1}(\varphi)\big)\right\|_{L^2_{\tc}}\leq C\|\varphi\|_{H^{k}(I_0)},\\
&&\left\|\tc^k\pa_c^k\big(\Lambda_{1,1}(\varphi)\big)\right\|_{L^{2}_{\tc}}\leq C\|\varphi\|_{H^{k}(I_0)},
\eeno here we used $\pa_{\tc}\big(\Lambda_{1,1}(\varphi)\big)\big|_{\tc=0}=0 $ for $\varphi\in{H^{2}(I_0)}, $ which follows from the facts that $\Lambda_{1,1,2}(\varphi) $ is an even function in $H^{2}_{\tc}(-v(1)/2,v(1)/2)$ and
$$\pa_{\tc}\big(\Lambda_{1,1,1}(\varphi)\big)=\f{v(1)}{2v'(1)v'(y_c)}\left(u'''(y_c)\int_{1}^{y_c}\varphi(z)dz-u'(y_c)\varphi'(y_c)\right). $$

Next, we consider $\Lambda_{2,1}(\varphi)$.
Let $\widetilde{f_1}(\tc)=\big(\Int(u''\varphi)\circ v^{-1}\big)\chi_{[-v(1),v(1)]}(\tc)$. Then we find that
\begin{align*}
&\Lambda_{2,1}(\varphi)(y_c)=\rmA_1\varphi(y_c)+\rho(c)\mathrm{II}_{1,1}(u''\varphi)(c)\\
&=\rho(- u'(y_c)\varphi(y_c)+\widetilde{f_1}(\tc))\partial_c\Big(\frac{1}{2\tc}H\big(f_2\big)(\tc)\Big)-\rho(c)\pa_{c}\Big(\f{1}{2\tc}
H\big(\widetilde{f_1}f_2\big)\Big)\\
&=\f{v(1)}{2v'(1)}\int_{1}^{y_c}u''(z)\varphi(z)dz-\f{v(1)}{2v'(1)}u'(y_c)\varphi(y_c)\\
&\quad+\f{u(1)-c}{4}\left(-H(Z(\widetilde{f_1}\widetilde{f_2}))(\tc)+\f{\widetilde{f_1}(\tc)-\varphi(y_c)u'(y_c)}{\tc}H(zZ(\widetilde{f_2}))\right).
\end{align*}
Let $\varphi_1(\tc)=Z(\widetilde{f_1}{f_2})(\tc)$, $\varphi_2(\tc)=\tc Z({f_2})(\tc)$, $\varphi_3(\tc)=\f{\widetilde{f_1}(\tc)-\varphi(y_c)u'(y_c)}{\tc}$ and $\varphi_4(\tc)=-\chi_{[-v(1),v(1)]}(\tc)$. We also have
\beno
\varphi_1(\tc)\varphi_4(\tc)+\varphi_3(\tc)\varphi_2(\tc)\big|_{\tc=\pm v(1)}=0.
\eeno
Then as above, we can deduce from Lemma \ref{lem:boundaryvanish} that
\beno
&&\left\|\pa_{\tc}^k\big(\Lambda_{2,1}(\varphi)\big)\right\|_{L^2_{\tc}}\leq C\|\varphi\|_{H^k(I_0)},\\
&&\left\|\tc^k\pa_c^k\big(\Lambda_{2,1}(\varphi)\big)\right\|_{L^{2}_{\tc}}\leq C\|\varphi\|_{H^k(I_0)},
\eeno
and $\pa_{\tc}\big(\Lambda_{2,1}(\varphi)\big)\big|_{\tc=0}=0 $ for $\varphi\in{H^{2}(I_0)}.$

If $\varphi(0)=0$ and $\varphi\in H^{1}(I_0)$, we get by Remark \ref{rem:Hg} and Hardy's inequality that
\beq\label{eq:estL^p1}
\Lambda_{j,1}(\varphi)\big|_{\tc=0}=0,\ \ \ \ \left\|{\tc}^{-1}\Lambda_{j,1}(\varphi)\right\|_{L^2_{\tc}}\leq C\left\|\pa_{\tc}\Lambda_{j,1}(\varphi)\right\|_{L^2_{\tc}}\leq C\|\varphi\|_{H^{1}(I_0)}.
\eeq

If $\varphi(1)=0$, due to $\Lambda_{j,1}(\varphi)\big|_{y_c=1}=0$, we get by Hardy's inequality and \eqref{eq:estL^p1} that
\begin{align*}
&\left\|(u(1)-c)^{-1}\Lambda_{j,1}(\varphi)\right\|_{L^2_{y_c}}\leq C\left\|\pa_{\tc}\Lambda_{j,1}(\varphi)\right\|_{L^2_{\tc}}\leq C\|\varphi\|_{H^{1}(I_0)}.
\end{align*}

If $\varphi(0)=0$ and $\varphi\in H^2(I_0)$, then $\Lambda_{j,1}(\varphi)\big|_{\tc=0}=\pa_{\tc}\big(\Lambda_{2,1}(\varphi)\big)\big|_{\tc=0}=0 $, thus Hardy's inequality gives
\begin{align}
\left\|{\tc^{-2}}\Lambda_{j,1}(\varphi)\right\|_{L^2_{y_c}}+\left\|\pa_c\big(\Lambda_{j,1}(\varphi)\big)\right\|_{L^2_{y_c}}\leq C\left\|\pa_{\tc}^2\Lambda_{j,1}(\varphi)\right\|_{L^2_{\tc}}\leq C\|\varphi\|_{H^{2}(I_0)}.\nonumber
\end{align}
This gives our result by noting $\|\cdot\|_{L^2_{y_c}}\sim \|\cdot\|_{L^2_{\tc}}$ due to $v'(y)\ge c_1$.
\end{proof}

\subsection{Estimate of the operator $\mathrm{II}_{1,2}$}
Recall that
\beno
\mathrm{II}_{1,2}(\varphi)(c)=\int_0^1\int_{y_c}^{z}{\varphi}(y)
\left(\f{1}{(u(z)-c)^2}\left(\f{\phi_1(y,c)\,}{\phi_1(z,c)^2}-1\right)\right)\,dydz.
\eeno
We introduce for $k=0,1,2$,
\begin{align*}
\mathcal{L}_k(\varphi)(c)=\int_0^1\int_{y_c}^{z}{\varphi}(y)\left(\frac{\partial_z+\partial_y}{u'(y_c)}+\partial_c\right)^k
\left(\f{1}{(u(z)-c)^2}\left(\f{\phi_1(y,c)\,}{\phi_1(z,c)^2}-1\right)\right)\,dydz,
\end{align*}
and for $k,j=0,1$,
\begin{align*}
I_{k,j}(\varphi)(c)=\int_{y_c}^{z}{\varphi}(y)\left(\frac{\partial_z+\partial_y}{u'(y_c)}+\partial_c\right)^k
\left(\f{1}{(u(z)-c)^2}\left(\f{\phi_1(y,c)\,}{\phi_1(z,c)^2}-1\right)\right)\,dy\bigg|_{z=j}.
\end{align*}
It is easy to see that $\textrm{II}_{1,2}({\varphi})=\mathcal{L}_0(\varphi)$ and for $k=0,1$,
\ben
\partial_c\mathcal{L}_k(\varphi)=\frac{1}{u'(y_c)}\left(\mathcal{L}_k(\varphi')-I_{k,1}(\varphi)+I_{k,0}(\varphi)\right)+\mathcal{L}_{k+1}(\varphi).\label{eq:Lk-d}
\een
By Proposition \ref{prop:phi1}, we deduce that for $k=0,1,2,$ and $y\in [z,y_c]$ or $y\in [y_c,z]$,
\begin{align*}
\left|\left(\frac{\partial_z+\partial_y}{u'(y_c)}+\partial_c\right)^k
\left(\f{1}{(u(z)-c)^2}\left(\f{\phi_1(y,c)\,}{\phi_1(z,c)^2}-1\right)\right)\right|\leq C\f{\min\{\al^2|z-y_c|^2,1\}}{(u(z)-c)^2u'(y_c)^{2k}}.
\end{align*}
Using the fact that $|u(z)-c|\ge C^{-1}(z+y_c)|z-y_c|$, we get
\begin{align*}
&\int_0^1\f{|z-y_c|^{1-\f{1}{p}}\min\{\al^2|z-y_c|^2,1\}}{|u(z)-c|^2}\,dz\\
&\leq C\int_0^1\f{\min\{\al^2|z-y_c|^2,1\}}{(z+y_c)^2|z-y_c|^{1+\f1p}}\,dz\\
&\leq C\int_{|z-y_c|\leq \f{1}{\al},0\leq z\leq 1}\f{\al^2|z-y_c|^{1-\f1p}}{(z+y_c)^2}\,dz
+C\int_{|z-y_c|> \f{1}{\al},0\leq z\leq 1}\f{|z-y_c|^{-1-\f1p}}{(z+y_c)^2}\,dz\\
&\leq C\min\Big\{\frac{\al^{1+\f{1}{p}}}{u'(y_c)},
\frac{\al^{\f{1}{p}}}{u'(y_c)^2}\Big\},
\end{align*}
which implies that for $k=0,1,2,$ and $p\in [1,+\infty)$,
\begin{align}\label{eq:Lk-Lp}
\big|u'(y_c)^{2k}\mathcal{L}_k(\varphi)(c)\big|
&\leq C\int_0^1\f{|z-y_c|^{1-\f{1}{p}}\min\{\al^2|z-y_c|^2,1\}}{|u(z)-c|^2}dz\left\|\varphi\right\|_{L^{p}(I_0)}\\
&\leq C\min\Big\{\frac{\al^{1+\f{1}{p}}}{u'(y_c)},\frac{\al^{\f{1}{p}}}{u'(y_c)^2}\Big\}\left\|\varphi\right\|_{L^{p}(I_0)}.\nonumber
\end{align}
Similarly, we have for $k=0,1$,
\begin{align}
|I_{k,1}(\varphi)(c)|
&\leq C\f{|1-y_c|\min\{\al^2|1-y_c|^2,1\}}{(u(1)-c)^2u'(y_c)^{2k}}\left\|\varphi\right\|_{L^{\infty}(I_0)},\label{eq:Ik1}\\
|I_{k,0}(\varphi)(c)|
&\leq C\f{y_c\min\{\al^2y_c^2,1\}}{(c-u(0))^2u'(y_c)^{2k}}\left\|\varphi\right\|_{L^{\infty}(I_0)}.\label{eq:Ik0}
\end{align}

\begin{lemma}\label{lem:Pi12-Linfty}
It holds that for $p\in [1,\infty)$,
\beno
&&\|u'(y_c)^2\mathrm{II}_{1,2}({\varphi})\|_{L^{\infty}(I_0)}\leq C\al^{\f{1}{p}}\left\|\varphi\right\|_{L^{p}(I_0)},
\eeno
and for $p\in [1,\infty]$,
\beno
&&\|u'(y_c)\mathrm{II}_{1,2}(\varphi)\|_{L^{\infty}(I_0)}\leq C\al^{1+\f{1}{p}}\left\|\varphi\right\|_{L^{p}(I_0)}.
\eeno
\end{lemma}
\begin{proof}
The lemma follows from \eqref{eq:Lk-Lp} except the case of $p=\infty$. We have
\begin{align*}
&\int_0^1\f{|z-y_c|\min\{\al^2|z-y_c|^2,1\}}{|u(z)-c|^2}\,dz\\
&\leq C\int_{|z-y_c|\leq \f{1}{\al},0\leq z\leq 1}\f{\al^2|z-y_c|}{(z+y_c)^2}\,dz
+C\int_{|z-y_c|> \f{1}{\al},0\leq z\leq 1}\f{|z-y_c|^{-1}}{(z+y_c)^2}\,dz\\
&\leq C\int_0^1\f{|\al|}{(z+y_c)^2}\,dz
\leq \frac{C\al}{u'(y_c)},
\end{align*}
which gives
\beq\label{eq:Lk-Linftynew}
\|u'(y_c)^{2k+1}\mathcal{L}_{k}(\varphi)\|_{L^{\infty}(I_0)}\leq C|\al|\|\varphi\|_{L^{\infty}(I_0)}.
\eeq
Thus, the second inequality holds for $p=+\infty$.
\end{proof}

\begin{lemma}\label{lem:Pi12-Linfty-d1}
It holds that
\beno
&&\big\|u'(y_c)^3\pa_c\mathrm{II}_{1,2}(\varphi)\big\|_{L^{\infty}(I_0)}\leq C\al\|\varphi\|_{L^{\infty}(I_0)}+C\al^{\f12}\|\varphi'\|_{L^2(I_0)},\\
&&\|u'(y_c)^3\pa_c\mathrm{II}_{1,2}(\varphi)\|_{L^{\infty}(I_0)}\leq C\al\|\varphi\|_{L^{\infty}(I_0)}+C\al u'(y_c)\|\varphi'\|_{L^{\infty}(I_0)}.
\eeno
If $\varphi(0)=0$, then for $p\in (1,\infty)$,
\begin{align*}
&\|u'(y_c)\mathrm{II}_{1,2}(\varphi)\|_{L^{\infty}(I_0)}
\leq C\al^{\f1p}\|\varphi\|_{W^{1,p}(I_0)},\\
&\|\mathrm{II}_{1,2}(\varphi)\|_{L^{\infty}(I_0)}
\leq C\al\|\varphi\|_{W^{1,\infty}(I_0)},\\
&\|u'(y_c)\rho\pa_c\mathrm{II}_{1,2}(\varphi)\|_{L^{\infty}(I_0)}\leq C\al^{\f12}\|\varphi\|_{H^1(I_0)},\\
&\|u'(y_c)^2\pa_c\mathrm{II}_{1,2}(\varphi)\|_{L^{\infty}(I_0)}\leq C\al\|\varphi'\|_{L^\infty(I_0)}.
\end{align*}
\end{lemma}
\begin{proof}
By \eqref{eq:Lk-d}, we have
\begin{align*}
\pa_c\mathrm{II}_{1,2}(\varphi)
=\frac{1}{u'(y_c)}\left(\mathcal{L}_0(\varphi')-I_{0,1}(\varphi)+I_{0,0}(\varphi)\right)+\mathcal{L}_{1}(\varphi).
\end{align*}
It follows from \eqref{eq:Lk-Lp} and \eqref{eq:Lk-Linftynew} that
\beno
&&\|u'(y_c)^2\mathcal{L}_0(\varphi')\|_{L^{\infty}(I_0)}
\leq C\al^{\f12}\|\varphi'\|_{L^2(I_0)},\\
&&\|u'(y_c)^3\mathcal{L}_1(\varphi)\|_{L^{\infty}(I_0)}
\leq C\al\|\varphi\|_{L^{\infty}(I_0)},
\eeno
which along with \eqref{eq:Ik0} and \eqref{eq:Ik1} gives
\beno
\|u'(y_c)^3\pa_c\mathrm{II}_{1,2}(\varphi)\|_{L^{\infty}(I_0)}\leq C\al\|\varphi\|_{L^{\infty}(I_0)}+C\al^{\f12}\|\varphi'\|_{L^2(I_0)}.
\eeno
On the other hand, we have by \eqref{eq:Lk-Linftynew} that
\beno
\|u'(y_c)\mathcal{L}_0(\varphi')\|_{L^{\infty}(I_0)}
\leq C\al\|\varphi'\|_{L^{\infty}(I_0)},
\eeno
which gives
\beno
\|u'(y_c)^3\pa_c\mathrm{II}_{1,2}(\varphi)\|_{L^{\infty}(I_0)}\leq C\al\|\varphi\|_{L^{\infty}(I_0)}+C\al u'(y_c)\|\varphi'\|_{L^{\infty}(I_0)}.
\eeno

If $\varphi(0)=0$, then we have
\begin{align}
\label{eq:estL_kL^p0}|\mathcal{L}_k(\varphi)(c)|
&\leq C\int_0^1\f{\min\{\al^2|z-y_c|^2,1\}}{u'(y_c)^{2k+1}|z-y_c|^{1+\f1p}}dz\left\|\f{\varphi}{y}\right\|_{L^{p}(I_0)}\\
&\leq C\f{\al^{\f1p}}{u'(y_c)^{2k+1}}\|\varphi'\|_{L^p(I_0)}.\nonumber
\end{align}
Similar to \eqref{eq:Lk-Linftynew}, we have
\begin{align}\label{eq:LK-linfty-2}
\|u'(y_c)^{2k}\mathcal{L}_k(\varphi)\|_{L^{\infty}(I_0)}&\leq C\int_0^1\f{\min\{\al^2|z-y_c|^2,1\}}{|u(z)-c|}dz\left\|\f{\varphi}{y}\right\|_{L^{\infty}(I_0)}\\
&\leq C\int_0^1\f{\min\{\al^2|z-y_c|^2,1\}}{(z-y_c)^2}dz\left\|\f{\varphi}{y}\right\|_{L^{\infty}(I_0)}\nonumber\\
&\leq C\al\left\|\f{\varphi}{y}\right\|_{L^{\infty}(I_0)}\le C\al\left\|{\varphi'}\right\|_{L^{\infty}(I_0)}.\nonumber
\end{align}
Similarly, we can deduce that for $k=0,1$ and $j=0,1$,
\begin{align}
&\|u'(y_c)^{2k}\rho I_{k,j}(\varphi)\|_{L^\infty(I_0)}\leq C\al^{\f1p}\|\varphi\|_{W^{1,p}(I_0)},\label{eq:Ikj-W1p}\\
&\|u'(y_c)^{2k+1}I_{k,j}(\varphi)\|_{L^\infty(I_0)}\leq C\al\left\|{\varphi'}\right\|_{L^{\infty}(I_0)}.\label{eq:Ikj-Linfty}
\end{align}
Then the last two inequalities follow from \eqref{eq:Lk-d},\eqref{eq:Lk-Lp} and the above estimates.
\end{proof}

\begin{lemma}\label{lem:Pi12-d2}
It holds that
\begin{align*}
\|u'(y_c)^3\rho\pa_c^2\mathrm{II}_{1,2}(\varphi)\|_{L^\infty(I_0)}
\leq& C\al^{\f12}u'(y_c)\|\varphi''\|_{L^2(I_0)}+C\al\|\varphi\|_{L^{\infty}(I_0)}\\
&+C\al u'(y_c)\|\varphi'\|_{L^{\infty}(I_0)}+C\al^{\f12}\|\varphi'\|_{L^2(I_0)}.
\end{align*}
If $\varphi(0)=0$, then we have
\beno
\|u'(y_c)^2\rho^{}\partial_c^2\mathrm{II}_{1,2}(\varphi)\|_{L^{\infty}(I_0)}\leq C\al\left\|{\varphi}'\right\|_{L^{\infty}(I_0)}+C\al^{\f{1}{2}}\left\|{\varphi}''\right\|_{L^{2}(I_0)}.
\eeno
\end{lemma}
\begin{proof}
First of all, we have
\begin{align*}
\partial_cI_{0,j}(\varphi)(c)
&=-\f{\varphi(y_c)}{u'(y_c)}\f{1}{(u(z)-c)^2}\left(\f{1\,}{\phi_1(z,c)^2}-1\right)\bigg|_{z=j}\\
&\quad+\int_{y_c}^{j}\partial_c
\left(\f{{\varphi}(y)}{(u(z)-c)^2}\left(\f{\phi_1(y,c)\,}{\phi_1(z,c)^2}-1\right)\right)\,dy\bigg|_{z=j}.
\end{align*}
By Proposition \ref{prop:phi1}, we have for $y\in [z,y_c]$ or $y\in [y_c, z]$,
\begin{align*}
&\left|\partial_c\left(\f{1}{(u(z)-c)^2}\left(\f{\phi_1(y,c)\,}{\phi_1(z,c)^2}-1\right)\right)\right|\\
&\leq C\left|\left(\f{1}{|u(z)-c|^3}\left(\f{\phi_1(y,c)\,}{\phi_1(z,c)^2}-1\right)\right)\right|
+C\left|\f{1}{|u(z)-c|^2}\pa_c\left(\f{\phi_1(y,c)\,}{\phi_1(z,c)^2}\right)\right|\\
&\leq C\f{\min\{\al^2|z-y_c|^2,1\}}{|u(z)-c|^3}
+C\f{\min\{\al^2|z-y_c|,\al\}}{u'(y_c)|u(z)-c|^2}.
\end{align*}
Thus, we have
\begin{align*}
|\pa_cI_{0,1}(\varphi)(c)|
&\leq C\|\varphi\|_{L^{\infty}(I_0)}\f{\min\{\al^2(1-y_c)^2,1\}}{(1-y_c)^2u'(y_c)^3}\\
&\quad+C\left|\partial_c\left(\f{1}{(u(1)-c)^2}\left(\f{\phi_1(y,c)\,}{\phi_1(1,c)^2}-1\right)\right)\right||1-y_c|\|\varphi\|_{L^{\infty}}\\
&\leq C\|\varphi\|_{L^{\infty}(I_0)}\f{\min\{\al^2(1-y_c)^2,1\}}{(1-y_c)^2u'(y_c)^3}
+C\|\varphi\|_{L^{\infty}(I_0)}\f{\min\{\al^2(1-y_c),\al\}}{|1-y_c|u'(y_c)^3},
\end{align*}
and
\begin{align*}
|\pa_cI_{0,0}(\varphi)(c)|\leq C\|\varphi\|_{L^{\infty}(I_0)}\Big(\f{\min\{\al^2 y_c^2,1\}}{y_c^2u'(y_c)^3}+\f{\min\{\al^2 y_c,\al\}}{y_cu'(y_c)^3}\Big).
\end{align*}
This implies that
\begin{align}
\|\pa_cI_{0,1}(\varphi)\|_{L^\infty(I_0)}+\|\pa_cI_{0,0}(\varphi)\|_{L^\infty(I_0)}
\leq C\|\varphi\|_{L^{\infty}(I_0)}\f{\min\{\al^2 \rho_0^2,\al\rho_0\}}{u'(y_c)^3\rho_0^2}.\label{eq:Ikj-d}
\end{align}

A direct calculation with \eqref{eq:Lk-d} gives
\begin{align*}
\pa_c^2\mathrm{II}_{1,2}(\varphi)
&=\pa_c\left(\frac{1}{u'(y_c)}\left(\mathcal{L}_0(\varphi')-I_{0,1}(\varphi)+I_{0,0}(\varphi)\right)+\mathcal{L}_{1}(\varphi)\right)\\
&=\frac{1}{u'(y_c)^2}\left(\mathcal{L}_0(\varphi'')-I_{0,1}(\varphi')+I_{0,0}(\varphi')\right)+\f{1}{u'(y_c)}\mathcal{L}_{1}(\varphi')\\
&\quad-\f{1}{u'(y_c)}\pa_c(I_{0,1}(\varphi)(c)-I_{0,0}(\varphi)(c))\\
&\quad-\f{u''(y_c)}{u'(y_c)^3}\left(\mathcal{L}_0(\varphi')-I_{0,1}(\varphi)+I_{0,0}(\varphi)\right)\\
&\quad+\frac{1}{u'(y_c)}\left(\mathcal{L}_1(\varphi')-I_{1,1}(\varphi)+I_{1,0}(\varphi)\right)+\mathcal{L}_{2}(\varphi),
\end{align*}
which along with \eqref{eq:Lk-Lp}, \eqref{eq:Lk-Linftynew}, \eqref{eq:Ik0},\eqref{eq:Ik1} and \eqref{eq:Ikj-d} gives the first inequality of the lemma.

If  $\varphi(0)=0$, then we have
\beno
|\pa_cI_{0,0}(\varphi)(c)|\leq C\left\|\f{\varphi(y)}{y}\right\|_{L^{\infty}(I_0)}\Big(\f{\min\{\al^2 y_c^2,1\}}{y_c^2u'(y_c)^2}+\f{\min\{\al^2 y_c,\al\}}{y_cu'(y_c)^2}\Big),
\eeno
which gives
\beno
\|u'(y_c)^{2}y_c\pa_cI_{0,0}(\varphi)\|_{L^\infty(I_0)}\leq C\al \left\|{\varphi'}\right\|_{L^{\infty}(I_0)}.
\eeno
Similarly, we have
\beno
\|u'(y_c)^{2}(1-y_c)\pa_cI_{0,1}(\varphi)\|_{L^\infty(I_0)}\leq C\al\left\|{\varphi'}\right\|_{L^{\infty}(I_0)}.
\eeno
Then the second inequality follows by combining the above inequalities with \eqref{eq:Lk-Lp}, \eqref{eq:LK-linfty-2} and \eqref{eq:Ikj-Linfty}.
\end{proof}

\section{$W^{2,1}$ estimate of the kernel}

For the sake of simplicity, we introduce the following notations:

\begin{itemize}
\item We denote by ${\mathcal L}^p$ a function $f$ which satisfies $\|f\|_{L_{y_c}^p}\leq C$.

\item We denote by $a{\mathcal L}^p\cap b{\mathcal L}^q$  a function $f$ which satisfies $\|f/a\|_{L_{y_c}^p}+\|f/b\|_{L_{y_c}^q}\leq C$.
\end{itemize}
Here the constant $C$ is independent of $\al$, and may be different from line to line.

\subsection{$W^{2,1}$ estimate of $K_o(c,\al)$}
Recall that
\begin{align*}
K_o(c,\al)=\f{\Lambda_1(\widehat{\om_{o}})(c)\Lambda_2(g)(c)}{(\mathrm{A}(c)^2+\mathrm{B}(c)^2){u'(y_c)}}.
\end{align*}

\begin{proof}[Proof of Proposition \ref{prop:Ko}]

{\bf Step 1.} $L^1$ estimate. We normalize $\|\widehat{\om}_{o}\|_{L^2}\leq 1,\ \|g\|_{L^2}\leq 1$. By Lemma \ref{lem:Pi12-Linfty} and Lemma \ref{lem:II3}, we get
\beno
&&\Lambda_{1,2}(\widehat{\om}_o)=\al\rho_0\cL^2+\al^{\f12}\cL^{\infty},\\
&&\Lambda_{2,2}(g)=\al\rho_0\cL^2+\al^{\f12}\cL^{\infty},
\eeno
here $\rho_0=\f{\rho}{u'(y_c)}\sim y_c(1-y_c)$. Then we get by Lemma \ref{lem:Lambda-11} that
\beno
&&\Lambda_{1}(\widehat{\om}_o)=(1+\al\rho_0)\cL^2+\al^{\f12}\cL^{\infty},\\
&&\Lambda_{2}(g)=(1+\al\rho_0)\cL^2+\al^{\f12}\cL^{\infty},
\eeno
which gives
\beno
\Lambda_{1}(\widehat{\om}_o)\Lambda_{2}(g)=(1+\al\rho_0)^2\cL^1+(1+\al\rho_0)\al^\f12\cL^2+\al\cL^\infty.
\eeno
By Lemma \ref{lem:A2+B2}, we have $\f{1}{\rmA^2+\rmB^2}=\f{CL^{\infty}}{(1+\al\rho_0)^2}$.
Thus, we obtain
\beno
\f{\Lambda_1(\widehat{\om_{o}})(c)\Lambda_2(g)(c)}{(\mathrm{A}(c)^2+\mathrm{B}(c)^2)}=\cL^1+\f {\al^\f12\cL^2} {1+\al\rho_0}=\cL^1,
\eeno
which gives
\beno
\big\|K_o(c,\al)\big\|_{L_c^1}
=\big\|u'(y_c)K_o(c,\al)\big\|_{L_{y_c}^1}\leq C.
\eeno

{\bf Step 2}. $W^{1,1}$ estimate. We normalize $\|\widehat{\om}_o\|_{H^1}\leq 1$ and $\|g'\|_{L^2}+\al\|g\|_{L^2}\leq 1$. Then we have
\ben
\|g\|_{L^{\infty}}\leq \|g\|_{L^2}^{\f12}\|g\|_{H^1}^{\f12}\leq C\al^{-\f12}.\label{eq:g-linfty}
\een
Thanks to $\widehat{\om}_o(\al, 0)=g(0)=g(1)=0$, we get by Hardy's inequality that
\ben
\big\|{\widehat{\om}_o}/{y}\big\|_{L^2}+\big\|{g}/{\rho_0}\big\|_{L^2}\leq C.\label{eq:gw-hardy}
\een
By Lemma \ref{lem:Pi12-Linfty-d1}, Lemma \ref{lem:II3}, \eqref{eq:g-linfty} and \eqref{eq:gw-hardy}, we have
\begin{align*}
\Lambda_{1,2}(\widehat{\om}_o)&=\al\rho_0\cL^{\infty}\cap\al\rho\cL^2+\al^{\f12}\rho_0\cL^{\infty},\\
\Lambda_{2,2}(g)&=\al^{\f12}\rho_0\cL^{\infty},\\
\pa_c\Lambda_{1,2}(\widehat{\om}_o)&=\al^{\f12}(1-y_c)\cL^{\infty}+\al^{\f12}\cL^{\infty}/u'+\al \cL^2\\
&=\f{(\al^{\f12}\cL^{\infty}+\al u'(y_c)\cL^2)}{u'(y_c)},\\
\pa_c\Lambda_{2,2}(g)&=\f{(\al^{\f12}\cL^{\infty}+\al {\rho_0}\cL^2)}{u'(y_c)}.
\end{align*}
from which and Lemma \ref{lem:Lambda-11}, we infer that
\begin{align*}
\Lambda_{1}(\widehat{\om}_o)&=\al\rho_0\cL^{\infty}\cap\al\rho\cL^2+\al^{\f12}\rho_0\cL^{\infty}+{y_c}\cL^2\\
&=\big(u'(y_c)(1+\al\rho_0)\cL^2+\al^{\f12}\rho_0\cL^{\infty}\big),\\
\Lambda_{2}(g)&=\al^{\f12}\rho_0\cL^{\infty}+\rho_0\cL^2,\\
\pa_c\Lambda_{1}(\widehat{\om}_o)&=\f{(\al^{\f12}\cL^{\infty}+(1+\al u'(y_c))\cL^2)}{u'(y_c)},\\
\pa_c\Lambda_{2}(g)&=\f{(\al^{\f12}\cL^{\infty}+(1+\al {\rho_0})\cL^2)}{u'(y_c)}.
\end{align*}
Then we obtain
\begin{align*}
\Lambda_{1}(\widehat{\om}_o)\Lambda_2(g)
&=\big(u'(y_c)(1+\al\rho_0)\cL^2+\al^{\f12}\rho_0\cL^{\infty}\big)\big(\al^{\f12}\rho_0\cL^{\infty}+\rho_0\cL^2\big)\\
&=\big(\al\rho_0^2 \cL^{\infty}+(\al^{\f12}\rho_0^2+\al^{\f12}\rho(1+\al\rho_0))\cL^{2}+\rho(1+\al\rho_0)\cL^1\big)\\
&=\big(\al\rho_0^2 \cL^{\infty}+\al^{\f12}\rho(1+\al\rho_0)\cL^{2}+\rho(1+\al\rho_0)\cL^1\big),
\end{align*}
and
\begin{align*}
\pa_c\big(\Lambda_{1}(\widehat{\om}_o)\Lambda_2(g)\big)
&=\Lambda_2(g)\pa_c\Lambda_{1}(\widehat{\om}_o)+\Lambda_{1}(\widehat{\om}_o)\pa_c\Lambda_2(g)\\
&=\rho_0\big(\al^{\f12}\cL^{\infty}+\cL^2\big)\f{C(\al^{\f12}\cL^{\infty}+(1+\al u'(y_c))\cL^2)}{u'(y_c)}\\
&\quad+\big(u'(y_c)(1+\al\rho_0)\cL^2+\al^{\f12}\rho_0\cL^{\infty}\big)\f{(\al^{\f12}\cL^{\infty}+(1+\al \rho_0)\cL^2)}{u'(y_c)}\\
&=\big(\al\cL^{\infty}+\al^{\f12}(1+\al\rho_0)\cL^2+(1+\al\rho_0)\cL^1\big)\\
&\quad+\big(\al\cL^{\infty}+\al^{\f12}(1+\al\rho_0)\cL^2+(1+\al\rho_0)^2\cL^1\big)\\
&=\big(\al\cL^{\infty}+\al^{\f12}(1+\al\rho_0)\cL^2+(1+\al\rho_0)^2\cL^1\big).
\end{align*}
By Lemma \ref{lem:A2+B2}, we have
\beno
\pa_c\Big(\f{1}{(\rmA^2+\rmB^2)u'(y_c)}\Big)=\f{\cL^{\infty}(1+\al u'(y_c))}{(1+\al\rho_0)^3{u'(y_c)^3}}
+\f{(\cL^2\cap (\cL^{\infty}/\rho))}{(1+\al\rho_0)^3u'(y_c)}.
\eeno

Summing up, we deduce that
\begin{align*}
&u'(y_c)\pa_cK_o(\al,c)\\
&=u'(y_c)\Lambda_{1}(\widehat{\om}_o)\Lambda_2(g)\pa_c\Big(\f{1}{(\rmA^2+\rmB^2)u'(y_c)}\Big)
+\f{\pa_c\big(\Lambda_{1}(\widehat{\om}_o)\Lambda_2(g)\big)}{(\rmA^2+\rmB^2)}\\
&=\big(\al\rho_0^2 \cL^{\infty}+\al^{\f12}\rho(1+\al\rho_0)\cL^{2}+\rho(1+\al\rho_0)\cL^1\big)\left(\f{\cL^{\infty}(1+\al u'(y_c))}{(1+\al\rho_0)^3{u'(y_c)^2}}
+\f{(\cL^2\cap (\cL^{\infty}/\rho))}{(1+\al\rho_0)^3}\right)\\
&\quad+\f{\big(\al\cL^{\infty}+\al^{\f12}(1+\al\rho_0)\cL^2+(1+\al\rho_0)^2\cL^1\big)}{(1+\al\rho_0)^2}\\
&=\f{\al\cL^{\infty}}{(1+\al\rho_0)^2}+\f{\al^{\f12}\cL^{2}}{1+\al\rho_0}+\cL^1+\f{\cL^{2}}{1+\al\rho_0}+\f{\al^{\f12}\cL^{2}}{(1+\al\rho_0)^2}\\
&=\f{\al\cL^{\infty}}{(1+\al\rho_0)^2}+\f{\al^{\f12}\cL^{2}}{1+\al\rho_0}+\cL^1=\cL^1,
\end{align*}
which gives
\beno
\big\|\pa_cK(c,\al)\big\|_{L_c^1}=\big\|u'(y_c)\pa_cK_o(c,\al)\big\|_{L_{y_c}^1}\leq C.
\eeno

{\bf Step 3.} $W^{2,1}$ estimate. We normalize $\|\widehat{\om}_o\|_{H^2}\leq 1$ and $\|g''\|_{L^2}+\al^2\|g\|_{L^2}\leq 1$. Then we get
\ben\label{eq:g-linfty2}
\|g\|_{L^2}\leq \f{C}{\al^2},\quad \|g'\|_{L^2}\leq \f{C}{\al},\quad \|g'\|_{L^{\infty}}\leq C\|g'\|_{L^2}^{\f12}\|g''\|_{L^2}^{\f12}\leq \f{C}{\sqrt{\al}}.
\een
Thanks to $\widehat{\om}_o(\al,0)=g(0)=g(1)=0$, we have
\ben\label{eq:gw-h2}
\|g/\rho_0\|_{L^\infty}\leq C\|{g'}\|_{L^\infty}\leq \f{C}{\sqrt{\al}},\quad
\|\widehat{\om}_o/y\|_{L^\infty}\le C.
\een
By Lemma \ref{lem:Pi12-Linfty-d1}, Lemma \ref{lem:II3}, Lemma \ref{lem:Pi12-d2}, \eqref{eq:g-linfty2} and \eqref{eq:gw-h2}, we have
\begin{align*}
&\Lambda_{1,2}(\widehat{\om}_o)=\al\rho \cL^{\infty},\quad \Lambda_{2,2}(g)=\al^{\f12}\rho \cL^{\infty},\\
&\pa_c\Lambda_{1,2}(\widehat{\om}_o)=\al \cL^{\infty}, \quad
\pa_c\Lambda_{2,2}(g)=\al^{\f12} \cL^{\infty},\\
&\pa_c^2\Lambda_{1,2}(\widehat{\om}_o)=\rho^{-1}\al \cL^{\infty}+(u'(y_c))^{-1}\al \cL^2,\\
&\pa_c^2\Lambda_{2,2}(g)=u'(y_c)^{-2}(\al^{\f12} \cL^{\infty}+\al\rho_0 \cL^2),
\end{align*}
from which and Lemma \ref{lem:Lambda-11}, we infer that
\begin{align*}
&\Lambda_1({\widehat{\om}_o})=(\rho \cL^2\cap u'(y_c)\cL^{\infty})+(u'(y_c)^2+\al\rho)\cL^{\infty},\\
&\partial_c\big(\Lambda_1({\widehat{\om}_o})\big)=(\cL^2+\al \cL^{\infty}),\\
&\partial_c^2\big(\Lambda_1({\widehat{\om}_o})\big)=\rho^{-1}( (1+\al\rho_0)\cL^2+\al \cL^{\infty}),
\end{align*}
and
\begin{align*}
&\Lambda_2(g)
=\rho \cL^2+\al^{\f{1}{2}}\rho \cL^{\infty},\\
&\partial_c\big(\Lambda_2(g)\big)
=(\cL^2+\al^{\f{1}{2}}\cL^{\infty}),\\
&\partial_c^2\big(\Lambda_2(g)\big)
=u'(y_c)^{-2}( (1+\al\rho_0)\cL^2+\al^{\f{1}{2}}\cL^{\infty}).
\end{align*}
By Lemma \ref{lem:A2+B2} , we have
\beno
\pa_c^2\Big(\f{1}{(\rmA^2+\rmB^2)u'(y_c)}\Big)
=\f{\cL^{\infty}(1+\al^2u'(y_c)^2)}{(1+\al\rho_0)^4{u'(y_c)^5}}
+\f{(\cL^2+\al \cL^{\infty})}{(1+\al\rho_0)^3{\rho}u'(y_c)}.
\eeno
Thus, we can deduce that
\begin{align*}
\Lambda_1(\widehat{\om}_o)\Lambda_2(g)
&=\rho\big(\rho \cL^2\cap u'(y_c)\cL^{\infty}+(u'(y_c)^2+\al\rho)\cL^{\infty}\big)\big( \cL^2 +\al^{\f{1}{2}} \cL^{\infty}\big)\\
&=\rho\big(\rho\cL^1\cap u'\cL^2+(u'(y_c)^2+\al\rho)\cL^{2}+\al^{\f12}(\rho \cL^2\cap u'(y_c)\cL^{\infty})+\al^\f32\rho\cL^\infty\big),
\end{align*}
and
\begin{align*}
&\pa_c\big(\Lambda_1(\widehat{\om}_o)\Lambda_2(g)\big)\\
&=\Lambda_2(g)\pa_c\Lambda_1(\widehat{\om}_o)+\Lambda_1(\widehat{\om}_o)\pa_c\Lambda_2(g)\\
&=\rho \big(\cL^2 +\al^{\f{1}{2}}\cL^{\infty}\big)( \cL^2+\al \cL^{\infty})\\
&\quad+\big(\rho \cL^2\cap u'(y_c)\cL^{\infty}+(u'(y_c)^2+\al\rho)\cL^{\infty}\big)(\cL^2+\al^{\f{1}{2}} \cL^{\infty})\\
&=\rho(\cL^1+\al\cL^2+\al^{\f32}\cL^{\infty})+\rho\cL^1\cap u'(y_c)\cL^2\\
&\quad+(u'(y_c)^2+\al\rho)\cL^2+\al^{\f12}(\rho\cL^2\cap u'(y_c)\cL^{\infty})+\al^{\f12}(u'(y_c)^2+\al\rho)\cL^{\infty}\\
&=\rho\cL^1+(u'(y_c)^2+\al\rho)\cL^2+\al^{\f12}(u'(y_c)^2+\al\rho)\cL^{\infty},
\end{align*}
and
\begin{align*}
&\pa_c^2\big(\Lambda_1(\widehat{\om}_o)\Lambda_2(g)\big)=\Lambda_2(g)\pa_c^2\Lambda_1(\widehat{\om}_o)+\Lambda_1(\widehat{\om}_o)\pa_c^2\Lambda_2(g)+2\pa_c\Lambda_2(g)\pa_c\Lambda_1(\widehat{\om}_o)\\
&=\rho (\cL^2 +\al^{\f{1}{2}}\cL^{\infty})\rho^{-1}( (1+\al\rho_0)\cL^2+\al \cL^{\infty})\\
&\quad+\big((\rho \cL^2\cap u'(y_c)\cL^{\infty})+(u'(y_c)^2+\al\rho)\cL^{\infty}\big)u'(y_c)^{-2}( (1+\al\rho_0)\cL^2+\al^{\f{1}{2}} \cL^{\infty})\\
&\quad+( \cL^2+\al \cL^{\infty})( \cL^2+\al^{\f{1}{2}} \cL^{\infty})\\
&=\big((1+\al\rho_0)\cL^1+\al^{\f12}(1+\al\rho_0)\cL^2+\al\cL^2+\al^{\f32}\cL^{\infty}\big)\\
&\quad+\big((1+\al\rho_0)\cL^1+(1+\al(1-y_c))(1+\al\rho_0)\cL^2+\al^{\f12}\cL^2+\al^{\f32}\cL^{\infty}\big)\\
&\quad+(\cL^1+\al\cL^2+\al^{\f32}\cL^{\infty})\\
&=\big((1+\al\rho_0)\cL^1+\al(1+\al\rho_0)\cL^2+\al^{\f32}\cL^{\infty}\big).
\end{align*}

With the above estimates, it is easy to deduce that
\begin{align*}
u'(y_c)\pa_c^2K_o(c,\al)=&\f{\pa_c^2\big(\Lambda_1(\widehat{\om}_o)\Lambda_2(g)\big)}{\rmA^2+\rmB^2}+2u'(y_c)\pa_c\big(\Lambda_1(\widehat{\om}_o)\Lambda_2(g)\big)\pa_c\Big(\f{1}{(\rmA^2+\rmB^2)u'(y_c)}\Big)\\
&\quad+u'(y_c)\Lambda_1(\widehat{\om}_o)\Lambda_2(g)\pa_c^2\Big(\f{1}{(\rmA^2+\rmB^2)u'(y_c)}\Big)=\al^{\f12}\cL^1,
\end{align*}
which gives
\beno
\|\pa_c^2K_o(c,\al)\|_{L^1_c}\leq C\|u'(y_c)\pa_c^2K_o(c,\al)\|_{L^1_{y_c}}\leq C\al^{\f12}.
\eeno

The above estimates also show that
\beno
K_o(c,\al)=\f{\rho(\cL^2+\al^{\f12}\cL^{\infty})\Big((\rho \cL^2\cap u'(y_c)\cL^{\infty})+(u'(y_c)^2+\al\rho)\cL^{\infty}\Big)}{(1+\al\rho_0)^2u'(y_c)}=C_\al\rho \cL^2,
\eeno
which along with the fact $K_o$ is continuous implies that
\beno
 K_o(u(0),\al)=K_o(u(1),\al)=0.
\eeno

This completes the proof of the proposition.
\end{proof}

\subsection{$W^{2,1}$ estimate of $K_e(c,\al)$}

Let us first present the following estimates of $\Lambda_{3,1}(\widehat{\om}_e)$  and $\Lambda_{4,1}(g)$, which are defined as follows
\beno
&&\Lambda_{3,1}(\widehat{\om}_e)(c)=J(c)\Big(\frac{u''(y_c)}{u'(y_c)}\rmE(\widehat{\om}_e)(c)+{\widehat{\omega}_{e}(y_c)}\Big),\\
&&\Lambda_{4,1}(g)(c)=J(c)\Big(\frac{\rmE(gu'')(c)}{u'(y_c)}+g(y_c)\Big),
\eeno
where $\mathrm{E}(\varphi)(c)=\int_{y_c}^0\varphi\phi_1(y,c)\,dy$.

\begin{lemma}\label{lem:Lambda_34,1}
1. If $\|\widehat{\om}_e\|_{L^2}\leq 1$ and $\|g\|_{L^2}\leq 1$, then
\begin{align*}
\Lambda_{3,1}(\widehat{\om}_e)=\al^{-2}\cL^2,\quad
\Lambda_{4,1}(g)=\al^{-2}\cL^2.
\end{align*}
2. If  $\|\widehat{\om}_e\|_{H^1}\leq 1$ and $|\al|\|g\|_{L^2}+\|g'\|_{L^2}\leq 1$, then
\begin{align*}
&\Lambda_{3,1}(\widehat{\om}_e)=\frac{\rho_0(\cL^2+\al \cL^{\infty})}{\al^2},\quad
\Lambda_{4,1}(g)=\frac{\rho_0(\cL^2+\al^{\frac{1}{2}}\cL^{\infty})}{\al^2},\\
&\pa_c\big(\Lambda_{3,1}(\widehat{\om}_e)\big)
=\frac{\cL^2+\al \cL^{\infty} }{\al^2 y_c},\quad
\pa_c\big(\Lambda_{4,1}(g)\big)
=\frac{\cL^2+\al^{\frac{1}{2}}\cL^{\infty}}{\al^2 y_c}.
\end{align*}
3. If  $\|\widehat{\om}_e\|_{H^2}\leq 1$ and $\|g''-\al^2g\|_{L^2}\leq 1$ and $\widehat{\om}_e'(\al,0)=g'(0)=0$, then
\beno
&&\Lambda_{3,1}(\widehat{\om}_e)=\rho \cL^{\infty}+C\al^{-2}(\rho_0\cL^{\infty}\cap \rho \cL^2),\\
&&\pa_c\Lambda_{3,1}(\widehat{\om}_e)=\cL^{\infty}+\al^{-2}\cL^2,\\
&&\pa_c^2\Lambda_{3,1}(\widehat{\om}_e)=\rho^{-1}(\cL^{\infty}+\al^{-2}\cL^2),
\eeno
and
\beno
&&\Lambda_{4,1}(g)=\rho\al^{-\f32}\cL^{\infty}+C\al^{-2}\rho \cL^2,\\
&&\pa_c\Lambda_{4,1}(g)=\al^{-\f32}\cL^{\infty}+C\al^{-2}\cL^2,\\
&&\pa_c^2\Lambda_{4,1}(g)=y_c^{-2}\big(\al^{-\f32}\cL^{\infty}+\al^{-2}\cL^2\big).
\eeno
\end{lemma}
\begin{proof}
{\bf Case 1.} $\|\widehat{\om}_e\|_{L^2}\leq 1$ and $\|g\|_{L^2}\leq 1$.

Due to $\phi_1(y,c)\leq \phi_1(0,c)$, we have
\begin{align*}
\Big\|\f{\rmE(\varphi)}{\phi_1(0,c)y_c}\Big\|_{L^2}\leq C\|\varphi\|_{L^2},
\end{align*}
which along with Lemma \ref{lem:J} implies that
\begin{align*}
\Lambda_{3,1}(\widehat{\om}_e)=\al^{-2}\cL^2,\quad
\Lambda_{4,1}(g)=\al^{-2}\cL^2.
\end{align*}

{\bf Case 2.} $\|\widehat{\om}_e\|_{H^1}\leq 1$ and $|\al|\|g\|_{L^2}+\|g'\|_{L^2}\leq 1$.
Thus, by Sobolev embedding, $\|\widehat{\om}_e\|_{L^\infty}\le 1$ and $\|g\|_{L^{\infty}}\leq C\al^{-\frac{1}{2}}$.

We write
\begin{align}
\nonumber&u''(y_c)\rmE(\widehat{\om}_e)+u'(y_c)\widehat{\om}_e(y_c)\\
\nonumber&=\int_{y_c}^{0}\big(u''(y_c)\phi_1(y,c)\widehat{\om}_e(y)-u''(y)\widehat{\om}_e(y_c)\phi_1(y_c,c)\big)dy\\
\nonumber&=\int_{y_c}^{0}\big(u''(y_c)-u''(y)\big)\phi_1(y,c)\widehat{\om}_e(y)dy\\
\nonumber&\quad+\int_{y_c}^{0}\big(\phi_1(y,c)-1\big)u''(y)\widehat{\om}_e(y)dy+\int_{y_c}^{0}u''(y)\big(\widehat{\om}_e(y)-\widehat{\om}_e(y_c)\big)dy\\
&=I_1+I_2+I_3.\nonumber
\end{align}
Using the facts that $\phi_1(y,c)-1\leq C\min\{\al^2|y-y_c|^2,1\}\phi_1(y,c)$ and $|u''(y)-u''(y_c)|\le C|y-y_c|y_c$ for $y\le y_c$(since $u$ is even), it is easy to get
\begin{align*}
&|I_1|\leq C|y_c|^3\phi_1(0,c)\|\widehat{\om}_e\|_{L^{\infty}},\\
&|I_2|\leq C|\al| |y_c|^2\min\big\{1,|\al||y_c|\big\}\phi_1(0,c)\|\widehat{\om}_e\|_{L^{\infty}},\\
&\big\|{I_3}/{y_c^2}\big\|_{L^2}\leq C\|\widehat{\om}_e'\|_{L^2}.\nonumber
\end{align*}
This shows that
\begin{align*}
&{u''(y_c)}\rmE(\widehat{\om}_e)+u'(y_c){\widehat{\om}_e(y_c)}
=\phi_1(0,c)y_c^2\al\cL^{\infty}+y_c^2\cL^2,
\end{align*}
which along with Lemma \ref{lem:J} gives
\begin{align}
\Lambda_{3,1}(\widehat{\om}_e)=\frac{\rho_0(\cL^2+\al \cL^{\infty})}{\al^2}.\label{eq:Lambda31-est1}
\end{align}
We write
\begin{align}
\nonumber&\rmE(u''g)+u'(y_c)g(y_c)\\
\nonumber&=\int_{y_c}^{0}\big(u''(y)\phi_1(y,c)g(y)-u''(y)g(y_c)\phi_1(y_c,c)\big)dy\\
\nonumber&=\int_{y_c}^{0}\big(\phi_1(y,c)-1\big)u''(y)g(y)dy+\int_{y_c}^{0}u''(y)\big(g(y)-g(y_c)\big)dy\\
&=I_1+I_2.\nonumber
\end{align}
It is easy to see that
\begin{align*}
&|I_1| \leq C|\al|^{\f12} |y_c|^2\min\big\{1,|\al| |y_c|\big\}\phi_1(0,c),\nonumber\\
&\big\|{I_2}/{y_c^2}\big\|_{L^2}\leq C\|g'\|_{L^2}.\nonumber
\end{align*}
This shows that
\begin{align*}
{\rmE(gu'')}+{u'(y_c)}g(y_c)=\phi_1(0,c)y_c^2\al^{\frac{1}{2}}\cL^{\infty}+y_c^2\cL^2,
\end{align*}
which along with Lemma \ref{lem:J} gives
\begin{align}
\Lambda_{4,1}(g)=\frac{\rho_0(\cL^2+\al^{\frac{1}{2}}\cL^{\infty})}{\al^2}.\label{eq:Lambda41-est1}
\end{align}

A direct calculation gives
\begin{align*}
&\partial_c\left({u''(y_c)}\rmE(\widehat{\om}_e)+u'(y_c){\widehat{\om}_e(y_c)}\right)\\
&=\f{u'''(y_c)}{u'(y_c)}\rmE(\widehat{\om}_e)
-\f{u''(y_c)\widehat{\om}_e(y_c)}{u'(y_c)}
+{u''(y_c)}\int_{y_c}^0\widehat{\om}_e\partial_c\phi_1(y,c)dy
+\f{({u'\widehat{\om}_e})'(y_c)}{u'(y_c)}\\
&=\f{u'''(y_c)}{u'(y_c)}\rmE(\widehat{\om}_e)
+\widehat{\om}_e'(y_c)+{u''(y_c)}\int_{y_c}^0\widehat{\om}_e\partial_c\phi_1(y,c)dy,
\end{align*}
and
\begin{align*}
\partial_c\left({\rmE_0(gu'')}+{u'(y_c)}g(y_c)\right)
&=-\f{gu''(y_c)}{u'(y_c)}+\int_{y_c}^0gu''\partial_c\phi_1(y,c)dy
+\f{(u'g)'(y_c)}{u'(y_c)}\\
&=g'(y_c)+\int_{y_c}^0gu''\partial_c\phi_1(y,c)dy.
\end{align*}
By Proposition \ref{prop:phi1}, we know that for $0\leq y\leq y_c$,
\beno
\Big|\f{\pa_c\phi_1(y,c)}{\phi_1(0,c)}\Big|\leq \f{C\al\min\{\al y_c,1\}}{y_c}.
\eeno
Thus, we can deduce that
\begin{align*}
&\partial_c\left({u''(y_c)}\rmE(\widehat{\om}_e)+u'(y_c){\widehat{\om}_e(y_c)}\right)
=\phi_1(0,c)\al\cL^{\infty}+\cL^2,\\
&\partial_c\left({\rmE(gu'')}+{u'(y_c)}g(y_c)\right)=\phi_1(0,c)\al^{\frac{1}{2}}\cL^{\infty}+\cL^2,
\end{align*}
from which and Lemma \ref{lem:J}, we infer that
\begin{align}
&\pa_c\big(\Lambda_{3,1}(\widehat{\om}_e)\big)
=\frac{\cL^2+\al \cL^{\infty} }{\al^2 y_c},\quad
\pa_c\big(\Lambda_{4,1}(g)\big)
=\frac{\cL^2+\al^{\frac{1}{2}}\cL^{\infty}}{\al^2 y_c}.\label{eq:Lambda3141-est2}
\end{align}

{\bf Case 3}. $\|\widehat{\om}_e\|_{H^2}\leq 1$ and $\|g''-\al^2g\|_{L^2}\leq 1$ and $\widehat{\om}_e'(\al,0)=g'(0)=0$. Thus, we have
$\|\widehat{\om}_e\|_{L^\infty}+\|\widehat{\om}_e'\|_{L^\infty}\le C$ and $\|g\|_{L^\infty}\le C\al^{-\f32}, \|g'\|_{L^\infty}\le C\al^{-\f12}$.

First of all, following the proof of \eqref{eq:Lambda31-est1}, \eqref{eq:Lambda41-est1} and \eqref{eq:Lambda3141-est2}, we can deduce that
\begin{align*}
&u''(y_c)\rmE(\widehat{\om}_e)+u'(y_c)\widehat{\om}_e(y_c)
=\phi_1(0,c)y_c^3\al^2\cL^{\infty}+(y_c^3\cL^2\cap y_c^2\cL^{\infty}),\\
&\partial_c\left({u''(y_c)}\rmE(\widehat{\om}_e)+u'(y_c){\widehat{\om}_e(y_c)}\right)
=\phi_1(0,c)\al^2y_c\cL^{\infty}+y_c\cL^2,\\
&\rmE(gu'')+u'(y_c)g(y_c)
=\phi_1(0,c)y_c^3\al^{\frac{1}{2}}\cL^{\infty}+y_c^3\cL^2,\\
&\partial_c\left(\rmE(gu'')+u'(y_c)g(y_c)\right)
=\phi_1(0,c)\al^{\frac{1}{2}}y_c\cL^{\infty}+y_c\cL^2,
\end{align*}
from which and Lemma \ref{lem:J}, we infer that
\beno
&&\Lambda_{3,1}(\widehat{\om}_e)=\rho \cL^{\infty}+C\al^{-2}(\rho_0\cL^{\infty}\cap \rho \cL^2),\quad\pa_c\Lambda_{3,1}(\widehat{\om}_e)=\cL^{\infty}+\al^{-2}\cL^2,\\
&&\Lambda_{4,1}(g)=\rho\al^{-\f32}\cL^{\infty}+C\al^{-2}\rho \cL^2,\quad\pa_c\Lambda_{4,1}(g)=\al^{-\f32}\cL^{\infty}+C\al^{-2}\cL^2.
\eeno

We have
\begin{align*}
\partial_c^2\left({\rmE(gu'')}+{u'(y_c)}g(y_c)\right)
=\f{g''(y_c)}{u'(y_c)}+\int_{y_c}^0gu''\partial_c^2\phi_1(y,c)dy,
\end{align*}
and
\begin{align*}
&\partial_c^2\left({u''(y_c)}\rmE(\widehat{\om}_e)+u'(y_c){\widehat{\om}_e(y_c)}\right)\\
&=\left(\f{u''''(y_c)}{u'(y_c)^2}-\f{(u'''u'')(y_c)}{u'(y_c)^3}\right)\rmE(\widehat{\om}_e)
+\f{\widehat{\om}_e''(y_c)}{u'(y_c)}+2\f{u'''(y_c)}{u'(y_c)}\int_{y_c}^0\widehat{\om}_e\partial_c\phi_1(y,c)dy\\
&\quad+{u''(y_c)}\int_{y_c}^0\widehat{\om}_e\partial_c^2\phi_1(y,c)dy.
\end{align*}
By Proposition \ref{prop:phi1}, we know that for $0<y<y_c$,
\beno
\left|\dfrac{\partial_c^2\phi_1(y,c)}{\phi_1(y,c)}\right|\leq \dfrac{C\al^2}{u'(y_c)^2}.
\eeno
Thus, we can obtain
\begin{align*}
&y_c\partial_c^2\left({u''(y_c)}\rmE(\widehat{\om}_e)+u'(y_c){\widehat{\om}_e(y_c)}\right)
=\phi_1(0,c)\al^2\cL^{\infty}+\cL^2,\\
&y_c\partial_c^2\left(\rmE(gu'')+u'(y_c)g(y_c)\right)
=\phi_1(0,c)\al^{\frac{1}{2}}\cL^{\infty}+\cL^2.
\end{align*}
from which and Lemma \ref{lem:J}, we infer that
\beno
\pa_c^2\Lambda_{3,1}(\widehat{\om}_e)=\rho^{-1}(\cL^{\infty}+\al^{-2}\cL^2),\quad \pa_c^2\Lambda_{4,1}(g)=y_c^{-2}\big(\al^{-\f32}\cL^{\infty}+\al^{-2}\cL^2\big).
\eeno
The proof is completed.
\end{proof}

Now we are in position to prove Proposition \ref{prop:Ke}. Let us recall
\begin{align*}
K_e(c,\al)=\f{\Lambda_3(\widehat{\omega}_{e})(c)\Lambda_4(g)(c)}{u'(y_c)(\rmA_2^2+\rmB_2^2)(c)}.
\end{align*}

\begin{proof}
{\bf Step 1.} $L^1$ estimate. We normalize $\|\widehat{\om}_e\|_{L^2}\leq 1$ and $\|g\|_{L^2}\leq 1$. Similar to Step 1 in Proposition \ref{prop:Ko}, we have
\beno
&&\Lambda_{1}(\widehat{\om}_e)=(1+\al\rho_0)\cL^2+\al^{\f12}\cL^{\infty},\\
&&\Lambda_{2}(g)=(1+\al\rho_0)\cL^2+\al^{\f12}\cL^{\infty},
\eeno
from which and Lemma \ref{lem:Lambda_34,1}, we infer that
\beno
&&\Lambda_{3}(\widehat{\om}_e)=(\al^{-2}+y_c^2)(1+\al\rho_0)\cL^2+\al^{\f12}y_c^2 \cL^{\infty},\\
&&\Lambda_{4}(g)=(\al^{-2}+y_c^2)(1+\al\rho_0)\cL^2+\al^{\f12}y_c^2 \cL^{\infty}.
\eeno
By Lemma \ref{lem:A2+B2-2}, we have $\f{1}{\rmA_2^2+\rmB_2^2}=\f{\al^4\cL^{\infty}}{(1+\al\rho_0)^2(1+\al y_c)^4}$. Thus,
\begin{align*}
u'(y_c)K_e(c,\al)&=\f{\Lambda_{3}(\widehat{\om}_e)\Lambda_{4}(g)}{\rmA_2^2+\rmB_2^2}\\
&=\f{\al^4\big((\al^{-2}+y_c^2)(1+\al\rho_0)\cL^2+\al^{\f12}y_c^2 \cL^{\infty}\big)\big(C(\al^{-2}+y_c^2)(1+\al\rho_0)\cL^2+\al^{\f12}y_c^2 \cL^{\infty}\big)}{(1+\al\rho_0)^2(1+\al y_c)^4}\\
&=\cL^1+\f{\al^{\f12}}{(1+\al\rho_0)}\cL^2+\f{\al }{(1+\al\rho_0)^2}\cL^{\infty}\\
&=\cL^1.
\end{align*}
which gives
\beno
\big\|K_e(c,\al)\big\|_{L_c^1}
=\big\|u'(y_c)K_e(c,\al)\big\|_{L_{y_c}^1}\leq C.
\eeno

{\bf Step 2.} {$W^{1,1}$ estimate. } We normalize $\|\widehat{\om}_e\|_{H^1}\leq 1$ and $\|g'\|_{L^2}+\al\|g\|_{L^2}\leq 1$. Then we have
\beno
&&\|g\|_{L^{\infty}}\leq C\|g\|_{L^2}^{1/2}\|g\|_{H^1}^{1/2}\leq C\al^{-\f12},\\
&&\big\|{g}/{(1-y_c)}\big\|_{L^2}\leq C\|g'\|_{L^2}\le C.
\eeno
By Lemma \ref{lem:Pi12-Linfty}, Lemma \ref{lem:Pi12-Linfty-d1} and Lemma \ref{lem:II3}, we have
\beno
&&\Lambda_{1,2}(\widehat{\om}_e)=\al\rho_0 \cL^{\infty},\quad\Lambda_{2,2}(g)=\al^{\f12}\rho_0 \cL^{\infty},\\
&&\pa_c\Lambda_{1,2}(\widehat{\omega}_e)
=y_c^{-1}\al \cL^{\infty}+C\al \cL^2,\\
&&\pa_c\Lambda_{2,2}(g)=y_c^{-1}\al^{\f12}\cL^{\infty}+\al \cL^2,
\eeno
from which, Lemma \ref{lem:Lambda-11} and Remark \ref{rem:Lambda}, we infer that
\beno
&&\Lambda_1(\widehat{\omega}_e)=\al\rho_0 \cL^{\infty}+\cL^{\infty},\\
&&\Lambda_2(g)=\al\rho_0 \cL^{\infty}+\cL^{\infty},\\
&&\pa_c\Lambda_{1}(\widehat{\omega}_e)
=y_c^{-1}\al \cL^{\infty}+\al \cL^2+y_c^{-1}\cL^2,\\
&&\pa_c\Lambda_{2}(g)=y_c^{-1}\al^{\f12} \cL^{\infty}+\al \cL^2+y_c^{-1}\cL^2.
\eeno
Then by Lemma \ref{lem:Lambda_34,1}, we obtain
\beno
&&\Lambda_3(\widehat{\omega}_e)
=\al^{-2}y_c(1+\al y_c)^2(\cL^{2}+\al \cL^{\infty}),\\
&&\Lambda_4(g)=\al^{-2}\rho_0(1+\al y_c)^2(\cL^{2}+\al^{\frac{1}{2}} \cL^{\infty}),\\
&&\partial_c\Lambda_3(\widehat{\omega}_e)
=\al^{-2}y_c^{-1}(1+\al y_c)^2 (\al \cL^{\infty}+(1+\al y_c)\cL^{2}),\\
&&\partial_c\Lambda_4(g)=\al^{-2}y_c^{-1}(1+\al y_c)^2( (1+\al \rho_0)\cL^{2}+\al^{\frac{1}{2}}\cL^{\infty}).
\eeno
By Lemma \ref{lem:A2+B2-2}, we have
\beno
\partial_c\left(\f{1}{(\rmA_2^2+\rmB_2^2){u'(y_c)}}\right)
=\f{\al^4\cL^{\infty}}{(1+\al\rho_0)^3(1+\al y_c)^3{u'(y_c)^3}}
+\f{\al^4(\cL^2\cap (\cL^{\infty}/\rho_0))}{(1+\al\rho_0)^3(1+\al y_c)^4u'(y_c)}.
\eeno
With these estimates, we deduce that
\begin{align*}
\Lambda_3(\widehat{\omega}_e)\Lambda_4(g)&=\al^{-4}\rho(1+\al y_c)^4( \cL^{2}+\al \cL^{\infty})( \cL^{2}+\al^{\frac{1}{2}} \cL^{\infty})\\
&=\al^{-4}\rho(1+\al y_c)^4( \cL^{1}+\al L^{2}+\al^{\f32}L^{\infty}),
\end{align*}
and
\begin{align*}
&\pa_c\big(\Lambda_3(\widehat{\omega}_e)\Lambda_4(g)\big)=\Lambda_3(\widehat{\omega}_e)\pa_c\Lambda_4(g)+\Lambda_4(g)\pa_c\Lambda_3(\widehat{\omega}_e)\\
&=\al^{-2}y_c(1+\al y_c)^2(\cL^{2}+\al \cL^{\infty})\al^{-2}y_c^{-1}(1+\al y_c)^2( (1+\al \rho_0)L^{2}+\al^{\frac{1}{2}}\cL^{\infty})\\
&\quad+C\al^{-2}\rho_0(1+\al y_c)^2(\cL^{2}+\al^{\frac{1}{2}}\cL^{\infty})\big(\al^{-2}y_c^{-1}(1+\al y_c)^2 (\al \cL^{\infty}+(1+\al y_c)\cL^{2})\big)\\
&=\al^{-4}(1+\al y_c)^4((1+\al \rho_0)\cL^{1}+\al(1+\al\rho_0)\cL^2+\al^{\frac{3}{2}}\cL^{\infty}).
\end{align*}
Thus, we obtain
\begin{align*}
&u'(y_c)\pa_cK_e(c,\al)\\
&=\f{\pa_c\big(\Lambda_3(\widehat{\omega}_e)\Lambda_4(g)\big)}{\rmA_2^2+\rmB_2^2}+u'(y_c)\Lambda_3(\widehat{\omega}_e)\Lambda_4(g)\partial_c\left(\f{1}{(\rmA_2^2+\rmB_2^2){u'(y_c)}}\right)\\
&=\f{( (1+\al \rho_0)\cL^{1}+\al(1+\al\rho_0)\cL^2+\al^{\frac{3}{2}}\cL^{\infty})}{(1+\al\rho_0)^2}\\
&\quad+\rho(1+\al y_c)( \cL^{1}+\al {\cL^{2}}+\al^{\f32}\cL^{\infty})\left(\f{\cL^{\infty}}{(1+\al\rho_0)^3{u'(y_c)^2}}
+\f{(\cL^2\cap (\cL^{\infty}/\rho_0))}{(1+\al\rho_0)^3(1+\al y_c)}\right)\\
&=\al^{\f12}\cL^1,
\end{align*}
which gives
\beno
\big\|\pa_cK_e(c,\al)\big\|_{L_c^1}
=\big\|u'(y_c)\pa_cK_e(c,\al)\big\|_{L_{y_c}^1}\leq C\al^{\f12}.
\eeno

{\bf Step 3.} {$W^{2,1}$ estimate. } We normalize $\|\widehat{\om}_e\|_{H^2}\leq 1$ and $\|g''-\al^2 g\|_{L^2}\leq 1$. That is, $\|g''\|_{L^2}^2+\al^2\|g'\|_{L^2}^2+\al^4\|g\|_{L^2}^2\leq 1$. Thus, we have
\beno
&&\|g\|_{L^{\infty}}\leq C\|g\|_{L^2}^{\f12}\|g\|_{H^{1}}^{\f12}\leq C\al^{-\f32},\\
&&\|g'\|_{L^{\infty}}\leq C\|g'\|_{L^2}^{\f12}\|g'\|_{H^1}^{\f12}\leq \f{C}{\sqrt{\al}}.
\eeno
Thanks to $\widehat{\om}_e'(0)=0$ and $g'(0)=g(1)=0$, we get by Hardy's inequality that
\beno
&&\big\|{\widehat{\om}_e'}/{y_c}\big\|_{L^2}\leq C\|\widehat{\om}_e''\|_{L^2}\leq C,\\
&&\big\|{g}/{(1-y_c)}\big\|_{L^{\infty}}\leq C\|g'\|_{L^{\infty}}\leq C\al^{-\f12},\\
&&\big\|{g'}/{y_c}\big\|_{L^2}\leq C\|g''\|_{L^2}\leq C.
\eeno
By Lemma \ref{lem:Pi12-Linfty}-Lemma \ref{lem:Pi12-d2} and Lemma \ref{lem:II3}, we have
\beno
&&\Lambda_{1,2}(\widehat{\om}_e)=\al\rho_0 \cL^{\infty},\quad \Lambda_{2,2}(g)=\al^{-\f12}\rho_0 \cL^{\infty},\\
&&\pa_c\Lambda_{1,2}(\widehat{\omega}_e)=y_c^{-1}\al\cL^{\infty}+\al y_cL^2\cap \al \cL^{\infty},\\
&&\pa_c\Lambda_{2,2}(g)=(y_c^{-1}\al^{-\f12}+\al^{\f12})\cL^{\infty},\\
&&\pa_c^2\Lambda_{1,2}(\widehat{\omega}_e)=y_c^{-3}(\al^{\f12}y_c+\al)\cL^{\infty}+\al y_c^{-1}\cL^2,\\
&&\pa_c^2\Lambda_{2,2}(g)=y_c^{-3}(\al^{\f12}y_c+\al^{-\f12})\cL^{\infty}+\al y_c^{-1}\cL^2,
\eeno
from which, Lemma \ref{lem:Lambda-11} and  Remark \ref{rem:Lambda},  we infer that
\beno
&&\Lambda_{1}(\widehat{\om}_e)=\al\rho_0\cL^{\infty}+\cL^{\infty},\\
&&\Lambda_{2}(g)=\al^{-\f12}\rho_0\cL^{\infty}+C\al^{-\f12}(1-y_c)\cL^{\infty},\\
&&\pa_c\Lambda_{1}(\widehat{\omega}_e)=y_c^{-1}\al \cL^{\infty}+\al y_c\cL^2\cap \al \cL^{\infty}+y_c^{-1}\cL^{\infty},\\
&&\pa_c\Lambda_{2}(g)=(y_c^{-1}\al^{-\f12}+\al^{\f12})\cL^{\infty},\\
&&\pa_c^2\Lambda_{1}(\widehat{\omega}_e)=y_c^{-3}(\al^{\f12}y_c+\al)\cL^{\infty}+\al y_c^{-1}\cL^2+y_c^{-2}\cL^2,\\
&&\pa_c^2\Lambda_{2}(g)=y_c^{-3}(\al^{\f12}y_c+\al^{-\f12})\cL^{\infty}+\al y_c^{-1}\cL^2+y_c^{-2}\cL^2.
\eeno
Then we infer from Lemma \ref{lem:Lambda_34,1} that
\begin{align*}
&\Lambda_3(\widehat{\omega}_e)=(\al^{-2}+y_c^2)(\rho \cL^{2}\cap \cL^{\infty})+y_c(y_c+\al \rho)\cL^{\infty},\\
&\Lambda_4(g)=\al^{-2}\rho(\cL^2+\al^{\frac{3}{2}}\cL^{\infty}),\\
&\partial_c\Lambda_3(\widehat{\omega}_e)
=(1+\al y_c)\cL^{\infty}+(\al^{-2}+y_c^2)\cL^{2},\\
&\partial_c\Lambda_4(g)=\al^{-2}(1+\al y_c)^2(\cL^2+\al^{\frac{1}{2}}\cL^{\infty}),\\
&\partial_c^2\Lambda_3(\widehat{\omega}_e)
=\rho^{-1}((1+\al y_c)\cL^{\infty}+(\al^{-2}+y_c^2)(1+\al \rho_0)\cL^2),\\
&\partial_c^2\Lambda_4(g)
=(1+\al^{-2} y_c^{-2}) ((1+\al \rho_0)\cL^2+\al^{\frac{1}{2}}\cL^{\infty}).
\end{align*}
By Lemma \ref{lem:A2+B2-2}, we have
\begin{align*}
&\pa_c^2\left(\f{1}{(\rmA_2^2+\rmB_2^2){u'(y_c)}}\right)\\
&=\f{\al^4\cL^{\infty}}{(1+\al\rho_0)^4(1+\al y_c)^2{u'(y_c)^5}}
+\f{\al^4(\cL^2+\al^2\cL^{\infty})}{(1+\al\rho_0)^3(1+\al y_c)^4{\rho}u'(y_c)}.
\end{align*}
With these estimates, we can deduce that
\begin{align*}
\Lambda_3(\widehat{\omega}_e)\Lambda_4(g)
&=\big((\al^{-2}+y_c^2)(\rho \cL^{2}\cap \cL^{\infty})+y_c(y_c+\al \rho)\cL^{\infty}\big)\al^{-2}\rho(\cL^2+\al^{\frac{3}{2}}\cL^{\infty})\\
&=\al^{-4}\rho(1+\al y_c)^2(\rho\cL^1\cap \cL^2)+\al^{-4}\rho(1+\al y_c)^2\al^{\f32}(\rho\cL^2\cap \cL^{\infty})\\
&\quad+y_c^2(1+\al\rho_0)\al^{-2}\rho\cL^2+\al^{\f32}\al^{-2}y_c^2(1+\al\rho_0)\rho\cL^{\infty},
\end{align*}
and
\begin{align*}
&\pa_c\big(\Lambda_3(\widehat{\omega}_e)\Lambda_4(g)\big)=\Lambda_4(g)\pa_c\Lambda_3(\widehat{\omega}_e)+\Lambda_3(\widehat{\omega}_e)\pa_c\Lambda_4(g)\\
&=\al^{-2}\rho(\cL^2+\al^{\frac{3}{2}}\cL^{\infty})\big((1+\al y_c)\cL^{\infty}+(\al^{-2}+y_c^2)\cL^{2}\big)\\
&\quad+\big((\al^{-2}+y_c^2)(\rho\cL^{2}\cap \cL^{\infty})+y_c(y_c+\al \rho)\cL^{\infty}\big)\al^{-2} (1+\al y_c)^2(\cL^2+\al^{\frac{1}{2}}\cL^{\infty})\\
&=\al^{-4}\rho (1+\al y_c)^2\cL^1+\al^{-2}\rho (1+\al y_c)^3(\cL^2+\al^{\f32}\cL^{\infty})\\
&\quad+\al^{-4} (1+\al y_c)^4\big((\rho\cL^1\cap \cL^2)+\al^{\f12}(\rho\cL^2\cap \cL^{\infty})\big)+C\al^{-2}y_c^2(1+\al\rho_0)(1+\al y_c)^2(\cL^2+\al^{\f12}\cL^{\infty})\\
&=\al^{-4}\rho (1+\al y_c)^4(\cL^1+\al^{\f12}\cL^2)+\al^{-2}\rho (1+\al y_c)^3(\cL^2+\al^{\f12}\cL^{\infty})\\
&\quad+\al^{-2}y_c^2(1+\al\rho_0)(1+\al y_c)^2(\cL^2+\al^{\f12}\cL^{\infty}),
\end{align*}
and
\begin{align*}
&\pa_c^2\big(\Lambda_3(\widehat{\omega}_e)\Lambda_4(g)\big)=\Lambda_4(g)\pa_c^2\Lambda_3(\widehat{\omega}_e)+\Lambda_3(\widehat{\omega}_e)\pa_c^2\Lambda_4(g)+2\pa_c\Lambda_4(g)\pa_c\Lambda_3(\widehat{\omega}_e)\\
&=\al^{-2}\rho(\cL^2+\al^{\frac{3}{2}}\cL^{\infty})((1+\al y_c)\cL^{\infty}+(\al^{-2}+y_c^2)(1+\al \rho_0)\cL^2)\\
&\quad+\big((\al^{-2}+y_c^2)(\rho\cL^{2}\cap\cL^{\infty})+y_c^2(1+\al \rho_0)\cL^{\infty}\big)(1+\al^{-2} y_c^{-2})((1+\al \rho_0)\cL^2+\al^{\frac{1}{2}}\cL^{\infty})\\
&\quad+\big((1+\al y_c)\cL^{\infty}+(\al^{-2}+y_c^2)\cL^{2}\big)
\al^{-2}(1+\al y_c)^2(\cL^2+\al^{\frac{1}{2}}\cL^{\infty})\\
&=\al^{-2} (1+\al y_c)^3(\cL^{2}+\al^{-2}(1+\al y_c)(1+\al \rho_0)\cL^1+\al^{\f12}\cL^{\infty}+\al^{-2}(1+\al y_c)(1+\al \rho_0)\al^{\f12}\cL^2)\\
&\quad+\al^{-2}(1+\al y_c)^2(1+\al \rho_0)\big(\al^{-2}(1+\al y_c)(1+\al \rho_0) \cL^{1}+(1+\al \rho_0)\cL^{2}+\al^{-2}(1+\al y_c) \al^{\f12}\cL^{2}+\al^{\f12}\cL^{\infty}\big)\\
&\quad+\al^{-2} (1+\al y_c)^3\big(\al^{\f12}\cL^{\infty}+\al^{\f12}\al^{-2}(1+\al y_c)\cL^{2}+\cL^{2}+\al^{-2}(1+\al y_c)\cL^{1}\big)\\
&=\al^{-2} (1+\al y_c)^3(\cL^{2}+\al^{-2}(1+\al y_c)(1+\al \rho_0)\cL^1+ \al^{\f12}\cL^{\infty}+\al^{-2}(1+\al y_c)(1+\al \rho_0)\al^{\f12}\cL^2)\\
&\quad+\al^{-2}(1+\al y_c)^2(1+\al \rho_0)^2\cL^{2},
\end{align*}
from which, we can deduce that
\begin{align*}
&u'(y_c)\pa_c^2K_e(c,\al)=\al^{\f32}\cL^1.
\end{align*}
which gives
\beno
\big\|\pa_c^2K_e(c,\al)\big\|_{L_c^1}
=\big\|u'(y_c)\pa_c^2K_e(c,\al)\big\|_{L_{y_c}^1}\leq C\al^{\f32}.
\eeno

Moreover, we have
\beno
K_e(c,\al)=C_\al\rho_0\cL^2,
\eeno
which along with the fact $K_e$ is continuous implies that $K_e(u(0),\al)=K_e(u(1),\al)=0$.
\end{proof}

\section{Appendix}

In this appendix, we present various estimates for some singular integral operators. Recall that the Hilbert transform $Hf(x)$ is defined by
\beno
H(f)(x)=p.v.\int\f{f(y)}{x-y}dy.
\eeno
In what follows, we assume that $u$ satisfies $(S)$, and let $v$ be defined by \eqref{def:v} and $\tc=v(y_c)$. Let $I_v=(-v(1), v(1)), I_v^+=(0,v(1)), I_v^-=(-v(1),0)$.

For a function $f$ defined in $I_v$, $H(f)(\tc)$ denotes the Hilbert transform of $\widetilde{f}(\tc)=f\chi_{I_v}(\tc)$. We denote by
$\|\cdot\|_{L^p}$ the norm of $L^p(I_v)$, and $\|\cdot\|_{W^{k,p}}$ or $\|\cdot\|_{H^k}$ the Sobolev norm of $W^{k,p}(I_v)$
or $H^k(I_v)$.

\begin{lemma}\label{lem:Hilbert}
Let $g\in W^{1,p}(I_v)$ for $p\in (1,+\infty)$.
If $g$ is even, then we have
\begin{align*}
\f{d}{d\tc}H(g)(\tc)&=H(g')(\tc)+\frac{2g(v(1))v(1)}{v(1)^2-\tc^2}.
\end{align*}
If $g$ is odd, then we have
\begin{align*}
&\f{d}{d\tc}H(g)(\tc)=H(g')(\tc)+\frac{2g(v(1))\tc}{v(1)^2-\tc^2},\\
&\tc H(g')(\tc)=H(zg')(\tc)+2g(v(1)),\\
&\tc H(g)(\tc)=H(zg)(\tc).
\end{align*}
\end{lemma}
\begin{proof}
It is easy to check that
\begin{align}
\f{d}{d\tc}H(g)(\tc)=H(g')(\tc)
-\f{g(v(1))}{\tc-v(1)}+\f{g(-v(1))}{\tc+v(1)},\label{eq:H-d}
\end{align}
and for $g(z)=-g(-z)$,
\begin{align*}
\tc H(g)(\tc)-H(zg)(\tc)=\int_{-v(1)}^{v(1)}\f{\tc g(z)-zg(z)}{\tc-z}dz=0,
\end{align*}
which give the lemma.
\end{proof}

We introduce the average operator $A^{(z)}(g)(\tc)$ defined by
\beno
A^{(z)}(g)(\tc)=\f 1{z-\tc}\int_{\tc}^zg(z')dz'.
\eeno
Notice that  for any $k\in \N$,
\begin{align*}
\f{d^k}{d\tc^k}A^{(z)}(g)(\tc)=\f{1}{z-\tc}\int_{\tc}^{z}g^{(k)}(z')\Big(\f{z'-z}{\tc-z}\Big)^k\,dz'.
\end{align*}
We infer that for any $p\in (1,+\infty]$
\beq\label{eq:aver}
\|A^{(z)}(g)\|_{W^{k,p}}\leq C\|g\|_{W^{k,p}}.
\eeq

\begin{lemma}\label{lem:Z}
We define
\beno
Z(g)(\tc)=g'(\tc)-\f{g(\tc)}{\tc}.
\eeno
For any even function $g\in W^{k,p}\big(I_v^\pm\big)\cap C(I_v)$ for $p\in (1,+\infty), k=2,3$ and $g(0)=0$, we have $Z(g)\in W^{k-1,p}(I_v)$
and
\beno
\|Z(g)\|_{W^{k-1,p}}\leq C\|g\|_{W^{k,p}(I_v^+)}.
\eeno
\end{lemma}

\begin{proof}
Since $g$ is even and $g(0)=0$, we have
\ben
\|Z(g)\|_{L^p}\le C\|g\|_{W^{1,p}(I_v^+)}.\label{eq:Z-est1}
\een
Take $f\in C_0^\infty(I_v)$ with $\|f\|_{L^{p'}}\le 1$.
Thanks to $Z(g)(z)=g'(z)-A^{(0)}(g')(z)$ and $\lim_{z\to 0\pm}Z(g)(z)=0$, we get by integration by parts that
\begin{align}
\left|\int_{I_v}Z(g)(z)f'(z)\,dz\right|=&\bigg|
\int_{I_v^-}\Big(g'(z)-A^{(0)}(g')(z)\Big)f'(z)\,dz+\int_{I_v^+}\Big(g'(z)-A^{(0)}(g')(z)\Big)f'(z)\,dz\bigg|\nonumber\\
&\leq \bigg(
\left|\int_{I_v^-}g''(z)f(z)dz\right|
+\left|\int_{I_v^+}g''(z)f(z)dz\right|\nonumber\\
&\qquad+\left|\int_{I_v^-}\big(A^{(0)}(g')\big)'(z)f(z)dz\right|
+\left|\int_{I_v^+}\big(A^{(0)}(g')\big)'(z)f(z)dz\right|\bigg)\nonumber\\
&\leq C\Big(\|g''\|_{L^p(I_v^+)}+\big\|\big(A^{(0)}(g')\big)'\big\|_{L^p(I_v^+)}\Big)
\leq C\|g\|_{W^{2,p}(I_v^+)}.\label{eq:Z-est2}
\end{align}

Thanks to
\beno
\lim_{z\to 0+}A^{(0)}(g')'(z)=\lim_{z\to 0-}A^{(0)}(g')'(z)=\f{g''(0)}{2},
\eeno
we get by integration by parts that
\begin{align}
\left|\int_{I_v}Z(g)(z)f''(z)\,dz\right|\leq&\bigg|\int_{I_v^-}\big(g''(z)-A^{0}(g')'(z)\big)f'(z)dz+\int_{I_v^+}\big(g''(z)-A^{0}(g')'(z)\big)f'(z)dz\bigg|\nonumber\\
&\leq \bigg|\f12{g''(0)}f(0)-\int_{I_v^-}\big(g'''(z)-A^{0}(g')''(z)\big)f(z)dz\nonumber\\
&\quad-\f12g''(0)f(0)-\int_{I_v^+}\big(g'''(z)-A^{0}(g')''(z)\big)f(z)dz\bigg|\nonumber\\
\leq& C\Big(\|g'''\|_{L^p(I_v^+)}+\big\|\big(A^{(0)}(g')\big)''\big\|_{L^p(I_v^+)}\Big)
\leq C\|g\|_{W^{3,p}(I_v^+)}.\label{eq:Z-est3}
\end{align}

Then the lemma follows from \eqref{eq:Z-est1}-\eqref{eq:Z-est3}.
\end{proof}

\begin{proposition}\label{prop:SIO}\
Let $g\in C(I_v)$ be an even function with $g(0)=0$.
Then it holds that

\begin{itemize}

\item[1.] if  $\varphi\in W^{1,p}(I_v)\cap W^{2,p}(I_v^\pm)$ for $p\in (1,+\infty)$, then we have
\beno
\left\|\tc\partial_c\left(\rho(c) \partial_c\Big(\frac{1}{2\widetilde{c}}H(g)(\widetilde{c})\Big)\right)\right\|_{L^p}
\leq C\|g\|_{W^{2,p}(I_v^+)};
\eeno

\item[2.] if $g\in H^3\big(I_v^\pm\big)$,  then we have
\beno
\left\|\rho(c)\partial_c^2\left(\rho(c) \partial_c\Big(\frac{1}{2\widetilde{c}}H(g)(\widetilde{c})\Big)\right)\right\|_{L^2}\leq C\|g\|_{H^3(I_v^+)},
\eeno
and
\begin{align*}
&\left|\tc\partial_c\left(\rho(c) \partial_c\Big(\frac{1}{2\widetilde{c}}H(g)(\widetilde{c})\Big)\right)\right|
\leq C\|g\|_{H^2(I_v^+)}\big|\ln (u(1)-c)\big|
+C\|g\|_{H^3(I_v^+)},\\
&\left|\rho(c)\partial_c\Big(\frac{1}{2\widetilde{c}}H(g)(\widetilde{c})\Big)\right|
\leq C|\tc|\|g\|_{H^3(I_v^+)}.
\end{align*}

\item[3.] if $g\in C^2(I_v)\cap C^3\big(I_v^\pm\big)$ and $\tc g(\tc)\in C^3(I_v)$, then for any $p\in (1,+\infty)$,  we have
\begin{align*}
&\left\|\tc\rho(c)\partial_c^2\left(\rho(c) \partial_c\Big(\frac{1}{2\widetilde{c}}H(g)(\widetilde{c})\Big)\right)\right\|_{L^p}
\leq C\|g\|_{C^2(I_v)}+C\|zg\|_{C^3(I_v)},\\
&\left|\tc\partial_c\left(\rho(c)\partial_c\Big(\frac{1}{2\widetilde{c}}H(g)(\widetilde{c})\Big)\right)\right|
\leq C\left(\|g\|_{C^2(I_v)}+\|zg\|_{C^3(I_v)}\right)\left|\ln (u(1)-c)\right|,\\
&\left|\rho(c) \partial_c\Big(\frac{1}{2\widetilde{c}}H(g)(\widetilde{c})\Big)\right|\leq C|\tc|\|g\|_{C^2(I_v)}+C|\tc|\|zg\|_{C^3(I_v)}.
\end{align*}

\end{itemize}
\end{proposition}

\begin{proof}
Let us prove the first case. If $g$ is even, it follows from Lemma \ref{lem:Hilbert} that
\begin{align*}
&(v(1)^2-\tc^2)(\tc (Hg)'(\tc)-Hg(\tc))\\
&=(v(1)^2-\tc^2)(\tc H(g')(\tc)-Hg(\tc))+2\tc g(v(1))v(1)\\
&=(v(1)^2-\tc^2) H(zg'-g)(\tc)+2\tc g(v(1))v(1),
\end{align*}
which gives
\begin{align*}
\rho(c) \partial_c\Big(\frac{1}{2\widetilde{c}}H(g)(\widetilde{c})\Big)
&=\f{\tc}{2}(v(1)^2-\tc^2)\partial_{\tc}\Big(\frac{1}{2\widetilde{c}}H(g)(\widetilde{c})\Big)\\
&=\f{v(1)^2-\tc^2}{4\tc}\big(\tc (Hg)'(\tc)-Hg(\tc)\big)\\
&=\f{u(1)-c}{4} H\left(Z(g)\right)(\tc)+ \f{g(v(1))v(1)}{2}.
\end{align*}
As $Z(g)$ is odd, we get by Lemma \ref{lem:Hilbert} again that
\begin{align}
\partial_c\left(\rho(c) \partial_c\Big(\frac{1}{2\widetilde{c}}H(g)(\widetilde{c})\Big)\right)
&=-\f{1}{4} H(Z(g))(\tc)+\f{v(1)^2-\tc^2}{8\tc} (H Z(g))'(\tc)\nonumber\\
&=-\f{1}{4} H(Z(g))(\tc)+\f{v(1)^2-\tc^2}{8\tc} H( Z(g)')(\tc)+\f{1}{4}Z(g)(v(1)),
\label{eq:Hg-d1}
\end{align}
from which and Lemma \ref{lem:Z}, we infer that
\begin{align*}
\left\|\tc\partial_c\left(\rho(c)\partial_c\Big(\frac{1}{2\widetilde{c}}H(g)(\widetilde{c})\Big)\right)\right\|_{L^p}&\leq C\|H(Z(g))\|_{L^p}+C\|H( Z(g)')\|_{L^p}
+C\|Z(g)\|_{L^{\infty}}\\
&\leq C\|Z(g)\|_{W^{1,p}}\le C\|g\|_{W^{2,p}(I_v^+)}.
\end{align*}

Now we consider the second case. As  $Z(g)$ is odd and $Z(g)'$ is even, we deduce from  Lemma \ref{lem:Hilbert} and \eqref{eq:Hg-d1} that
\begin{align*}
&\partial_c^2\left(\rho(c) \partial_c\Big(\frac{1}{2\widetilde{c}}H(g)(\widetilde{c})\Big)\right)\\
&=-\f{1}{8\tc}\left( H(Z(g)')(\tc)
+\frac{2Z(g)(v(1))\tc}{v(1)^2-\tc^2}\right)
-\f{v(1)^2+\tc^2}{16\tc^3} H( Z(g)')(\tc)\\
&\quad+\f{v(1)^2-\tc^2}{16\tc^2}\left( H( Z(g)'')(\tc)
+\frac{2Z(g)'(v(1))v(1)}{v(1)^2-\tc^2}\right)\\
&=-\f{v(1)^2+3\tc^2}{16\tc^3} H( Z(g)')(\tc)
+\f{v(1)^2-\tc^2}{16\tc^2} H( Z(g)'')(\tc)
-\frac{Z(g)(v(1))}{4(v(1)^2-\tc^2)}
+\frac{Z(g)'(v(1))v(1)}{8\tc^2}.
\end{align*}
Due to $Z(g)'(0)=\f{g''(0)}{2}$, we have
\begin{align*}
H( Z(g)')(\tc)&=p.v.\int_{-v(1)}^{v(1)}\f{Z(g)'(z)-\f{g''(0)}{2}}{\tc-z}dz+\f{g''(0)}{2}\ln \f{v(1)-\tc}{v(1)+\tc}\\
&=\tc H\Big(A^{(0)}\big(Z(g)''\big)\Big)(\tc)+\f{g''(0)}{2}\ln \f{v(1)-\tc}{v(1)+\tc}.
\end{align*}
Thus, we get by Lemma \ref{lem:Z} that
\begin{align*}
&\left\|\rho(c)\partial_c^2\left(\rho(c) \partial_c\Big(\frac{1}{2\widetilde{c}}H(g)(\widetilde{c})\Big)\right)\right\|_{L^2}\\
&\leq C\left\|H\Big(A^{(0)}\big(Z(g)''\big)\Big)\right\|_{L^2}+C|g''(0)|+C\|H( Z(g)'')\|_{L^2}
+C\|Z(g)\|_{W^{1,\infty}}\\
&\leq C\|g\|_{H^3(I_v^+)}.
\end{align*}

For the pointwise estimate, we need to use the following estimate
\begin{align}
|H(g)(\tc)|
&=\left|p.v.\int_{-v(1)}^{v(1)}\f{g(z)}{\tc-z}dz\right|\nonumber\\
&\leq \left|g(z)\ln|\tc-z|\big|_{z=-v(1)}^{v(1)}\right|+\left|\int_{-v(1)}^{v(1)}g'(z)\ln|\tc-z|dz\right|\nonumber\\
&\leq C\|g\|_{L^{\infty}}{(|\ln (u(1)-c)|+1)}+C\|g'\|_{L^2}\nonumber\\
&\leq C\|g\|_{H^1}{(|\ln (u(1)-c)|+1)}.\label{eq:H-p-est}
\end{align}
Then we infer from \eqref{eq:Hg-d1} and \eqref{eq:H-p-est} that
\begin{align*}
\left|\tc\partial_c\left(\rho(c) \partial_c\Big(\frac{1}{2\widetilde{c}}H(g)(\widetilde{c})\Big)\right)\right|&\leq C\|Z(g)\|_{H^1}{(|\ln (u(1)-c)|+1)}
+C\|Z(g)'\|_{H^1}
+C\|Z(g)\|_{L^{\infty}}\\
&\leq C\|g\|_{H^2(I_v^+)}{(|\ln (u(1)-c)|+1)}
+C\|g\|_{H^3(I_v^+)},
\end{align*}
which along with the fact that $\rho(c) \partial_c\Big(\frac{1}{2\widetilde{c}}H(g)(\widetilde{c})\Big)\bigg|_{\tc=0}=0$(see Remark \ref{rem:Hg}) gives
\beno
\left|\rho(c) \partial_c\Big(\frac{1}{2\widetilde{c}}H(g)(\widetilde{c})\Big)\right|\leq C|\tc|\|g\|_{H^3(I_v^+
)}.
\eeno

Next we consider the third case.  By the above proof, we have
\begin{align*}
&\tc\rho(c)\partial_c^2\left(\rho(c) \partial_c\Big(\frac{1}{2\widetilde{c}}H(g)(\widetilde{c})\Big)\right)\\
&=-\f{(v(1)^2+3\tc^2)(v(1)^2-\tc^2)}{16} H( Z(g)')(\tc)
+\f{\tc(v(1)^2-\tc^2)^2}{16} H(Z(g)'')(\tc)\\
&\quad-\frac{\tc^3Z(g)(v(1))}{4}
+\frac{\tc(v(1)^2-\tc^2)Z(g)'(v(1))v(1)}{8}.
\end{align*}
As $Z(g)''$ is odd, we have $\tc H\big( Z(g)''\big)(\tc)=H\big(zZ(g)''\big)(\tc)$, thus
\begin{align*}
\|\tc H\big( Z(g)''\big)\|_{L^p}
&\leq C\|\tc Z(g)''\|_{L^p}\\
&\leq C\left\|\tc Z(g)\right\|_{W^{2,p}}+C\left\|Z(g)\right\|_{W^{1,p}}\\
&\leq C\|\tc g\|_{C^3(I_v)}+C\|g\|_{C^2(I_v)}.
\end{align*}
As $Z(g)'$ is even, we have
\begin{align*}
&(v(1)^2-\tc^2)H(Z(g)')(\tc)\\
&=(v(1)^2-\tc^2)Z(g)'(v(1))\ln \left|\f{v(1)-\tc}{v(1)+\tc}\right|
-(v(1)^2-\tc^2)\int_{-v(1)}^{v(1)}Z(g)''(z)\ln |\tc-z|dz\\
&=(v(1)^2-\tc^2)Z(g)'(v(1))\ln \left|\f{v(1)-\tc}{v(1)+\tc}\right|
-\f12(v(1)^2-\tc^2)\int_{-v(1)}^{v(1)}zZ(g)''(z)\f{1}{z}\ln \left|\f{\tc-z}{\tc+z}\right|dz,
\end{align*}
which implies that for any $p\in (1,\infty]$,
\begin{align}
\left\|(v(1)^2-\tc^2)H(Z(g)')\right\|_{L^p}
\leq C\|\tc g\|_{C^3(I_v)}+C\|g\|_{C^2(I_v)}.\label{eq:Hg-est-p}
\end{align}
Thus, we obtain
\beno
\left\|\tc\rho(c)\partial_c^2\left(\rho(c) \partial_c\Big(\frac{1}{2\widetilde{c}}H(g)(\widetilde{c})\Big)\right)\right\|_{L^p}
\leq C\|g\|_{C^2(I_v)}+C\|\tc g\|_{C^3(I_v)}.
\eeno

By \eqref{eq:Hg-d1}, we have
\begin{align*}
\tc\partial_c\left(\rho(c) \partial_c\Big(\frac{1}{2\widetilde{c}}H(g)(\widetilde{c})\Big)\right)
=-\f{\tc}{4} H(Z(g))(\tc)+\f{v(1)^2-\tc^2}{8} H( Z(g)')(\tc)+\f{\tc}{4}Z(g)(v(1)),
\end{align*}
which along with \eqref{eq:H-p-est} and \eqref{eq:Hg-est-p} gives
\begin{align*}
&\left|\tc\partial_c\left(\rho(c) \partial_c\Big(\frac{1}{2\widetilde{c}}H(g)(\widetilde{c})\Big)\right)\right|\\
&\leq C\|H(Z(g))\|_{L^{\infty}}+C\|\tc g\|_{C^3(I_v)}+C\|g\|_{C^2(I_v)}\\
&\leq C\|g\|_{C^2(I_v)}\big|\ln (u(1)-c)\big|
+C\|\tc g\|_{C^3(I_v)}+C\|g\|_{C^2(I_v)}.
\end{align*}
This along with the fact that $\rho(c) \partial_c\Big(\frac{1}{2\widetilde{c}}H(g)(\widetilde{c})\Big)\bigg|_{\tc=0}=0$ gives
\beno
\left|\rho(c) \partial_c\Big(\frac{1}{2\widetilde{c}}H(g)(\widetilde{c})\Big)\right|\leq C|\tc|\|g\|_{C^2(I_v)}+C|\tc|\|\tc g\|_{C^3(I_v)}.
\eeno

This completes the proof of the proposition.
\end{proof}

\begin{remark}\label{rem:Hg}
Using the formula
\beno
p.v.\int_{-v(1)}^{v(1)}\f{g(z)}{\tc-y}dy=p.v.\int_{-v(1)}^{v(1)}\f{g(z)-g(0)}{\tc-z}dz+\f{g(0)}{2}\ln \left|\f{\tc-v(1)}{\tc+v(1)}\right|,
\eeno
and  the fact that $\xi(\tc)=\f{1}{\tc}(g(\tc)-g(0))$ is odd, we get by Lemma \ref{lem:Hilbert}  that
\begin{align*}
p.v.\int_{-v(1)}^{v(1)}\f{g(z)}{\tc-z}dz=p.v. \tc\int_{-v(1)}^{v(1)}\f{\xi(z)}{\tc-z}dz+\f{g(0)}{2}\ln \left|\f{\tc-v(1)}{\tc+v(1)}\right|.
\end{align*}
Therefore,
\begin{align*}
&\f{1}{2}\tc (u(1)-c) \partial_{\tc}\Big(\frac{1}{2\widetilde{c}}H(g)(\widetilde{c})\Big)\\
&=\f{1}{2}\tc (u(1)-c) \partial_{\tc}\Big(\f12\int_{-v(1)}^{v(1)}\f{\xi(z)}{\tc-z}dz+\f{g(0)}{2\tc}\ln \left|\f{\tc-v(1)}{\tc+v(1)}\right|\Big)\\
&=\f{1}{2}\tc (u(1)-c) \Big(\f12\int_{-v(1)}^{v(1)}\f{\xi'(z)}{\tc-z}dz+\f{\xi(v(1))\tc}{u(1)-c}\Big)+\partial_{\tc}\Big(\f{g(0)}{2\tc}\ln \left|\f{\tc-v(1)}{\tc+v(1)}\right|\Big),
\end{align*}
which implies that $\rho(c) \partial_c\Big(\frac{1}{2\widetilde{c}}H(g)(\widetilde{c})\Big)\bigg|_{\tc=0}=0$.

\end{remark}

\begin{lemma}\label{lem:boundaryvanish}
Assume that $\varphi_k\in W^{1,p}(I_v), p\in (1,+\infty)$ for $k=1,2,3,4$ satisfies
\ben
\varphi_1(\tc)\varphi_4(\tc)+\varphi_3(\tc)\varphi_2(\tc)\big|_{\tc=\pm v(1)}=0.\label{ass:bc}
\een
Let $G(\tc)=\varphi_4(\tc)H(\varphi_1)(\tc)+\varphi_3(\tc)H(\varphi_2)(\tc)$. Then it holds that
\beno
&&\left\|\pa_{\tc}\left((u(1)-c)G\right)\right\|_{L^{p}}
\leq C\|\varphi_1\|_{W^{1,p}}\|\varphi_4\|_{W^{1,p}}+C\|\varphi_2\|_{W^{1,p}}\|\varphi_3\|_{W^{1,p}},\\
&&\left\|\pa_{\tc}^2\left((u(1)-c)G\right)\right\|_{L^{2}}
\leq C\|\varphi_1\|_{H^2}\|\varphi_4\|_{H^2}+\|\varphi_2\|_{H^2}\|\varphi_3\|_{H^2}.
\eeno
\end{lemma}

\begin{proof}
By \eqref{eq:H-d} and \eqref{ass:bc}, it is easy to deduce that
\begin{align*}
&\f{d}{d\tc}\big(\varphi_4(\tc)H(\varphi_1)(\tc)+\varphi_3(\tc)H(\varphi_2)(\tc)\big)\\
&=\varphi_4'(\tc)H(\varphi_1)(\tc)+\varphi_3'(\tc)H(\varphi_2)(\tc)
+\varphi_4(\tc)\f{d}{d\tc}H(\varphi_1)(\tc)+\varphi_3(\tc)\f{d}{d\tc}H(\varphi_2)(\tc)\\
&=\varphi_4'(\tc)H(\varphi_1)(\tc)+\varphi_3'(\tc)H(\varphi_2)(\tc)
+\varphi_4(\tc)\Big(H(\varphi_1')(\tc)-\frac{\varphi_1(v(1))}{\tc-v(1)}+\frac{\varphi_1(-v(1))}{\tc+v(1)}\Big)\\
&\quad+\varphi_3(\tc)\Big(H(\varphi_2')(\tc)-\frac{\varphi_2(v(1))}{\tc-v(1)}+\frac{\varphi_2(-v(1))}{\tc+v(1)}\Big)\\
&=\varphi_4'(\tc)H(\varphi_1)(\tc)+\varphi_3'(\tc)H(\varphi_2)(\tc)+\varphi_4(\tc)H(\varphi_1')(\tc)
+\varphi_3(\tc)H(\varphi_2')(\tc)\\
&\quad-\varphi_1(v(1))A^{(v(1))}(\varphi_4')(\tc)+\varphi_1(-v(1))A^{(-v(1))}(\varphi_4')(\tc)\\
&\quad-\varphi_2(v(1))A^{(v(1))}(\varphi_3')(\tc)+\varphi_2(-v(1))A^{(-v(1))}(\varphi_3')(\tc),
\end{align*}
which along with \eqref{eq:aver} and \eqref{eq:H-p-est} gives the first inequality of the lemma. The proof of the second inequality is similar.
\end{proof}

\section*{Acknowledgement}
Z. Zhang is partially supported by NSF of China under Grant 11425103.

\end{document}